\RequirePackage[leqno]{amsmath}

\documentclass[pdflatex,sn-mathphys-num]{sn-jnl}

\usepackage[title]{appendix}

\usepackage{anyfontsize,bold-extra}
\usepackage{scalerel}

\usepackage{graphicx}
\usepackage{amsthm,amssymb}
\usepackage{placeins}
\usepackage{xcolor}
\usepackage{multirow}
\usepackage[utf8]{inputenc}

\usepackage[linesnumbered, boxed, english, onelanguage, vlined]{algorithm2e}
\usepackage{mathtools}
\usepackage{mathrsfs}

\usepackage[linesnumbered, boxed, english, onelanguage, vlined]{algorithm2e}
\usepackage{mathtools}
\usepackage{mathrsfs}
\usepackage{subfigure}

\allowdisplaybreaks

\vbadness=\maxdimen
\newtheorem{theorem}{Theorem}
\newtheorem{corollary}[theorem]{Corollary}
\newtheorem{lemma}[theorem]{Lemma}
\newtheorem{proposition}[theorem]{Proposition}
\theoremstyle{remark}
\newtheorem{remark}[theorem]{Remark}

\DeclareMathOperator{\ldet}{ldet}
\DeclareMathOperator{\Trace}{Tr}
\DeclareMathOperator{\rank}{rank}
\DeclareMathOperator*{\argmin}{argmin}

\DeclareMathOperator{\diag}{\textstyle{diag}}
\DeclareMathOperator{\Diag}{\textstyle{Diag}}

\DeclareMathOperator{\vvec}{\textstyle{vec}}

\makeatletter
\renewcommand*{\top}{%
  {\mathpalette\@transpose{}}%
}
\newcommand*{\@transpose}[2]{%
  \raisebox{\depth}{$\m@th#1\scriptscriptstyle\mathsf{T}$}%
}
\makeatother

\begin{document}

\title[ADMM for 0/1 D-Opt and MESP relaxations]{ADMM for 0/1 D-optimality and Maximum-Entropy Sampling Relaxations}

\author[1]{\fnm{Gabriel} \sur{Ponte}}\email{gabponte@umich.edu}

\author*[2]{\fnm{Marcia} \sur{Fampa}}\email{fampa@cos.ufrj.br}

\author[1]{\fnm{Jon} \sur{Lee}}\email{jonxlee@umich.edu}

\author[3]{\fnm{Luze} \sur{Xu}}\email{xuluze@ust.hk}

\affil[1]{\orgdiv{University of Michigan}
}

\affil[2]{\orgdiv{Universidade Federal do Rio de Janeiro}
}

\affil[3]{\orgdiv{Hong Kong University of Science and Technology}
}

\abstract{The  0/1 D-optimality problem and the Maximum-Entropy Sampling problem are two well-known NP-hard discrete maximization problems in experimental design.
Algorithms for exact optimization (of moderate-sized instances)
are based on branch-and-bound. The best upper-bounding methods are based on convex relaxation. We present ADMM (Alternating Direction Method of Multipliers) algorithms for 
solving these relaxations and experimentally demonstrate their practical value.}

\date{\today}


\maketitle


\section{Introduction}

Some challenging families of discrete nonlinear-optimization problems come from the area of experimental design. One important 
problem is the Gaussian case of the
0/1 D-optimality problem (D-Opt).
Briefly, the problem aims to select a subset of 
$s$ design points, from a set of
$n$ given design points in $\mathbb{R}^m$, with the goal of 
minimizing the ``generalized variance'' of the least-squares 
parameter estimates; see, for example, \cite*{PonteFampaLeeMPB}
and the references therein. Another problem is the Gaussian case of the maximum-entropy sampling problem (MESP). Here we have an input covariance matrix of order $n$, and we wish to select a principal submatrix of order $s$, so as to maximize the 
``differential entropy'' (see, for example, \cite*{FL2022}).

Generally, the workhorse algorithm for ``convex MINLO'' (that is, mixed-integer nonlinear optimization, where the continuous relaxation are convex optimization problems) is B\&B (branch-and-bound); see, for example, \cite*{BONMIN} and \cite*{muriqui}.
Indeed, such approaches have been developed specifically for 0/1 D-Opt (see \cite*{PonteFampaLeeMPB}) and MESP
(see \cite*{FL2022}). In general, for convex MINLO, many convex relaxations must be solved very quickly.  For 0/1 D-Opt and MESP in particular, the computational effort to solve 
individual convex relaxations can be quite substantial (in contrast to, for example, 0/1 linear optimization, where convex relaxations are linear programs that can additionally be warm- or hot-started very effectively). 

In what follows, we present fast ADMM (Alternating Direction Method of Multipliers) algorithms to
solve some convex relaxations (from the literature) for 0/1 D-Opt and MESP,
with an eye toward their future use in B\&B.

We assume some familiarity with particular aspects of convex optimization, 
in particular, with the well-known ADMM 
and its development from the augmented Lagrangian function, 
for problems of the form $\min_{x,z} \{f(x)+g(z) ~|~ Ax + Bz = c\}$;
see, for example, \cite*{boyd2011distributed}.
When we can implement the iteration updates quickly,
ADMM (and other first-order methods), are a method of choice
for approximately solving large-scale convex-optimization problems.
We assume a bit of familiarity with B\&B for convex MINLO, but most of
what we present can be appreciated without any knowledge of that topic.
A particularly nice feature of ADMM versus \emph{some} primal methods, in the B\&B context, is that 
 warm-starting a child solve from
a parent solution is trivial. Of course, it remains to be seen if 
warm starting like this is effective within B\&B.

\medskip
\noindent{\bf Brief literature review.}  
D-optimality, whose criterion is
maximizing the (logarithm of the) determinant of an appropriate positive-definite matrix,
is a very well-studied topic in the experimental design literature. There are many variations, and we concentrate on the 0/1 version of the problem, which we carefully state in \S\ref{sec:d-opt}. A recent reference on the state-of-the-art for B\&B approaches is \cite*{PonteFampaLeeMPB}, with many references therein to  
background and previous work. A key upper bound
based on convex relaxation is the ``natural bound'',
and we propose herein an ADMM algorithm for its fast calculation. Related to this is \cite*{nagata2021data}, which proposes an ADMM for ``A-optimal design'' (which seeks to maximize a trace). More similar is  \cite*{Scheinberg2010},
which 
gives an ADMM algorithm for $\max\{\log \det (X) - \Trace(SX) - \tau \|X\|_1\}$, a convex relaxation for ``sparse inverse covariance selection''.

MESP is a closely related problem  
 in the experimental-design literature, 
  which we carefully state in \S\ref{sec:MESP}. 
  A recent reference on the state-of-the-art for B\&B approaches is \cite*{FL2022}, with many references therein to  
background and previous work. Key upper bounds based on convex relaxation are the ``linx bound'' (see \cite*{Kurt_linx}), the ``factorization bound'' (see 
\cite*{Nikolov,Weijun,FL2022,ChenFampaLee_Fact}),
and the ``BQP bound'' (see \cite*{Anstreicher_BQP_entropy}), and in the sequel, we propose new ADMM algorithms for their calculation. 
There is also an important ``factorization bound'' for 
0/1 D-optimality instances, but we can see it as applying the 
MESP ``factorization bound'' to an appropriately-constructed 
instance of MESP (see \cite*{PonteFampaLeeMPB} and \cite*{MESP2DOPT}). 

\smallskip
\noindent{\bf Organization and contributions.}
In \S\ref{sec:d-opt}, we present a new ADMM algorithm for the natural bound for D-Opt. In \S\ref{sec:MESP}, we present a new  ADMM algorithm for the factorization bound for MESP, which requires significant new theoretical results. We also present new ADMM algorithms for the linx and BQP bounds for MESP.
In \S\ref{sec:numexp}, we present results of numerical experiments, demonstrating the benefits of our approach. 
Specifically, we will see that our ADMM algorithm for the natural bound for D-Opt is significantly better for large instances than applying commercial (and other) solvers. Additionally, 
we will see that while our ADMM algorithm for the linx bound for MESP does not perform well compared to commercial solvers,
our ADMM algorithm for the factorization bound for MESP does perform quite well. 
Another highlight is that with our ADMM approach, we 
could calculate the BQP bound for MESP for much larger
instances than was previously possible.
In \S\ref{sec:concl}, we make some concluding remarks.
In Appendices A and B, we have some supporting material. 

All of our ADMM algorithms are for convex minimization problems.
Because our problems satisfy appropriate technical conditions (see, for example, \cite*[Section 3.2]{boyd2011distributed}), our (2-block) ADMMs are guaranteed
to globally converge (using any positive penalty parameter).\footnote{However, we cannot directly apply standard results to guarantee \emph{fast} convergence for our ADMMs, 
because $\log\,\det(\cdot)$ is not strongly concave
on all of $\mathbb{S}^n_{++}$\,; see \cite*{deng_yin}, and the references therein.}  Specifically, from \cite[Section 3.2.1]{boyd2011distributed},
for all of the ADMMs that we present:
(i) the primal iterates converge to a feasible (but not necessarily optimal) point, 
(ii) the objective iterates converge to the optimal value, and
(iii) The dual iterates converge to a dual optimal point.
We note that this behavior is sufficient for our B\&B use case,
as we provide additional techniques to obtain true dual-feasible solutions.
In any case, our numerical experiments
demonstrate the practical effectiveness of our ADMM algorithms.  We do note that for our experiments, we 
make modifications to ADMM to gain practical speed, and then
theoretical convergence proofs no longer apply.

\smallskip
\noindent{\bf Notation.}
Throughout, we denote any all-zero square matrix simply by $0$, while we denote any all-zero (column) vector by $\mathbf{0}$.
We denote any all-one vector
by $\mathbf{e}$, any \hbox{$i$-th} standard unit vector by $\mathbf{e}_i$\,, 
any all-one matrix by $J$,
 and the identity matrix of order $n$ by $I_n$\,.
 We let $\mathbb{S}^n$  (resp., $\mathbb{S}^n_+$~, $\mathbb{S}^n_{++}$)
 denote the set of symmetric (resp., positive-semidefinite, positive-definite) matrices of order $n$. 
 We let $\Diag(x)$ denote the $n\times n$ diagonal matrix with diagonal elements given by the components of $x\in \mathbb{R}^n$, and we let $\diag(X)$ denote the $n$-vector with elements given by the diagonal elements of $X\in\mathbb{R}^{n\times n}$.
 When $X$ is symmetric, we let $\lambda(X)$ denote its non-increasing
 list of real eigenvalues.
We let
$\ldet$ denote the natural logarithm of the determinant,
and we freely use the facts that $\ldet(\cdot)$ is a (strictly) concave function on $\mathbb{S}^n_{++}$ and that $\nabla \ldet(M)=M^{-\top}$, when
$M\in\mathbb{R}^{n\times n}$ has positive determinant
(see, for example, Theorem A.4.12 and Lemma A.4.13 in \cite*{FL2022}).
We let $\Trace$ denote the trace. 
 We denote Frobenius norm by $\|\cdot\|_F$ and 2-norm by $\|\cdot\|_2$\,.
For a matrix $M$, we denote row $i$ by $M_{i\cdot}$ and
column $j$ by $M_{\cdot j}$\,.
For compatible $M_1$ and $M_2$\,, 
$M_1\bullet M_2:=\Trace(M_1^\top M_2)$ is the matrix dot-product, and $M_1\circ M_2$ is the Hadamard (i.e., element-wise) product. 
For any symmetric matrix $M$, $\vvec_{\scaleto{\Delta}{4.0pt}}(M)$ is defined to be the vectorization of the lower-triangular matrix of $M$ with off-diagonal elements multiplied by $\Delta$. This notation is helpful
to transfer between dot products and norms on symmetric matrices and related quantities on vectors.
In particular, for a symmetric matrix $M$,
we have
$
\|M\|^2_{\rm F}=\sum_i \sum_j M_{ij}^2
=\sum_i M_{ii}^2 + 2\sum_{i<j} M_{ij}^2
= \sum_i M_{ii}^2 + \sum_{i<j} \left( \sqrt{2} M_{ij} \right)^2
= \|\vvec_{\scaleto{\sqrt{2}}{4.0pt}} (M)\|_2^2
\,.
$
Similarly, for symmetric $M$ and $X$,  we have
$
M\bullet X=\sum_i M_{ii}X_{ii} + 2\sum_{i<j} M_{ij}X_{ij} = \vvec_{\scaleto{2}{4.0pt}}(M)^\top \vvec_{\scaleto{1}{4.0pt}}(X)
\,.
$

In the different subsections, in presenting ADMM algorithms,
the primal variables $x$ and $Z$ (and the associated iterates $x^t$ and $Z^t$), the
Lagrange multiplier $\Psi$ (and the associated iterates $\Psi^t$), and the iterates $Y^t$
have similar uses but different meanings. Throughout, $\theta_\ell$ denotes the $\ell$-th greatest eigenvalue of $\rho Y^{t+1}$, and $\lambda_\ell$ denotes the $\ell$-th greatest eigenvalue of $Z^{t+1}$. 


\section{The 0/1 D-optimality Problem}\label{sec:d-opt}
The 0/1 \emph{D-optimality problem} is
\begin{align*}\label{dopt}\tag{D-Opt}
\textstyle
&\displaystyle\max_x \left\{ \ldet\textstyle\left(\sum_{\ell\in N} \left(v_\ell v_\ell^\top \right) x_\ell \right) ~:~ \mathbf{e}^\top x=s,~ x \in \{0,1\}^n
\right\},
\end{align*}
where $v_\ell \in \mathbb{R}^{m}$, 
for $\ell\in N:=\{1,\ldots,n\}$,
with
$s\geq m$.
The motivation for this model is that the  $n$ points $v_\ell\in\mathbb{R}^m$ are 
potential (costly) ``design points'' for a linear-regression model in $m$ ``factors''.
\ref{dopt} seeks to choose $s$ design points, from the full set of $n$ of them, so as to minimize the  determinant of the covariance matrix
(i.e., the ``generalized variance'')
of the parameter 
estimates in a linear model that would seek to linearly predict responses
based on the chosen $s$ experiments. It turns out that in the Gaussian case, the volume of the standard confidence ellipsoid for the true parameters is
inversely proportional to the determinant of the sum of $v_\ell v_\ell^\top$\,,
over the chosen design points. So, we can see that \ref{dopt} is a truly fundamental problem in the design of experiments.

It is very useful to define  $A:= (v_1, v_2, \dots, v_n)^\top$
(which we always assume has full column rank), and so we have  
$
\sum_{\ell\in N} \left(v_\ell v_\ell^\top \right) x_\ell = A^\top \Diag(x) A$.
Relative to \ref{dopt}, we consider the 
\emph{natural bound}
\begin{align*}\label{prob}\tag{$\mathcal{N}$}
\textstyle
&\displaystyle\max_x \left\{ \ldet \left(A^\top \Diag(x) A \right)~:~ \mathbf{e}^\top x=s,~ 
x\in[0,1]^n
\right\};
\end{align*}
see \cite*{PonteFampaLeeMPB}, and the references therein.
Toward developing an ADMM algorithm for \ref{prob}, we
introduce a variable $Z \in \mathbb{S}^{m}$, and we rewrite \ref{prob}	as  
\begin{alignat}{1}\label{prob:admmdopt1}
& \displaystyle\min_{x,Z}\left\{
-\ldet(Z) ~:~ 
-A^\top \Diag(x) A + Z = 0,~
\mathbf{e}^\top x = s,~
x\in[0,1]^n\right\}.
\end{alignat}
It is easy to see that \eqref{prob:admmdopt1} is basically in a standard form for 
ADMM, 
$\min_{x,z} \{f(x)+g(z) ~:~ Ax + Bz = c\}$ (see, for example, \cite*{boyd2011distributed}). Rather than explicitly put it into the standard form,
we prefer to stay with the form  \eqref{prob:admmdopt1} (here and later), which is 
closer to the natural formulation of the problem. 

The augmented Lagrangian function associated to \eqref{prob:admmdopt1} is 
\begin{align*}
 \mathcal{L}_\rho(x,Z,\Psi,\delta) :=& -\ldet(Z) \!+\! \frac{\rho}{2}\!\left\| -A^\top \Diag(x) A \!+\! Z\!+\! \Psi  \right\|^2_F \\
 &\qquad\!+\! \frac{\rho}{2}\!\left(-\mathbf{e}^\top x \!+\! s + \delta  \right)^2 - \frac{\rho}{2}\!\left\|\Psi\right\|_F^2 -  \frac{\rho}{2}\delta^2, 
\end{align*}
where $\rho >0$ is the penalty parameter and $\Psi \in \mathbb{S}^{m}$, $\delta \in \mathbb{R}$  are the scaled  Lagrangian multipliers. 
Similar to the development of \cite*{Scheinberg2010}
for ``sparse inverse covariance selection'', 
we will apply the ADMM algorithm to \eqref{prob:admmdopt1}\,, by iteratively solving, for $t=0,1,\ldots$, 
\begin{align}
    &x^{t+1}:=\textstyle\argmin_{x \in [0,1]^n}  
\mathcal{L}_\rho(x,Z^{t},\Psi^{t},\delta^t),\label{eq:xmindoptsubpa}\\
    &Z^{t+1}:=\textstyle\argmin_Z ~ \mathcal{L}_\rho(x^{t+1},Z,\Psi^{t},\delta^t),\label{eq:Zmindoptsubpa}\\
    &\textstyle\Psi^{t+1}:=\Psi^{t} - A^\top \Diag(x^{t+1}) A + Z^{t+1},\label{eq:updatepsi}\\
    &\textstyle\delta^{t+1}:=\delta^{t} - \mathbf{e}^\top x^{t+1} +  s.\nonumber
\end{align}
Next, we detail how to solve the subproblems above. 


\subsection{Update \texorpdfstring{$x$}{x}}\label{subsec:update_xdopt}
To update $x$,  we consider  subproblem \eqref{eq:xmindoptsubpa}, more specifically,
\begin{equation*} 
\begin{array}{rl}
x^{t+1}&:= \textstyle \argmin_{x \in [0,1]^n}\left\{
\left\|- A^\top \Diag(x) A \!+\! Z^t\!+\! \Psi^t  \right\|^2_F \!+\! \left(-\mathbf{e}^\top x \!+\! s + \delta^t  \right)^2\right\}\\[6pt] 
&=\textstyle \argmin_{x \in [0,1]^n}\left\{\left\| Hx - d^t\right\|^2_2\right\},
\end{array}
\end{equation*}
where  $d^t:=  \begin{bmatrix}
    \vvec_{\scaleto{\sqrt{2}}{5.5pt}}(Z^t+\Psi^t)\\
    s+\delta^t
\end{bmatrix}_{\strut}$ and $H := \begin{bmatrix}
    G\\
    \mathbf{e}^\top
\end{bmatrix}$, where $G \in \mathbb{R}^{\frac{m(m+1)}{2} \times n}$ is a matrix defined via
$G_{\cdot \ell} := \vvec_{\scaleto{\sqrt{2}}{5.5pt}}\left(v_\ell v_\ell^\top\right)$, for $\ell \in N$.
Then, we have $Gx = \vvec_{\scaleto{\sqrt{2}}{5.5pt}}\left(A^\top \Diag(x) A\right)$.

This is a particular case of  the well-known bounded-variable least-squares (BVLS) problem, and there are several efficient algorithms to solve it; see \cite*{stark1995bounded},  for example.


\subsection{Update \texorpdfstring{$Z$}{Z}}\label{subsec:update_Zdopt}

To update $Z$,   we consider  subproblem \eqref{eq:Zmindoptsubpa}, more specifically, 
\begin{align}
Z^{t+1}:=\textstyle\argmin_Z\left\{ \!-\ldet(Z) \!+\! \frac{\rho}{2}\!\left\|Z \!-\! Y^{t+1}  \right\|^2_F\right\},\label{eq:Zmindoptsubprob}
\end{align}
where $Y^{t+1} :=  A^\top \Diag(x^{t+1}) A \!-\! \Psi^t$. Then we update $Z$ following Proposition \ref{lem:closedformulaupdateZ}.

Using the same ideas as \cite*{Scheinberg2010} (see also \cite*[Section 6.5]{boyd2011distributed}), we have the following result and corollary.

\begin{proposition}
\label{lem:closedformulaupdateZ}
Given $Y^{t+1} \in \mathbb{S}^{m}$ and a positive scalar $\rho$. Let $\rho Y^{t+1} =: Q \Theta Q^\top$ be the eigendecomposition, where $\Theta:= \Diag(\theta_1,\dots,\theta_m)$ and $Q^\top Q = Q Q^\top = I_m$. Then a closed-form optimal solution to \eqref{eq:Zmindoptsubprob} is given by $Z^{t+1} := Q \Lambda Q^\top$ where $\Lambda:=\Diag(\lambda_1\,\ldots,\lambda_m)$ is an $m \times m$ diagonal matrix with 
\[
\lambda_\ell := \left(\theta_\ell + \sqrt{\theta_\ell^2 + 4\rho}\right) \Big/2\rho, \quad \mbox{ for } \ell = 1,\dots,m.
\]
\end{proposition}

\proof 
It suffices to show that $Z^{t+1}$ satisfies the  first-order optimality condition of $\min_{Z} \{ \!-\ldet(Z) \!+\! \frac{\rho}{2}\!\left\|Z \!-\! Y^{t+1}  \right\|^2_F\}$, which is obtained by setting the gradient of the objective function equal to zero, that  is, 
\begin{equation}\label{eq:optcondupdateWdopt}
  -Z^{-1} + \rho(Z - Y^{t+1}) = 0, 
\end{equation}
together with the implicit constraint $Z \succ 0$. We can rewrite  \eqref{eq:optcondupdateWdopt} as 
\begin{align*}
    &\rho Z - Z^{-1} = \rho Y^{t+1} ~\Leftrightarrow~ \rho Z - Z^{-1} =  Q \Theta Q^\top ~\Leftrightarrow~ \rho Q^\top ZQ - Q^\top Z^{-1} Q = \Theta. 
\end{align*}
From the orthogonality of $Q$, we can  verify that the last equation is satisfied by  $Z := Q \Lambda Q^\top$ where $\Lambda:=\Diag(\lambda_1\,\ldots,\lambda_m)$ is an $m \times m$ diagonal matrix such that  $\rho \lambda_{\ell} - 1/\lambda_{\ell} = \theta_\ell$ for $\ell = 1,\dots,m$. Thus, we have 
\[
\lambda_\ell = \frac{\theta_\ell + \sqrt{\theta_\ell^2 + 4\rho}}{2\rho}, \quad \mbox{ for } \ell = 1,\dots,m,
\]
which are always positive, because $\rho > 0$. The result follows.  
\endproof

\begin{corollary} 
\label{cor:psidoptpd}
    Given $x^{t+1}\in\mathbb{R}^n$ and $\Psi^t\in\mathbb{S}^m$, let $Y^{t+1} := A^\top \Diag(x^{t+1}) A - \Psi^t$\,.  For   
    $\rho > 0$, let  
      $\rho Y^{t+1} =: Q \Theta Q^\top$ be the eigendecomposition, where $\Theta := \Diag(\theta_1,\theta_2,\dots,\theta_m)$ with $\theta_1\geq\theta_2\geq \dots \geq \theta_m$ and $Q^\top Q = Q Q^\top = I_m$\,.  Construct $Z^{t+1}$ following Proposition \ref{lem:closedformulaupdateZ}. Then  $\Psi^{t+1}$ computed  by \eqref{eq:updatepsi} is positive definite,  and is given by $Q \Diag(\nu_1,\nu_2,\ldots,\nu_m) Q^\top$ where 
    \[
    \nu_\ell  := \left(-\theta_\ell + \sqrt{\theta_\ell^2 + 4\rho}\right)\Big/ 2\rho,\quad \ell = 1,\dots,m,
    \]
with $\nu_1 \leq \nu_2 \leq \dots \leq \nu_m\,.$
\end{corollary}

\proof
    From \eqref{eq:updatepsi}, we can directly obtain the eigendecomposition of $\Psi^{t+1}$, given the eigendecompositions of $Z^{t+1}$ and $\rho Y^{t+1}$.  Moreover, noticing that the function $f_\rho:\mathbb{R}\rightarrow\mathbb{R}$ defined by $f_\rho(a):=-a + \sqrt{a^2 + 4\rho}$ is decreasing in $a$, we can verify that $\nu_1 \leq \nu_2 \leq \dots \leq \nu_m\,$.   
\endproof


\section{The Maximum-Entropy Sampling Problem}\label{sec:MESP}

Let $C$ be a symmetric positive semidefinite matrix with rows/columns
indexed from $N:=\{1,2,\ldots,n\}$, with $n >1$.
For $0< s < n$,
we define the \emph{maximum-entropy sampling problem}

\begin{align}\tag{MESP}\label{MESP}
\begin{array}{ll}
&z(C,s):=\displaystyle\max_x \left\{\ldet \left(C[S(x),S(x)] \right)~:~ \mathbf{e}^\top x =s,~ x\in\{0,1\}^n\right\},
\end{array}
\end{align}
where $S(x)$ denotes the support of $x\in\{0,1\}^n$, 
 $C[S,S]$ denotes the principal submatrix indexed by $S$.
 For feasibility, we assume that $\rank(C)\geq s$. 
 \ref{MESP} was introduced by \cite*{SW}; also see \cite*{FL2022} and the many references therein. Briefly, in the Gaussian case,  $\ldet (C[S,S])$ is
 proportional to the ``differential entropy'' (see \cite*{Shannon}) of a vector of random variables 
 having covariance matrix $C[S,S]$. So \ref{MESP} seeks to
 find the ``most informative'' $s$-subvector from an $n$-vector following a joint 
 Gaussian distribution. \ref{MESP} finds application in many areas, for example environmental monitoring (see \cite*[Chapter 4]{FL2022}).

In the remainder of this section, we 
develop 
ADMM algorithms for three
well-known convex relaxations of \ref{MESP}: the linx bound, the factorization bound, and the BQP bound.
Although it is not necessary for following most of what we present,
we note that for \ref{MESP}, there are two important general principles that we wish to highlight now, as they are relevant to the bounding methods (see \cite*[Sections 1.5--1.6]{FL2022} for more details):

\smallskip

\begin{itemize}
\item Scaling: For $\gamma>0$, 
$z(C,s) = z(\gamma C,s) -s \ln \gamma$, leading to the equivalent ``scaled problem" (see \cite*{AFLW_Using}). 
\item Complementation: If $\rank(C)=n$, then $z(C,s)=z(C^{-1},n-s) + \ldet C$, leading to the equivalent ``complementary problem" (see \cite*{AFLW_Using}).
\end{itemize}
\smallskip

\noindent The linx bound is invariant under complementation, and the factorization bound is invariant under scaling. But for other combinations of principles and bounding techniques, we can get very different bounds, and this is 
why these two principles are useful. 


\subsection{An ADMM for the linx bound}\label{subsec:linx}

Relative to \ref{MESP},
we consider the 
\emph{(scaled) linx bound}
\begin{align*}\label{prob_linx}\tag{linx$_\gamma$}
\textstyle
&\max_x \left\{\textstyle \frac{1}{2}(\ldet\left(\gamma C\Diag(x) C + \Diag(\mathbf{e}\!-\!x) \right) -s\log(\gamma)) \!~:\!~ \mathbf{e}^\top x\!=\!s,~
x\in[0,1]^n
\right\}, 
\end{align*}
where $C \in \mathbb{S}^n_+$ and $\gamma >0$ is a given scaling parameter.
The linx bound was introduced by \cite*{Kurt_linx}; also see \cite*{FL2022,ChenFampaLee_Fact}. 
The linx bound is convex in $\log(\gamma)$. Exploiting this property, a quasi-Newton method has been proposed to optimize $\gamma$, i.e., to compute the scaling parameter that yields the best possible upper bound (see  \cite*[Section 3.3.5]{FL2022}).

Toward developing an ADMM algorithm for \ref{prob_linx}\,, we
introduce a variable $Z \in \mathbb{S}^{n}$, and we rewrite
\ref{prob_linx}	as 

\begin{alignat}{2}\label{prob:admmlinx1}
& \textstyle \frac{1}{2}\displaystyle\min_{x,Z} ~ && -(\ldet(Z) -s\log(\gamma))\\
& \mbox{\!\!\quad s.t.} \quad &&-\left(\gamma C \Diag(x) C +  \Diag(\mathbf{e}-x)\right)  + Z = 0,\nonumber\\
&&& \mathbf{e}^\top x = s,\nonumber\\
&&& 
x\in[0,1]^n.
\nonumber
\end{alignat}

The augmented Lagrangian function associated to \eqref{prob:admmlinx1} is 
\begin{align*}
\mathcal{L}_\rho(&x,Z,\Psi,\delta):=-\ldet(Z) \!+\! \frac{\rho}{2}\!\left\| -\gamma C \Diag(x) C- \Diag(\mathbf{e}-x)  \!+\! Z \!+\! \Psi    \right\|^2_F \!+\! \frac{\rho}{2}\!\left(-\mathbf{e}^\top x \!+\! s + \delta  \right)^2 \\
&\qquad\qquad\qquad - \frac{\rho}{2}\!\left\|\Psi\right\|_F^2 -  \frac{\rho}{2}\delta^2 +s\log(\gamma), 
\end{align*}
where $\rho >0$ is the penalty parameter and $\Psi \in \mathbb{S}^{n}$, $\delta \in \mathbb{R}$  are the scaled  Lagrangian multipliers. We will apply the ADMM algorithm to \eqref{prob:admmlinx1}\,, by iteratively solving, for $t=0,1,\ldots$, 
\begin{align}
  &x^{t+1}:=\textstyle\argmin_{x\in[0,1]^n} 
\mathcal{L}_\rho(x,Z^{t},\Psi^{t},\delta^t),\label{eq:xminlinxsubpa}\\
    &Z^{t+1}:=\textstyle\argmin_Z ~ \mathcal{L}_\rho(x^{t+1},Z,\Psi^{t},\delta^t),\label{eq:Zminlinxsubpa}\\
    &\textstyle\Psi^{t+1}:=\Psi^{t} - \gamma C \Diag(x^{t+1}) C- \Diag(\mathbf{e}-x^{t+1})  + Z^{t+1},\nonumber\\
    &\textstyle\delta^{t+1}:=\delta^{t} - \mathbf{e}^\top x^{t+1} +  s.\nonumber
\end{align}


\subsubsection{Update \texorpdfstring{$x$}{x}}\label{subsec:update_xlinx}
We consider  subproblem \eqref{eq:xminlinxsubpa}, more specifically,
\begin{align}
x^{t+1}&:= \textstyle\argmin_{x \in [0,1]^n}\!\left\{\!
\left\| -\gamma C \Diag(x) C\!- \!\Diag(\mathbf{e}-x)  \!+\! Z^t \!+\! \Psi^t    \right\|^2_F \!+\!\left(-\mathbf{e}^\top x \!+\! s + \delta^{t}  \right)^2\right\}\nonumber\\
&=\textstyle \argmin_{x \in [0,1]^n}\left\{\left\| Hx - d^t\right\|^2_2\right\},\label{eq:xminlinxsubprob}
\end{align}
where $d^t:=  \begin{bmatrix}
    \vvec_{\scaleto{\sqrt{2}}{5.5pt}}(Z^{t}+\Psi^{t}-I_{n})\\
    s+\delta^{t}
\end{bmatrix}_{\strut}$ 
and $H := \begin{bmatrix}
    G\\
    \mathbf{e}^\top
\end{bmatrix}$, 
where $G \in \mathbb{R}^{\frac{n(n+1)}{2} \times n}$ is a matrix defined via $G_{\cdot \ell} := \vvec_{\scaleto{\sqrt{2}}{5.5pt}}\left(\gamma C_{\ell \cdot}^\top C_{\ell\cdot} - \Diag(\mathbf{e}_\ell) \right)$, for $\ell \in N$.
Then we have $Gx = \vvec_{\scaleto{\sqrt{2}}{5.5pt}}\left(\gamma C \Diag(x) C- \Diag(x) \right)$. 
As mentioned in Section \ref{subsec:update_xdopt}, the update of $x$ consists of solving a BVLS problem.


\subsubsection{Update \texorpdfstring{$Z$}{Z}}\label{subsec:update_Zlinx}
We consider  subproblem \eqref{eq:Zminlinxsubpa}, more specifically, 
\begin{align}
Z^{t+1}:=\textstyle\argmin_Z\left\{ \!-\ldet(Z) \!+\! \frac{\rho}{2}\!\left\|Z \!-\! Y^{t+1}  \right\|^2_F\right\},\label{eq:Zminlinxsubprob}
\end{align}
where $Y^{t+1}:=\gamma C \Diag(x^{t+1}) C+ \Diag(\mathbf{e}-x^{t+1})  - \Psi^t$. Then we update $Z$ following Proposition \ref{lem:closedformulaupdateZ}.


\subsection{An ADMM for the factorization bound}\label{subsec:ddfact}
Relative to \ref{MESP},
we wish to consider the ``factorization bound''; see \cite*{Nikolov,Weijun,FL2022,ChenFampaLee_Fact} and also \cite*{ChenFampaLeeGenScaling,augmentedDDfact}.
The factorization bound has a rather complicated development, and 
we need to go into the details of it, toward developing one of the updates
in the ADMM that we will present. 
The factorization bound is based on
 a fundamental technical lemma of Nikolov. 

\begin{lemma}[\protect{\cite*[Lemma 13]{Nikolov}}]\label{Ni13}
 Let $\lambda\in\mathbb{R}_+^k$ satisfy $\lambda_1\geq \lambda_2\geq \cdots\geq \lambda_k$\,, define $\lambda_0:=+\infty$, and let $s$ be an integer satisfying
 $0<s\leq k$. Then there exists a unique integer $i$, with $0\leq i< s$, such that
 \begin{equation*} 
 \lambda_{i }>\textstyle\frac{1}{s-i }\textstyle\sum_{\ell=i+1}^k \lambda_{\ell}\geq \lambda_{i+1}~.
 \end{equation*}
\end{lemma}
\noindent Although we cannot give an intuition for this lemma, the 
proof in \cite*{Nikolov} is neither long nor hard to follow. 

Next, suppose that  $\lambda\in\mathbb{R}^k_+$ with  
$\lambda_1\geq\lambda_2\geq\cdots\geq\lambda_k$~. Let $\hat\imath$ be the unique integer defined by Lemma \ref{Ni13}. We define
\begin{equation*}
\phi_s(\lambda):=\textstyle\sum_{\ell=1}^{\hat\imath} \log\left(\lambda_\ell\right) + (s - \hat\imath)\log\left(\frac{1}{s-{\hat\imath}} \sum_{\ell=\hat\imath+1}^{k}
\lambda_\ell\right),
\end{equation*}
and, for $X\in\mathbb{S}_{+}^k$~, we define the \emph{$\Gamma$-function}
\begin{equation*} 
\Gamma_s(X):= \phi_s(\lambda(X)).
\end{equation*}

Now suppose that the rank of $C$ is $r\geq s$.
We factorize $C=FF^\top$,
with $F\in \mathbb{R}^{n\times k}$, for some $k$ satisfying $r\le k \le n$. 
This could be a 
Cholesky-type factorization, as in  \cite*{Nikolov} and \cite*{Weijun}, where $F$ is lower triangular and $k:=r$,   it could be derived from
a spectral decomposition  $C=\sum_{i=1}^r \mu_i v_i v_i^\top$\,, by selecting   $\sqrt{\mu_i}v_i$ as the column $i$ of $F$,
 $i=1,\ldots,k:=r$, 
 or it could be derived from 
  the matrix square root of $C$, where $F:=C^{1/2}$, and  $k:=n$.

Finally, we have the \emph{factorization bound} 
\begin{align*}\label{prob_ddfact}\tag{DDFact}
\textstyle
&\displaystyle\max_x \left\{ \Gamma_s(F^\top \Diag(x)F) \, : \, \mathbf{e}^\top x=s,~
x\in[0,1]^n
\right\}. 
\end{align*}
The name ``DDFact'' (from the literature) stems from the fact that it
can be derived from the Lagrangian dual of the Lagrangian dual of
a non-convex formulation, ``Fact''.
In fact, the optimal value of \ref{prob_ddfact} does not depend on which factorization is chosen;
see \cite*[Theorem 2.2]{ChenFampaLee_Fact}. 

Toward developing an ADMM algorithm for \ref{prob_ddfact}, we
introduce a variable $Z \in \mathbb{S}^{k}$, and we rewrite
 \ref{prob_ddfact} as  
\begin{alignat}{1}\label{prob:admmddfact1}
& \displaystyle\min_{x,Z} \left\{ -\Gamma_s(Z) ~:~
-F^\top\Diag(x)F + Z = 0,~
\mathbf{e}^\top x = s,~
x\in[0,1]^n\right\}.
\end{alignat}
The augmented Lagrangian function associated to \eqref{prob:admmddfact1} is 
\begin{align*}
&\mathcal{L}_\rho(x,Z,\Psi,\delta)\!:=
\!-\Gamma_s(Z) \!+\! \frac{\rho}{2}\!\left\|-F^\top\Diag(x)F \!+\! Z \!+\! \Psi  \right\|^2_F \!+\! \frac{\rho}{2}\!\left(-\mathbf{e}^\top x \!+\! s + \delta  \right)^2- \frac{\rho}{2}\!\left\|\Psi\right\|_F^2 -  \frac{\rho}{2}\delta^2\,, 
\end{align*}
where $\rho >0$ is the penalty parameter and $\Psi \in \mathbb{S}^{k}$, $\delta \in \mathbb{R}$  are the scaled  Lagrangian multipliers. We will apply the ADMM algorithm to \eqref{prob:admmddfact1}\,, by iteratively solving, for $t=0,1,\ldots$, 
\begin{align}
    &x^{t+1}:=\textstyle\argmin_{x\in[0,1]^n}  
\mathcal{L}_\rho(x,Z^{t},\Psi^{t},\delta^t),\label{eq:xminddfactsubpa}\\
    &Z^{t+1}:=\textstyle\argmin_Z ~ \mathcal{L}_\rho(x^{t+1},Z,\Psi^{t},\delta^t),\label{eq:Zminddfactsubpa}\\
    &\textstyle\Psi^{t+1}:=\Psi^{t} -F^\top\Diag(x^{t+1})F + Z^{t+1},\label{eq:Psiddfactsubpa}\\ 
    &\textstyle\delta^{t+1}:=\delta^{t} - \mathbf{e}^\top x^{t+1} +s .\nonumber
\end{align}
Next, we detail how to solve the subproblems above.
The $x$ update is rather straightforward, another BVLS problem.
In contrast, the $Z$ update has a closed form, under some technical conditions, and its derivation is rather complicated. 


\subsubsection{Update \texorpdfstring{$x$}{x}}\label{subsec:update_xddfact}
We consider  subproblem \eqref{eq:xminddfactsubpa}, more specifically,
\begin{equation}\label{eq:xminddfactsubprob}
\begin{array}{rl}
x^{t+1}&=\textstyle\argmin_{x \in [0,1]^n}\left\{ 
\left\|-F^\top\Diag(x)F \!+\! Z^t \!+\! \Psi^t  \right\|^2_F \!+\! \left(-\mathbf{e}^\top x \!+\! s + \delta^t  \right)^2\right\}\\[6pt]
&=\textstyle \argmin_{x \in [0,1]^n}\left\{\left\| Hx - d^t\right\|^2_2\right\},
\end{array}
\end{equation}
where  $d^t:= \begin{bmatrix}
   \vvec_{\scaleto{\sqrt{2}}{5.5pt}}(Z^t+\Psi^t)\\
    s+\delta^t
\end{bmatrix}_{\strut}$ and $H := \begin{bmatrix}
    G\\
    \mathbf{e}^\top
\end{bmatrix}$, where $G \in \mathbb{R}^{\frac{k(k+1)}{2} \times n}$ is a matrix defined via
$G_{\cdot \ell} := \vvec_{\scaleto{\sqrt{2}}{5.5pt}}\left( F_{\ell \cdot}^\top F_{\ell\cdot} \right)$, for $\ell \in N$. Then we have $Gx = \vvec_{\scaleto{\sqrt{2}}{5.5pt}}\left(F^\top \Diag(x) F\right)$.
As mentioned in Sections \ref{subsec:update_xdopt} and \ref{subsec:update_xlinx}, the update of $x$ consists of solving a BVLS problem.


\subsubsection{Update \texorpdfstring{$Z$}{Z}}\label{subsec:update_Zddfact}
We consider  subproblem \eqref{eq:Zminddfactsubpa}, more specifically, 
\begin{align}
Z^{t+1}=\textstyle\argmin_Z\left\{ 
\!-\Gamma_s(Z) \!+\! \frac{\rho}{2}\!\left\|Z \!-\! Y^{t+1}  \right\|^2_F\right\},\label{eq:Zminddfactsubprob}
\end{align}
where $Y^{t+1} := F^\top\Diag(x^{t+1})F \!-\! \Psi^t$. In Theorem \ref{thm:updateGamma}, we present a closed-form solution for \eqref{eq:Zminddfactsubprob} under some technical conditions, which we use to update $Z$.
Next, we construct the basis for its derivation.

\begin{proposition}[\protect{\cite*[Proposition 2]{Weijun}}]\label{prop:weijun_supgrad}
   Let $0<s\leq k$ and 
    $Z\in\mathbb{S}_+^k$ 
with rank $r\in[s,k]$.   
   Suppose that the eigenvalues of $Z$ are $\lambda_1 \geq \dots \geq \lambda_r > \lambda_{r+1} = \dots = \lambda_k = 0$ and $Z = Q \Diag(\lambda) Q^\top$ with an orthonormal matrix $Q$. Let 
   $\hat\imath$ be the unique integer  defined by Lemma \ref{Ni13}. Then the subdifferential of the function $\Gamma_s(\cdot)$ at $Z$ 
   is 
   \[\partial \Gamma_s(Z) = Q \Diag(\beta) Q^\top,\]
   where,
\begin{equation*} 
    \begin{array}{lll}
    \displaystyle
        \beta \in \mbox{\rm conv}\Big\{\beta\,:\, &\beta_\ell = 1/\lambda_\ell\,,  & \ell = 1,\dots,\hat\imath;\\
         &\beta_\ell = \displaystyle\frac{s-\hat\imath}{\sum_{j = \hat\imath+1}^k \lambda_j},\quad &\ell =\hat\imath+1,\dots,r;\\
         &\beta_\ell \geq \beta_r\,,&\ell = r+1,\dots,k\Big\}.
    \end{array}
\end{equation*}
\end{proposition}

\begin{lemma}\label{lem:theta_j}
 Let $\theta\in\mathbb{R}^k$ 
 satisfy $\theta_1\geq \theta_2 \geq \cdots\geq \theta_k$\,, define $\theta_0:=+\infty$, let $\rho > 0$, and let $s$ be an integer satisfying $0<s\leq k$. Suppose that
\begin{equation}\label{conditionj}\textstyle\sum_{\ell=s}^k \theta_{\ell} + \sqrt{\left(\textstyle\sum_{\ell=s}^k \theta_{\ell}\right)^2 + 4\rho(k-s+1)} \geq \theta_s + \sqrt{\theta_s^2 + 4\rho}\,.
\end{equation}     Then there exists a unique  
 integer $j$, with $0\leq j< s$, such that
 \begin{equation}\label{reslemiota_a}
 \begin{array}{ll}
 \theta_{j} \!+\! \sqrt{\theta_j^2 + 4\rho}&>\frac{1}{s-j}\left(\displaystyle\sum_{\ell=j+1}^k \theta_{\ell} + \sqrt{\left(\sum_{\ell=j+1}^k \theta_{\ell}\right)^2 + 4\rho(k-j)(s-j)}\,\right)\\[6pt]
 &\geq \theta_{j +1} + \sqrt{\theta_{j+1}^2 + 4\rho}\,.
 \end{array}
 \end{equation}
\end{lemma}

\proof 
Consider the function
\begin{equation*} 
    f_{\rho}(u):=u+\sqrt{u^2+4\rho}\,,
\end{equation*}
    which is increasing in $u$.
    Then
    \[
    f_{\rho}(u)=\frac{4\rho}{-u+\sqrt{u^2+4\rho}}\Rightarrow -u+\sqrt{u^2+4\rho}=\frac{4\rho}{f_{\rho}(u)}\Rightarrow u=\frac{1}{2}\left(f_{\rho}(u) - \frac{4\rho}{f_{\rho}(u)}\right).
    \]

 For $\tau>0$, let 
$u_\tau$ be the value such that 
    $f_{\rho}(u_\tau)=\tau f_{\rho}(u)$.
    Then
    \[
    u_\tau = \frac{1}{2}\left(f_{\rho}(u_\tau) - \frac{4\rho}{f_{\rho}(u_\tau)}\right) = \frac{1}{2}\left(\tau f_{\rho}(u) - \frac{4\rho}{\tau f_{\rho}(u)}\right).
    \]
    As $\frac{4\rho}{f_{\rho}(u)}=f_\rho(u)-2u$, we have 
\[
u_\tau = \frac{1}{2}\left(\tau f_\rho(u)-\frac{1}{\tau}f_\rho(u) +\frac{2u}{\tau}\right) = \frac{u}{\tau} - \frac{1-\tau^2}{2\tau} f_\rho(u).
\]
Then, middle term in the \eqref{reslemiota_a} can be written as
\begin{align*}
&\frac{1}{s-j}\sqrt{(k-j)(s-j)} \; f_\rho\left(\frac{\sum_{\ell=j+1}^{k} \theta_\ell}{\sqrt{(k-j)(s-j)}}\right)\\
& \qquad\qquad   = \sqrt{\frac{k-j}{s-j}} \;f_\rho\left(\frac{\sum_{\ell=j+1}^{k} \theta_\ell}{\sqrt{(k-j)(s-j)}}\right)   ~=~ \sqrt{\frac{k-j}{s-j}}\; f_\rho\left(\sqrt{\frac{k-j}{s-j}} \; \frac{\sum_{\ell=j+1}^{k} \theta_\ell}{k-j}\right).
\end{align*}

Now let $\tau:=1\Big/\sqrt{\frac{k-j}{s-j}}$, 
and consider $\tau f_{\rho}(u)$, for 
$u=\theta_j$ and $u=\theta_{j+1}$\,. We have 
\begin{align*}
    &\tau f_\rho(\theta_j) = f_\rho\left(\sqrt{\frac{k-j}{s-j}}\left(\theta_j - \frac{1-\frac{s-j}{k-j}}{2} f_\rho(\theta_j)\right)\right),\\
    &\tau f_\rho(\theta_{j+1}) = f_\rho\left(\sqrt{\frac{k-j}{s-j}}\left(\theta_{j+1} - \frac{1-\frac{s-j}{k-j}}{2} f_\rho(\theta_{j+1})\right)\right).
\end{align*}

The lemma asks for $j$ such that
\begin{align*}
    &f_\rho(\theta_j) > \frac{1}{\tau} \;f_\rho\left(\sqrt{\frac{k-j}{s-j}} \frac{\sum_{\ell=j+1}^k\theta_\ell}{k-j}\right) \geq f_\rho(\theta_{j+1})\\
    &\qquad\Leftrightarrow f_\rho\left(\sqrt{\frac{k-j}{s-j}} \left(\theta_{j} - \frac{1-\frac{s-j}{k-j}}{2} f_\rho(\theta_j)\right)\right)\\
    &\qquad\qquad\qquad > f_\rho\left(\sqrt{\frac{k-j}{s-j}}\; \frac{\sum_{\ell=j+1}^k\theta_\ell}{k-j}\right) \\
    &\qquad\qquad\qquad \geq f_\rho\left(\sqrt{\frac{k-j}{s-j}} \left(\theta_{j+1} - \frac{1-\frac{s-j}{k-j}}{2} f_\rho(\theta_{j+1})\right)\right).
\end{align*}
Because $f_\rho$ is increasing, this is if and only if
\begin{align}
    &\theta_j - \frac{k-s}{2(k-j)}f_\rho(\theta_j) > \frac{1}{k-j}\sum_{\ell=j+1}^{k}\theta_\ell\geq \theta_{j+1} - \frac{k-s}{2(k-j)}f_\rho(\theta_{j+1})\nonumber \\
    &\Leftrightarrow(k-j)\theta_j - \frac{k-s}{2}f_\rho(\theta_j) > \sum_{\ell=j+1}^{k}\theta_\ell \geq (k-j)\theta_{j+1} - \frac{k-s}{2} f_\rho(\theta_{j+1}). \label{reslemiota_ab}
\end{align}

Let 
\[
\mathcal J := \left\{0\leq j< s~:~\sum_{\ell=j+1}^{k}\theta_\ell \geq (k-j)\theta_{j+1} - \frac{k-s}{2} \; f_\rho(\theta_{j+1})\right\},
\]

Note that $\mathcal J$ is nonempty because the right-hand inequality in \eqref{reslemiota_ab} holds for some $j$  if and only if the right-hand inequality in \eqref{reslemiota_a} holds for the same $j$. As the right-hand inequality in \eqref{reslemiota_a} reduces to \eqref{conditionj} when $j=s-1$, we are assured that $s-1\in\mathcal J$.

Let 
\[
\hat\jmath:=\min\left\{j~:~j\in \mathcal J \right\}.
\]

Next, we show that $\hat \jmath$ is the unique integer, with $0\leq \hat\jmath < s$, for which \eqref{reslemiota_a} holds, or equivalently, for which \eqref{reslemiota_ab} holds.  

\noindent \underline{Case 1: $0\leq j<\hat \jmath$.} Then the right-hand inequality in \eqref{reslemiota_ab} does not hold.
\smallskip

\noindent \underline{Case 2: $j=\hat\jmath$.} If $\hat\jmath=0$, then, because $\theta_0:=+\infty$, the left-hand inequality in \eqref{reslemiota_ab} also holds;
if $\hat\jmath>0$, then
we have
\begin{align*}
&\sum_{\ell=\hat\jmath}^{k}\theta_\ell<\left(k-\left(\hat\jmath-1\right)\right)\theta_{\hat\jmath} - \frac{k-s}{2}f_\rho(\theta_{\hat\jmath})
~\Leftrightarrow~ \sum_{\ell=\hat\jmath+1}^{k}\theta_\ell<\left(k-\hat\jmath\right)\theta_{\hat\jmath} - \frac{k-s}{2}f_\rho(\theta_{\hat\jmath}). 
\end{align*}
So, the left-hand inequality in \eqref{reslemiota_ab} also holds.

\smallskip

\noindent \underline{Case 3:  $\hat\jmath<j<s$.} We will first show that $j\!-\!1\in \mathcal J$, and therefore the right-hand inequality in \eqref{reslemiota_ab} holds for $j\!-\!1$. Using this result, we finally show that the left-hand inequality in \eqref{reslemiota_ab} does not hold for $j$.

To show that $j\!-\!1 \in\mathcal J$, it suffices to show that $i\in\mathcal J$, for 
all $i$ such that $\hat\jmath \leq i <s$. Equivalently, we will show that  if $i\in\mathcal J$,   then $i\!+\!1 \in\mathcal J$, for all $0\leq i <s\!-\!1$. We note that if $i\in\mathcal J$, we have 
\[
\sum_{\ell=i+1}^{k}\theta_\ell \geq (k-i)\theta_{i+1} - \frac{k-s}{2} \; f_\rho(\theta_{i+1}) \Leftrightarrow \sum_{\ell=i+2}^{k}\theta_\ell \geq (k-(i+1))\theta_{i+1} - \frac{k-s}{2} \; f_\rho(\theta_{i+1}).
\]
Then, as $\theta_{i+1}\geq \theta_{i+2}$~,  it suffices to prove that 
\[
g_\rho(u) := (k-(i+1))u - \frac{k-s}{2} \; f_\rho(u) = \frac{k-2(i+1)+s}{2}u - \frac{k-s}{2} \sqrt{u^2+4\rho}
\]
is a non-decreasing function of $u$. 

Let $a := \frac{k + s -2(i+1) }{2}$ and $b := \frac{k-s}{2}$.  Then we have $g_\rho^\prime(u)= a - \frac{bu}{\sqrt{u^2+4\rho}}$. We can easily verify that $a\geq b\geq 0$, and then it is straightforward to see that $g_\rho^\prime(u)\geq 0$, for $u<0$. For $u\geq 0$, we have 
\begin{align*}
     &g_\rho^\prime(u)\geq 0 \Leftrightarrow a - \frac{bu}{\sqrt{u^2+4\rho}} \geq 0\Leftrightarrow
     a^2 ({u^2+4\rho}) \geq b^2u^2\Leftrightarrow (a+b)(a-b) u^2 + 4a^2 \rho \geq 0\,,
\end{align*}
where the two last inequalities also hold because $a\geq b\geq 0$.
We conclude that $g_\rho(u)$ is non-decreasing, and therefore $i\in\mathcal J$, for 
all $i$ such that $\hat\jmath \leq i <s$.
In  particular,  $j\!-\!1\in\mathcal J$,
so
\begin{align*}
&\sum_{\ell=j}^{k}\theta_\ell \geq (k-(j-1))\theta_{j} - \frac{k-s}{2} \; f_\rho(\theta_{j})
~\Leftrightarrow~ \sum_{\ell=j+1}^{k}\theta_\ell \geq (k-j)\theta_{j} - \frac{k-s}{2} \; f_\rho(\theta_{j}), 
\end{align*}
which shows that  the left-hand inequality in \eqref{reslemiota_a} does not hold for $j$.
\endproof

\begin{remark}\label{rem:jhat}
    Note that \eqref{conditionj} is satisfied when  $\textstyle\sum_{\ell=s+1}^k \theta_{\ell} \geq 0$,
because, in this case,
$\textstyle\sum_{\ell=s}^k \theta_{\ell} \geq \theta_s$\,, and because we also have
    $4\rho(k-s+1) \geq  4\rho$.
    In particular, $\theta\geq 0$ implies \eqref{conditionj}.
    Moreover, when $s=k$ and $\theta_{k-1}>\theta_k$, we can verify that \eqref{conditionj} holds as well; also see Lemma \ref{lem:sequalk}. 
\end{remark}

\begin{remark}
   Notice that Lemma \ref{lem:theta_j}  
   becomes Lemma \ref{Ni13} when $\rho=0$. We wish to emphasize that  Lemma  \ref{lem:theta_j} 
   does not follow from 
   Lemma \ref{Ni13} for any
   sequence of $\lambda_\ell$\,. 
   So  Lemma \ref{lem:theta_j} is
   a genuine and subtle extension of  Lemma \ref{Ni13}.
\end{remark}

\begin{lemma}\label{lem:sequalk}
Let $\theta\in\mathbb{R}^k$ 
 satisfy $\theta_1\geq \theta_2 \geq \cdots\geq \theta_k$\,, define $\theta_0:=+\infty$.
 Let $\xi$ $(0\leq \xi\leq k-1)$ be such that
 $\theta_\xi> \theta_{\xi +1}=\cdots=\theta_{k-1}=\theta_k$\,. Let  $\rho>0$.
 For $s=k$, there is a unique $j$ that satisfies \eqref{reslemiota_a}, which is precisely $\xi$.
\end{lemma}

\proof 
    When $s=k$ and $j=\xi$, the middle term in \eqref{reslemiota_a} reduces to $\theta_{\xi +1} +\sqrt{\theta_{\xi +1}^2+4\rho}$. Therefore, we can easily see that both inequalities in \eqref{reslemiota_a} hold. 
\endproof

In Lemma \ref{lem:lambda_j} we will define $\lambda\in \mathbb{R}^k$ such that $\lambda_1 \geq \lambda_2 \geq\dots \geq \lambda_{k}$\,. Later, we will see that $\lambda^{t+1}:=\lambda$ corresponds to the vector of eigenvalues of the closed-form optimal solution $Z^{t+1}$ for \eqref{eq:Zminddfactsubprob}, which we will construct.

\begin{lemma}\label{lem:lambda_j}
    Let $\theta \in \mathbb{R}^k$ with $\theta_1\geq\theta_2\geq \dots \geq \theta_k$\,, $\theta_0:=+\infty$, $\rho>0$, $0 < s \leq k$. Assume that there exists a
    unique $j$  called $\hat\jmath$ that satisfies  \eqref{reslemiota_a}. Define 
    \begin{equation*} 
        \eta:=\eta(\hat\jmath) := 
        {\sum\limits_{\ell= \hat\jmath+1}^{k}\theta_\ell + \sqrt{\left(\sum\limits_{\ell= \hat\jmath+1}^{k}\theta_\ell\right)^2 + 4\rho(k-\hat\jmath)(s-\hat\jmath)}}\,,
    \end{equation*}
    and  $\lambda:=\lambda(\hat\jmath) \in \mathbb{R}^k$ with
\begin{equation}\label{def:lambda}
    \lambda_\ell  := \begin{cases}
    \displaystyle
        \frac{\theta_\ell + \sqrt{\theta_\ell^2 + 4\rho}}{2\rho},\quad &\ell = 1,\dots,\hat\jmath;\\[8pt]
        \displaystyle
        \frac{\theta_\ell}{\rho} + \frac{2(s-\hat\jmath)}{\eta},\quad& \ell = \hat\jmath+1,\dots,k.
    \end{cases}
\end{equation}
Let $\lambda_0:=+\infty$. Then, we have that $\lambda_1 \geq \lambda_2 \geq\dots\geq  \lambda_k\,,$  and $\lambda_{\hat\jmath}  > \frac{\eta}{2\rho(s-\hat\jmath)} \geq \lambda_{\hat\jmath+1} \, .$
\end{lemma}

\proof 
    Because $\theta_1\geq\theta_2\geq \dots\geq\theta_k$\,, we have $\lambda_1 \geq \lambda_2 \geq\dots\geq \lambda_{\hat\jmath}$ and  $\lambda_{\hat\jmath+1} \geq \lambda_{\hat\jmath+2}\geq \dots \geq \lambda_{k}$\,. Now we just need to show that $\lambda_{\hat\jmath} >{\frac{\eta}{2\rho(s-\hat\jmath)} \geq } \lambda_{\hat\jmath+1}$\,. From the left-hand inequality in \eqref{reslemiota_a}, we have that $\theta_{\hat\jmath} + \sqrt{\theta_{\hat\jmath}^2 + 4\rho} > \frac{\eta}{s-{\hat\jmath}}$, then
    $$\lambda_{\hat\jmath} =\frac{\theta_{\hat\jmath} + \sqrt{\theta_{\hat\jmath}^2 + 4\rho}}{2\rho} >  \frac{\eta}{2\rho(s-{\hat\jmath})}.$$
     Now, define $w := \frac{\eta}{s-{\hat\jmath}}$ and note that $w>0$. Then, to show that 
    $$ \frac{\eta}{2\rho(s-{\hat\jmath})} \geq \lambda_{{\hat\jmath}+1}\,,$$ it suffices to verify that
    \begin{align*}
        &\frac{w}{2\rho} \geq \frac{\theta_{{\hat\jmath}+1}}{\rho} + \frac{2}{w} \; \Leftrightarrow \;  w^2 - 2\theta_{{\hat\jmath}+1}w - {4\rho} \geq 0\, \Leftrightarrow\\
        &\left(w- {(\theta_{{\hat\jmath}+1} + \sqrt{\theta_{{\hat\jmath}+1}^2 + 4\rho})}\right)   \left(w- {(\theta_{{\hat\jmath}+1} - \sqrt{\theta_{{\hat\jmath}+1}^2 + 4\rho})}\right)\geq 0.
    \end{align*}
    From the right-hand inequality in \eqref{reslemiota_a}, we have $w \geq\theta_{{\hat\jmath} +1} + \sqrt{\theta_{{\hat\jmath}+1}^2 + 4\rho}$\,. Then
    \[
     w \geq {\theta_{{\hat\jmath}+1} + \sqrt{\theta_{{\hat\jmath}+1}^2 + 4\rho}} \geq {\theta_{{\hat\jmath}+1} - \sqrt{\theta_{{\hat\jmath}+1}^2 + 4\rho}}\,.
    \]
     Therefore, the result follows.   
\endproof

\begin{lemma}\label{lem:zpsd}
Let $\theta\in\mathbb{R}^k$ 
 satisfy $\theta_1\geq \theta_2 \geq \cdots\geq \theta_k$\,, define $\theta_0:=+\infty$.
 Let $\xi$ $(0\leq \xi\leq k-1)$ be such that
 $\theta_\xi> \theta_{\xi+1}=\cdots=\theta_{k-1}=\theta_k$\,. Let  $\rho>0$. Assume that $\theta:=\rho\tilde{\theta}$, that is, $\theta$ varies linearly with $\rho$. 
 Then, for $s=k$, the vector $\lambda$ constructed in Lemma \ref{lem:lambda_j} is nonnegative for all $\rho>0$. 
\end{lemma}

\proof 
If $\theta_{\xi+1}\geq 0$, then $\theta\in \mathbb{R}^k_+$ and the result trivially follows. Therefore, in the following, we consider that $\theta_{\xi+1}<0$.

From Lemma \ref{lem:sequalk}, we know that $\xi$ is the unique  integer  that  satisfies  \eqref{reslemiota_a}, i.e., $\hat{\jmath}=\xi$ in Lemma \ref{lem:lambda_j}.  It is straightforward to see from \eqref{def:lambda} that $\lambda_\ell > 0$ for $\ell=1,\ldots,\hat{\jmath}$. So, it remains to prove that $\lambda_\ell  > 0$ for $\ell=\hat{\jmath}+1,\ldots,k$.  We have that 
\begin{align*}
\lambda_{\ell} & ~=~ \frac{\theta_{\ell}}{\rho} + \frac{2(s-\hat{\jmath})}{\eta} ~=~ \frac{\theta_{\ell}}{\rho} + \frac{2(k-\xi)}{(k-\xi)\theta_{\xi+1} + \sqrt{(k-\xi)^2\theta_{\xi+1}^2 +4\rho(k-\xi)^2}}\\
& ~=~ \frac{\theta_{\ell}}{\rho} + \frac{2}{\theta_{\xi+1} + \sqrt{\theta_{\xi+1}^2 +4\rho}}
~=~ \tilde{\theta}_{\ell} + \frac{2}{\rho \tilde{\theta}_{\xi+1} + \sqrt{\rho^2\tilde{\theta}_{\xi+1}^2 +4\rho}}\,.
\end{align*}
We see from the last expression that $\lambda_{\ell}=\lambda_{\ell}(\rho)$ is a decreasing function of $\rho$.  Also, as  $\hat{\jmath}= \xi$, we have   
from our assumption that $\theta_{\ell} = \theta_{\xi +1}$\,, for all $\ell\geq \hat{\jmath}+1$.
Then, it suffices to show that 
$\lim_{\rho\rightarrow +\infty} \lambda_{\xi+1}(\rho)=0$, which holds because 
\begin{align*}
&\lim_{\rho\rightarrow +\infty} \tilde{\theta}_{\xi+1} + \frac{2}{\sqrt{\rho^2\tilde{\theta}_{\xi+1}^2 +4\rho} + \rho \tilde{\theta}_{\xi+1}}= \lim_{\rho\rightarrow +\infty} \tilde{\theta}_{\xi+1} + \frac{2\left(\sqrt{\rho^2\tilde{\theta}_{\xi+1}^2 +4\rho} - \rho \tilde{\theta}_{\xi+1}\right)}{\left(\sqrt{\rho^2\tilde{\theta}_{\xi+1}^2 +4\rho}\right)^2 - \rho^2 \tilde{\theta}^2_{\xi+1}}\\
&\qquad=\lim_{\rho\rightarrow +\infty} \tilde{\theta}_{\xi+1} + \frac{2\sqrt{\rho^2\tilde{\theta}_{\xi+1}^2 +4\rho} - 2\rho \tilde{\theta}_{\xi+1}}{4\rho} =\lim_{\rho\rightarrow +\infty} \frac{\tilde{\theta}_{\xi+1} + \sqrt{\tilde{\theta}_{\xi+1}^2 +4/\rho} }{2}\\
&\qquad = \frac{\tilde{\theta}_{\xi+1} + |\tilde{\theta}_{\xi+1}| }{2} = 0.
\end{align*}
\vskip-20pt
\endproof

\medskip

In Lemma \ref{lem:sum_lambda}, we show that the $\hat\imath$  defined by Lemma \ref{Ni13} for the $\lambda$ constructed in Lemma \ref{lem:lambda_j}, is precisely the ${\hat\jmath}$ defined by Lemma \ref{lem:theta_j}. This is a key result for the construction of a closed-form solution for \eqref{eq:Zminddfactsubprob} in Theorem \ref{thm:updateGamma}.

\begin{lemma}\label{lem:sum_lambda}
 Let $\theta\in\mathbb{R}^k$ 
 satisfy $\theta_1\geq \theta_2 \geq \cdots\geq \theta_k$\,, define $\theta_0:=+\infty$, let $\rho > 0$, and let $s$ be an integer satisfying $0<s\leq k$.
 Suppose that there exists a unique $j$  called $\hat\jmath$, that satisfies  \eqref{reslemiota_a}, and let  $\lambda$ be defined by Lemma \ref{lem:lambda_j}. 
      Then ${\hat\jmath}$ is the unique integer $\hat\imath$ defined by Lemma \ref{Ni13} for $\lambda$. 
\end{lemma}

\proof 
From Lemma \ref{lem:lambda_j}, we have 
\[
\lambda_{\hat\jmath}>\frac{\eta}{(s-{\hat\jmath})2\rho} \geq \lambda_{{\hat\jmath}+1}\,.
\]
    Now, let $\zeta := \sum_{\ell = {\hat\jmath}+1}^k \theta_\ell$\,. Then, we have  
    \begin{align*}
    \sum_{\ell = {\hat\jmath}+1}^k \lambda_\ell &= \frac{\zeta}{\rho} + \frac{2(k-{\hat\jmath})(s-{\hat\jmath})}{\zeta + \sqrt{\zeta^2 + 4\rho(k-{\hat\jmath})(s-{\hat\jmath})}}\\[5pt]
    &= \frac{\zeta^2 + \zeta\sqrt{\zeta^2 + 4\rho(k-{\hat\jmath})(s-{\hat\jmath})} + 2\rho(k-{\hat\jmath})(s-{\hat\jmath})}{\rho(\zeta + \sqrt{\zeta^2 + 4\rho(k-{\hat\jmath})(s-{\hat\jmath})})}\\[5pt]
    &= \frac{\zeta^2 + 2\zeta\sqrt{\zeta^2 + 4\rho(k-{\hat\jmath})(s-{\hat\jmath})} + \left(\zeta^2 + 4\rho(k-{\hat\jmath})(s-{\hat\jmath})\right)}{2\rho(\zeta + \sqrt{\zeta^2 + 4\rho(k-{\hat\jmath})(s-{\hat\jmath})})}\\[5pt]
    &= \frac{\left(\zeta + \sqrt{\zeta^2 + 4\rho(k-{\hat\jmath})(s-{\hat\jmath})}\right)^2}{2\rho(\zeta + \sqrt{\zeta^2 + 4\rho(k-{\hat\jmath})(s-{\hat\jmath})})} ~=~ \frac{\zeta + \sqrt{\zeta^2 + 4\rho(k-{\hat\jmath})(s-{\hat\jmath})}}{2\rho}~=~ \frac{\eta}{2\rho}.
\end{align*}
Therefore, 
we can see that the unique integer $\hat\imath$ defined  for $\lambda$ by Lemma \ref{Ni13} is exactly ${\hat\jmath}$.   
\endproof

 For clarity, we  omit the dependence of certain parameters on the ADMM iteration  in the remainder of this section. In particular, we adopt the simplified notation $Q:=Q^{t+1}$, $\Theta:=\Theta^{t+1}$, $\theta:=\theta^{t+1}$, $\lambda:=\lambda^{t+1}$, $\beta:=\beta^{t+1}$,  $\nu:=\nu^{t+1}$, $\eta:=\eta^{t+1}$, and $\hat\jmath:=\hat\jmath^{\,t+1}$.

\begin{theorem}\label{thm:updateGamma}
    Given $Y^{t+1}\in \mathbb{S}^k$\,,
    $0 < s \leq k$, and 
    $\rho > 0$. 
     Let $\rho Y^{t+1} =: Q \Theta Q^\top$ be the eigendecomposition, where $\Theta := \Diag(\theta_1,\theta_2,\allowbreak\ldots,\theta_k)$ with $\theta_1\geq\theta_2\geq \dots \geq \theta_k$ and $Q^\top Q = Q Q^\top = I_k$\,. Assume that there exists a unique 
     $j$  called $\hat\jmath$ that satisfies  \eqref{reslemiota_a}.
    Let 
   $\lambda$ be defined as in Lemma \ref{lem:lambda_j} and assume that $\lambda\geq 0$. 
    Then, 
    a closed-form optimal solution to \eqref{eq:Zminddfactsubprob} is given by $Z^{t+1} := Q \Diag(\lambda) Q^\top$.
\end{theorem}

\proof 
Let ${\hat\jmath}$ be the unique integer defined by Lemma \ref{lem:theta_j}.
In Lemma \ref{lem:sum_lambda}, we showed  that ${\hat\jmath}$ is $\hat\imath$  defined by Lemma \ref{Ni13} for $\lambda$.
Therefore, from Proposition \ref{prop:weijun_supgrad}, we have 
that $ Q \Diag(\beta) Q^\top\in \partial \Gamma_s(Z^{t+1}),$ where
\begin{align*} 
    \beta_\ell := \begin{cases}
        \displaystyle\frac{1}{\lambda_\ell},\quad &\ell = 1,\dots,{\hat\jmath};\\[10pt]
        \displaystyle\frac{2\rho(s-{\hat\jmath})}{\eta},\quad& \ell = {\hat\jmath}+1,\dots,k,
    \end{cases}
\end{align*}
where $\eta$ is defined in Lemma 
 \ref{lem:lambda_j}.
 
Let $f(Z) := -\Gamma_s(Z) + \frac{\rho}{2}\|Z-Y^{t+1}\|_F^2$\,. 
Note that 
\begin{align*}
    \partial f(Z^{t+1}) &\ni -Q\Diag(\beta) Q^\top + \rho(Z^{t+1}-Y^{t+1})\\
    &= -Q\Diag(\beta) Q^\top + \rho Q\Diag(\lambda)Q^\top - Q\Theta Q^\top\\
    &= Q\Diag(\rho \lambda - \beta - \theta) Q^\top.
\end{align*}
It suffices to show that $0 \in \partial f(Z^{t+1})$, and hence 
it suffices to show that $\rho \lambda_\ell - \beta_\ell - \theta_\ell = 0$ for $\ell = 1,\dots,k$. For  $\ell = 1,\dots,{\hat\jmath}$, we have
\begin{align*}
    \rho \lambda_\ell - \frac{1}{\lambda_\ell} - \theta_\ell &~=~ \rho\frac{\theta_\ell + \sqrt{\theta_\ell^2 + 4\rho}}{2\rho} - \frac{2\rho}{\theta_\ell + \sqrt{\theta_\ell^2 + 4\rho}} - \theta_\ell\\[10pt]
    &~=~ \frac{\left({\theta_\ell + \sqrt{\theta_\ell^2 + 4\rho}}\right)^2 - \left(\theta_\ell^2 + 2\theta_\ell \sqrt{\theta_\ell^2 + 4\rho} + (\theta_\ell^2+4\rho)\right) }{2\theta_\ell + 2\sqrt{\theta_\ell^2 + 4\rho}}\\[10pt]
    &~=~ \frac{\left({\theta_\ell + \sqrt{\theta_\ell^2 + 4\rho}}\right)^2 - \left({\theta_\ell + \sqrt{\theta_\ell^2 + 4\rho}}\right)^2}{2\theta_\ell + 2\sqrt{\theta_\ell^2 + 4\rho}} ~=~ 0.
\end{align*}
For $\ell = {\hat\jmath}+1,\dots,k$, we have
\begin{align*}
     \rho \lambda_\ell - \frac{2\rho(s-{\hat\jmath})}{\eta} - \theta_\ell &= \rho\frac{\theta_\ell}{\rho} + \rho\frac{2(s-{\hat\jmath})}{\eta}   - \frac{2\rho(s-{\hat\jmath})}{\eta} - \theta_\ell = 0,
\end{align*}
and therefore $0 \in \partial f(Z^{t+1})$.   
\endproof

\begin{corollary}\label{cor:psiddfactpd}
    Given $x^{t+1}\in\mathbb{R}^n$ and $\Psi^t\in\mathbb{S}^k$, let $Y^{t+1} := F^\top \Diag(x^{t+1}) F - \Psi^t$\,. For 
    $\rho > 0$, 
     let $\rho Y^{t+1} =: Q \Theta Q^\top$ be the eigendecomposition, where $\Theta := \Diag(\theta_1,\theta_2,\dots,\theta_k)$ with $\theta_1\geq\theta_2\geq \dots \geq \theta_k$ and $Q^\top Q = Q Q^\top = I_k$\,.  Assume that there exists a (unique) $j$ 
     called $\hat\jmath$ that satisfies  \eqref{reslemiota_a}, and construct $Z^{t+1}$ following Theorem \ref{thm:updateGamma}. Then  $\Psi^{t+1}$, computed by \eqref{eq:Psiddfactsubpa}, is positive definite  and is given by $Q \Diag(\nu) Q^\top$, where 
    \begin{equation*} 
    \nu_\ell  := \begin{cases}\displaystyle
        \frac{-\theta_\ell + \sqrt{\theta_\ell^2 + 4\rho}}{2\rho},\quad &\ell = 1,\dots,\hat\jmath;\\[7pt]
        \displaystyle
        \frac{2(s-\hat\jmath)}{\eta},\quad& \ell = \hat\jmath+1,\dots,k,
    \end{cases}
\end{equation*}
with $\nu_1 \leq \nu_2 \leq \dots \leq \nu_k\,,$ where $\eta$ is defined in Lemma 
 \ref{lem:lambda_j}.
\end{corollary}

\proof 
    From \eqref{eq:Psiddfactsubpa}, we have $ \Psi^{t+1} := \Psi^t - F^\top \Diag(x^{t+1})F +Z^{t+1} = Z^{t+1} -Y^{t+1}$. Following the construction of $Z^{t+1}$ using $\lambda$ defined in Lemma \ref{lem:lambda_j}, we have $\Psi^{t+1} = Q \Diag(\lambda - \frac{1}{\rho}\theta) Q^\top$,  then we define $\nu := \lambda - \frac{1}{\rho}\theta$.  Note that for $\ell = 1,\dots,\hat\jmath$, we have
    $$\nu_\ell = \frac{\theta_\ell + \sqrt{\theta_\ell^2 + 4\rho}}{2\rho} - \frac{\theta_\ell}{\rho} =  \frac{-\theta_\ell + \sqrt{\theta_\ell^2 + 4\rho}}{2\rho} \,,$$
   and for $\ell = \hat\jmath+1,\dots,k$. we have
    $$\nu_\ell = \frac{\theta_\ell}{\rho} + \frac{2(s-\hat\jmath)}{\eta} - \frac{\theta_\ell}{\rho} =  \frac{2(s-\hat\jmath)}{\eta}\,.$$
    Also, we note that because $\rho > 0$ and $0 \leq \hat\jmath < s \leq k$, then $\nu > 0$. Finally, we note that the function $f_\rho:\mathbb{R}\rightarrow\mathbb{R}$, defined by  $
    f_{\rho}(a):=-a+\sqrt{a^2+4\rho}$, 
    is decreasing in $a$, so $\nu_1 \leq \dots \leq \nu_{\hat\jmath}\,.$ Then, it suffices to show that  $\frac{-\theta_{\hat\jmath} + \sqrt{\theta_{\hat\jmath}^2 + 4\rho}}{2\rho} \leq 2(s-\hat\jmath)/\eta\,.$ Suppose instead that $$\frac{-\theta_{\hat\jmath} + \sqrt{\theta_{\hat\jmath}^2 + 4\rho}}{2\rho} > \frac{2(s-\hat\jmath)}{\eta}\,.$$
    From Lemma \ref{lem:lambda_j},  we have $\lambda_{\hat\jmath} = \frac{\theta_{\hat\jmath} + \sqrt{\theta_{\hat\jmath}^2 + 4\rho}}{2\rho} > \frac{\eta}{2\rho(s-\hat\jmath)} ~\Leftrightarrow~  \frac{2(s-\hat\jmath)}{\eta} > \frac{2}{\theta_{\hat\jmath} + \sqrt{\theta_{\hat\jmath}^2 + 4\rho}}\,.$
    Then, we have 
    \begin{align*}
        &\frac{-\theta_{\hat\jmath} + \sqrt{\theta_{\hat\jmath}^2 + 4\rho}}{2\rho} > \frac{2}{\theta_{\hat\jmath} + \sqrt{\theta_{\hat\jmath}^2 + 4\rho}}~\Leftrightarrow~
         -\theta^2_{\hat\jmath} + \theta_{\hat\jmath}^2 + 4\rho > 4\rho ~\Leftrightarrow~
          4\rho > 4\rho\,.
    \end{align*}
    This contradiction completes the proof.
\endproof


\subsection{An ADMM for the BQP bound}\label{subsec:BQP}
Relative to \ref{MESP},
we consider the 
\emph{(scaled) BQP bound}
\begin{align}\label{bqp_original}\tag{BQP$_\gamma$}
\textstyle
\displaystyle\max_{x,X} \{\ldet&\left(\gamma C\circ X + \Diag(\mathbf{e}-x) \right) - s\log(\gamma) \, :\\& \, \mathbf{e}^\top x\!=\!s,\,X\mathbf{e}\!=\!sx,\,x\!=\!\diag(X),\,X\!\succeq\!xx^\top\},\nonumber
\end{align}
where $C \in \mathbb{S}^n_+$\,,  and $\gamma>0$ is a given scaling parameter. 
 The BQP bound is convex in $\log(\gamma)$ and a quasi-Newton method can be used to find the optimal $\gamma$, similar to the approach used for the linx bound (see \cite*[Section 3.6.5]{FL2022}).
The BQP bound was introduced in \cite*{Anstreicher_BQP_entropy}; also see \cite*{FL2022}.  
Because of the matrix variable, 
experimentation with the BQP bound has been limited.
So a strong motivation of ours in developing an ADMM algorithm
for the BQP bound is
to be able to apply it to larger instances than were heretofore possible.

Toward developing an ADMM algorithm for  \ref{bqp_original}\,, we
introduce the variables $W,E,Z \in \mathbb{S}^{n+1}$, and we rewrite \ref{bqp_original} as 
\begin{alignat}{2}
& \displaystyle\min_{E,W,Z}\quad && -\ldet(Z) + s\log(\gamma)\nonumber\\
 & \mbox{\!\!\quad s.t.} \quad && -( \tilde C\circ W + I_{n+1}) + Z =0,\nonumber\\
 &&&  W - E = 0,\label{prob:bqp}\\
 &&& g_\ell  - G_\ell  \bullet  W = 0, \quad \ell  = 1,\dots,2n + 2,\nonumber\\
&&& W,Z \in \mathbb{S}^{n+1}, ~ E \in \mathbb{S}^{n+1}_+, \nonumber
\end{alignat}
where $\tilde C := \begin{bmatrix}
    0 ~&~ \mathbf{0}^\top\\
    \mathbf{0} ~&~ \gamma C-I_n
\end{bmatrix}\in \mathbb{S}^{n+1}$, $W := \begin{bmatrix}
    1 ~&~ x^\top\\
    x ~&~ {X}
\end{bmatrix}\in \mathbb{S}^{n+1}$, and 
for $\ell=1,\ldots,n$, 
\[
G_\ell:=\left[
\begin{array}{cc}
0 & -\frac{1}{2}\mathbf{e}_\ell^\top \\[3pt] 
 -\frac{1}{2}\mathbf{e}_\ell & \mathbf{e}_\ell\mathbf{e}_\ell^\top
\end{array}
\right]
,~ g_\ell:= 0; \quad  G_{\ell + n}:=\textstyle\frac{1}{2}\left[
\begin{array}{cc}
0 & -s\mathbf{e}_\ell^\top \\[3pt] 
 -s\mathbf{e}_\ell &  
 \mathbf{e}_\ell \mathbf{e}^\top + \mathbf{e}\mathbf{e}_\ell^\top
\end{array}
\right]
,~ g_{\ell + n}:= 0;  
\]
and
\[
G_{2n+1}:=\textstyle\frac{1}{2}\left[
\begin{array}{ll}
0 & \mathbf{e}^\top \\[3pt] 
 \mathbf{e} & 0
\end{array}
\right]
,~ g_{2n+1}:= s; \quad  G_{2n+2}:=\left[
\begin{array}{ll}
1 & \mathbf{0}^\top \\[3pt] 
 \mathbf{0} & 0
\end{array}
\right]
,~ g_{2n+2}:= 1.
\]
The correctness of reformulation \eqref{prob:bqp} of \ref{bqp_original} can be established by verifying the following equivalences:
\begin{align*}
    & \tilde C\circ W + I_{n+1} = \gamma C\circ X + \Diag(\mathbf{e}-x);\\
    &x-\diag(X)=0 \,\Leftrightarrow\, g_\ell  - G_\ell  \bullet  W =0, \mbox{ for } \ell=1,\ldots,n;\\
    &sx-X\mathbf{e}=0 \,\Leftrightarrow\, g_{\ell+n}  - G_{\ell+n}  \bullet  W =0, \mbox{ for } \ell=1,\ldots,n;\\
    &s-\mathbf{e}^\top x=0 \,\Leftrightarrow\, g_{2n+1}  - G_{2n+1}  \bullet  W =0; \\
     &1-W_{11}=0 \,\Leftrightarrow\, g_{2n+2}  - G_{2n+2}  \bullet  W =0.
\end{align*}

The augmented Lagrangian function associated to \eqref{prob:bqp} is
\begin{align*}
\mathcal{L}_\rho(&W,E,Z,\Psi,\Phi,\omega):=\!-\ldet(Z) \!+\! \frac{\rho}{2}\!\left\|Z \!-\!\tilde{C}\circ {W} \!- \!I_{n+1}  \!+\! \Psi  \right\|^2_F  \!+\! \frac{\rho}{2}\!\left\|W \!-\!{E}  \!+\! \Phi  \right\|^2_F  \\
&\qquad\qquad\qquad\qquad \!+\! \sum_{\ell =1}^{2n+2 }\!\frac{\rho}{2}\!\left(g_\ell  \!-\! G_\ell \bullet  W  \!+\! \omega_\ell   \right)^2 - \frac{\rho}{2}\!\left\|\Psi\right\|_F^2 - \frac{\rho}{2}\!\left\|\Phi\right\|_F^2 -  \frac{\rho}{2}\|\omega\|_2^2 +s\log(\gamma), 
\end{align*}
where $\rho >0$ is the penalty parameter and 
$\Psi, \Phi \in \mathbb{S}^{n+1}$, $\omega \in \mathbb{R}^{2n+2}$
are the scaled  Lagrangian multipliers. We will apply the ADMM method to \eqref{prob:bqp}\,, by iteratively solving, for $t=0,1,\ldots$, 
\begin{align}
&W^{t+1}:=\textstyle\argmin_{W} ~ \mathcal{L}_\rho(W,E^t,Z^{t},\Psi^{t},\Phi^t,\omega^t),\label{eq:Wminbqpsubpaa}\\
    &(E^{t+1},Z^{t+1}):=\textstyle\argmin_{E \succeq 0,Z} ~ \mathcal{L}_\rho(W^{t+1},E,Z,\Psi^{t},\Phi^t,\omega^t),\label{eq:EZminbqpsubpaa}\\
    &\textstyle\Psi^{t+1}:=\Psi^{t}+ Z^{t+1} -\tilde{C}\circ W^{t+1} - I_{n+1}\,,\nonumber\\
     &\textstyle\Phi^{t+1}:=\Phi^{t}+ W^{t+1} -E^{t+1},\nonumber\\
     &\textstyle\omega^{t+1}_\ell :=\omega^{t}_\ell + g_\ell  -G_\ell \bullet W^{t+1},\quad \ell  =1,\dots,2n+2.\nonumber
\end{align}


\subsubsection{Update \texorpdfstring{$W$}{W}}\label{subsec:update_Wbqp3}
To update $W$,  we consider  subproblem \eqref{eq:Wminbqpsubpaa}, more specifically, 
\begin{align}\label{eq:subprobsubabqp3}
W^{t+1}:=& \textstyle
\argmin_{W} \displaystyle\left\{\!
\left\|\tilde{C}\circ {W} \!-\! (Z^t \!+\! \Psi^t \!-\!I_{n+1})  \right\|^2_F
\!+\! 
\left\|{W} - (E^t - \Phi^t)  \right\|^2_F \right.\\
&\left.\qquad \textstyle \!+\! \sum_{\ell =1}^{2n+2 }\!\left(   g_\ell  \!-\! G_\ell \bullet  W  \!+\! \omega_\ell ^t\right)^2
\vphantom{\left\|\tilde{C}\circ {W} \!-\! (Z^t \!+\! \Psi^t \!-\!I_{n+1})  \right\|^2_F}
\right\}.
\nonumber
\end{align}
We can verify that \eqref{eq:subprobsubabqp3} is equivalent to 
the  least-squares problem $\min_{ u}\{\left\| H u - d^{t}\right\|^2_2\}$, where
\[
H := \begin{bmatrix}
    \Diag(\vvec_{\scaleto{\sqrt{2}}{5.5pt}}(\tilde{C}))\\
    \Diag(\vvec_{\scaleto{\sqrt{2}}{5.5pt}}(J))\\
    \vvec_{\scaleto{2}{4.0pt}}(G_1)^\top\\
    \vdots\\
    \vvec_{\scaleto{2}{4.0pt}}(G_{2n})^\top\\
    \vvec_{\scaleto{2}{4.0pt}}(G_{2n+1})^\top\\
        \vvec_{\scaleto{2}{4.0pt}}(G_{2n+2})^\top
\end{bmatrix}, \quad d^t:= \begin{bmatrix}
    \vvec_{\scaleto{\sqrt{2}}{5.5pt}}(Z^t + \Psi^t - I_{n+1})\\
   \vvec_{\scaleto{\sqrt{2}}{5.5pt}}(E^t - \Phi^t)\\
   \omega^t_{1}\\
   \vdots\\
   \omega^t_{2n}\\
    \omega^t_{2n+1} + s\\
    \omega^t_{2n+2} +1
\end{bmatrix}, \quad u:=\vvec_{\scaleto{1}{4.0pt}}(W).
\]
 
We note that the least-squares problem $\min_{ u}\{\left\| H u - d^{t}\right\|^2_2\}$ has a  closed-form solution, and 
that the solution is unique because $H$ is full-column rank;
moreover, we note that $H$ does not change during the ADMM iterations. Therefore, we compute the Cholesky factor of the coefficient matrix associated to the normal equations of the least-squares problem only once, and we use it at each iteration of the ADMM algorithm to solve \eqref{eq:subprobsubabqp3}.


\subsubsection{Update \texorpdfstring{$E$}{E} and \texorpdfstring{$Z$}{Z}}\label{subsec:update_EZbqp3}
To update $E$ and $Z$,  we consider  subproblem \eqref{eq:EZminbqpsubpaa}, more specifically,
\begin{equation}\label{eq:EZminsdpsubproba}
\begin{array}{rl}
(E^{t+1},Z^{t+1}):=& \textstyle\argmin_{E\succeq 0,Z}\left\{ \!-\ldet(Z) \!+\! \frac{\rho}{2}\!\left\|Z \!-\!  (\tilde{C}\circ W^{t+1}  + I_{n+1} - \Psi^t)  \right\|^2_F \right.\\[5pt]
&\qquad \left.+ \left\|E - ({W}^{t+1} +  \Phi^t)  \right\|^2_F\right\}.
\end{array}
\end{equation}

We note that the minimization problem in \eqref{eq:EZminsdpsubproba} is separable with respect to $E$ and $Z$, and so these updates can be done in parallel. 

\begin{itemize}

\item  To update $E$, we consider the  subproblem  
\begin{equation}\label{eq:Eminsdpsubproba} 
\begin{array}{rl}
E^{t+1}&:= \textstyle\argmin_{E\succeq 0}\left\{\!\left\|E - {Y}^{t+1}  \right\|^2_F\right\},
\end{array}
\end{equation}
where ${Y}^{t+1}:={W}^{t+1} +  \Phi^t$. Then, we update $E$ following Theorem \ref{thm:closedformulaupdateEa}.

\begin{theorem}[\protect{\cite*[Theorem 2.1]{higham1988computing}}]\label{thm:closedformulaupdateEa}
    Given ${Y}^{t+1} \in \mathbb{S}^{n+1}$. Let ${Y}^{t+1} =: Q \Theta Q^\top$ be the eigendecomposition, where $\Theta:= \Diag(\theta_1,\dots,\theta_{n+1})$ and $Q^\top Q = Q Q^\top = I_{n+1}$\,. Then a closed-form  solution to \eqref{eq:Eminsdpsubproba} is given by $E^{t+1} := Q \Lambda Q^\top$ where $\Lambda:=\Diag(\lambda_1\,\ldots,\lambda_{n+1})$ 
    and $\lambda_\ell := \max(\theta_\ell,0),$ for $\ell = 1,\dots,n+1.$
\end{theorem}

\item  To update $Z$, we consider  the subproblem
\begin{align}
Z^{t+1}:=\textstyle\argmin_Z\left\{ \!-\ldet(Z) \!+\! \frac{\rho}{2}\!\left\|Z \!-\! {Y}^{t+1}  \right\|^2_F\right\},\nonumber
\end{align}
where ${Y}^{t+1}:=  \tilde{C}\circ W^{t+1}  + I_{n+1} - \Psi^t$. Then, we update $Z$ following Proposition \ref{lem:closedformulaupdateZ}.

\end{itemize}


\section{Numerical Experiments}\label{sec:numexp}

In this section, we evaluate our proposed ADMM algorithms for the relaxations \ref{prob} of \ref{dopt}, and \ref{prob_ddfact} and \ref{bqp_original} of \ref{MESP}, comparing them with general-purpose solvers.
The choice of a good penalty parameter $\rho$, 
for augmented-Lagrangian methods like ADMM, 
is critical for practical performance.
For our experiments designed for ``proof of concept'', 
we found good values, which we tabulate in Appendix \ref{app:performance}. We can see that for each group of problems,
these good choices for $\rho$ trend in a predictable manner. This bodes well for us in our motivating context of B\&B; see Section \S\ref{sec:concl} for more extensive comments on this point. 


\subsection{Our computational framework}


\subsubsection{Solvers, parameter settings and computational aspects}

We selected the general-purpose solvers KNITRO (see \cite*{KNITRO}), MOSEK (see \cite*{MOSEK}), and SDPT3 (see \cite*{SDPT3}), which are commonly used in the literature for the kind of problems we solve.  All our algorithms were implemented in Julia v1.11.3, except the code that calls SDPT3, which was implemented in MATLAB R2023b. 
We used the parameter settings for the solvers aiming at their best performance, considering tolerances similar to those used in our ADMM algorithms. Next, we summarize the settings that we employed, so that it is possible to reproduce our experiments.
For KNITRO, we employed KNITRO 14.0.0 (via the Julia wrapper KNITRO.jl v0.14.4), using \texttt{CONVEX = true}, \texttt{FEASTOL} = $10^{-6}$ (feasibility tolerance), \texttt{OPTTOLABS} = $0.05$ (absolute optimality tolerance),
\texttt{ALGORITHM} = 1 (Interior/Direct algorithm), 
\texttt{HESSOPT} = 6 (KNITRO computes a limited-memory quasi-Newton BFGS Hessian; we used the default value of 
\texttt{LMSIZE} = 10 limited-memory pairs stored when approximating the Hessian). 
For MOSEK, we employed MOSEK 10.2.15 (via the Julia wrapper MOSEKTools.jl  v0.15.5), with  \texttt{MSK\_DPAR\_INTPNT\_CO\_TOL\_REL\_GAP} = $0.05$ (relative gap used by the interior-point optimizer for conic problems) and \texttt{MSK\_DPAR\_INTPNT\_CO\_TOL\_DFEAS}  = $0.05$ (dual-feasibility tolerance used by the interior-point optimizer for conic problems).  We note that we used the default primal feasibility tolerance of $10^{-8}$ for MOSEK, even though it is tighter than the one used for the other solvers and for our ADMM, because loosening the feasibility tolerance did not lead to good convergence behavior for MOSEK.  For SDPT3, we used SDPT3 4.0, with \texttt{gaptol} = $10^{-4}$, \texttt{inftol} = $10^{-5}$.

We also experimented with two open-source Julia implementations of first-order methods, namely FrankWolfe.jl (see \cite*{frankwolfeIJOC})  and COSMO.jl, an ADMM-algorithm for convex conic problems (see \cite*{COSMO}). 
For FrankWolfe.jl, we set the parameters \texttt{max\_iteration}=$10^4$ and \texttt{epsilon}=$5\cdot 10^{-2}$ (the ``Frank-Wolfe gap''). To handle the constraints, FrankWolfe.jl calls a generic solver from MathOptInterface.jl (MOI), which we select to be KNITRO. For COSMO.jl,  we set the maximum number of ADMM iterations to infinity, \texttt{eps\_abs} =$10^{-4}$ (absolute tolerance), \texttt{eps\_rel} =$10^{-5}$ (relative tolerance). Both of these first-order methods did not work well on our problems, as we can see with the  detailed results presented in Appendix \ref{app:performance}.

We set a time limit of 1 hour to solve each instance using each procedure tested.

We ran our experiments on ``zebratoo'', a
32-core machine (running Windows Server 2022 Standard):
two Intel Xeon Gold 6444Y processors running at 3.60GHz, with 16 cores each, and 128 GB of memory.


\subsubsection{Optimality tolerance and gap analysis}

In all of our experiments, we obtain solutions for the relaxations within the absolute optimality tolerance of $0.05$. We note that this is a sufficient precision for applying the upper bounds inside a B\&B algorithm, which is our motivating use case, as $0.05$ is not significant when compared to the differences between the upper bounds and the best known solution values for the instances considered of \ref{dopt} and \ref{MESP}. 
These differences (``D-Opt gap'' and ``MESP-gap'') are presented in Appendix \ref{app:performance}\footnote{Here and throughout, consistent with the literature (see, for example, \cite*[Proposition 1.1.1 and Remark 1.1.2]{FL2022}, we consider absolute gaps rather than relative gaps, because the $\ldet(\cdot)$ objectives are not generally nonnegative.}. The best known solutions for \ref{dopt} and \ref{MESP} were obtained with local-search procedures from  \cite*{PonteFampaLeeMPB} and \cite*{KLQ}, respectively. Of course, as a B\&B would proceed, we can expect to eventually see small gaps, and for such relevant B\&B subproblems, one could seek more accurate solutions. 


\subsubsection{Primal and dual feasibility}

We note that despite the optimality tolerance, the bounds computed for \ref{dopt} and \ref{MESP} are genuine bounds, as they are derived from the objective values of dual-feasible solutions. For \ref{prob}, a dual-feasible solution is derived in closed form (see, for example,  \cite*[Section 2]{PonteFampaLeeMPB}). Similarly, for \ref{prob_ddfact}, the dual-feasible solution admits a closed-form expression (see, for example,  \cite*[Section 3.4.4.1]{FL2022}). In the case of \ref{bqp_original}\,, the dual-feasible solution is obtained by solving a simple semidefinite program (see, for example, \cite*[Section 3.6.4]{FL2022}). To solve this semidefinite program, we use the HYPATIA solver (see \cite*{coey2022solving}), which we have found to be very efficient and convenient for these relatively simple semidefinite problems. Additional details on the construction of these solutions are provided in Appendix \ref{sec:appdual}. 

With the approaches adopted, dual-feasible solutions for \ref{prob} and the \ref{prob_ddfact} bound can be computed very efficiently. Accordingly, we begin computing them from the start of the ADMM execution. For \ref{prob} it is computed every 25 iterations, whereas for the \ref{prob_ddfact} bound, owing to its typically fast convergence, it is computed every 5 iterations. In contrast, for the \ref{bqp_original} bound we must solve a semidefinite problem to obtain dual-feasible solutions. Although these semidefinite problems are relatively simple, solving them from the beginning of the ADMM execution would be computationally too expensive and not worthwhile. For this reason, we delay solving them until iteration 1000, when the algorithm is more likely to have reached a sufficiently small gap to satisfy the stopping criterion. From that point on, we solve the semidefinite problem every 50 iterations.

All dual-feasible solutions are constructed from primal-feasible solutions of the corresponding relaxations. We obtain rigorous primal-feasible solutions by projecting the approximate primal solutions, produced either by ADMM or by the solvers used for the relaxations, onto the feasible sets. This is accomplished via an alternating projection algorithm (see, for example,  \cite*{cheney1959proximity}), which is applied until a feasibility tolerance of $10^{-5}$ is achieved.

Thus, throughout the ADMM iterations, we periodically  project the approximate primal solutions onto the feasible sets and  compute dual-feasible solutions. The ADMM procedure is terminated once the duality gap, defined as the difference between the objective values of the dual-feasible solution and the projected primal solution, falls below $0.05$.


\subsubsection{Implementation  details}\label{subsec:implementation}

We initialize the variables for the ADMM procedures as follows:
\begin{itemize}
    \item 
For \ref{prob}, we set $\Psi^0:=0$, $\delta^0:=0$, and $Z^0:=(s/n) A^\top A$ (equivalently, $Z^0:=A^\top \Diag(\bar{x}) A$, where $\bar{x}:= (s/n)\mathbf{e}$). 
\item 
For \ref{prob_ddfact}, we set $\Psi^0:=0$, $\delta^0:=0$, and $Z^0:= \Diag((\lambda_1(C),\ldots,\lambda_k(C)^\top))$, (equivalently, $Z^0:=F^\top \Diag(\bar{x})F$, with $C := Q \Diag((\lambda_1(C),\ldots,\lambda_k(C)^\top)) Q^\top$ (the spectral decomposition of $C$, where $k:=\rank(C)$, $Q\in\mathbb{R}^{n\times k}$),   $F := Q\Diag((\lambda_1(C),\ldots,\lambda_k(C)^\top)^{1/2})$, and $\bar{x} := (s/n)\mathbf{e}$.
\item 
For \ref{bqp_original}, we set $\omega^0:=\mathbf{e}$ and initialize all remaining variables to $0$.
\end{itemize}
  
For the BVLS subproblems (see \S\S\ref{subsec:update_xdopt}, \ref{subsec:update_xlinx}, and \ref{subsec:update_xddfact}), we took only one gradient-direction step, 
and then we projected the solution onto the domain $[0,1]^n$, which worked very well as a heuristic to speed up the iterations. Although not directly applicable to this heuristic, 
we note that there is some theory for convergence of inexact updates within ADMM (see \cite*{EckWang2018,EckYao}, for example). 
We note that the approximation used to solve the BVLS subproblems is motivated by empirical evidence obtained from extensive computational experiments with the ADMM procedures for both the natural bound and the DDFact bound across test instances of diverse origins. In particular, we compared the use of a single gradient step in the BVLS subproblems against exact solutions obtained with Gurobi, as well as against variants involving multiple gradient steps.
The numerical results consistently indicate that the best overall performance of the ADMM schemes is achieved when employing the adopted approximation.

A bottleneck of the ADMM algorithms that we propose is the eigendecomposition of a matrix $\rho Y^{t+1}$ at each iteration, to update a matrix variable $Z$ (see \S\S\ref{subsec:update_Zdopt}, \ref{subsec:update_Zlinx}, \ref{subsec:update_Zddfact}, and \ref{subsec:update_EZbqp3}
).  
We note that the dimension of $\rho Y^{t+1}$ varies for each relaxation. For the natural bound \ref{prob} for \ref{dopt}, the dimension is $m$, which is generally small compared to $n$ in applications of the problem. For the \ref{prob_ddfact} bound for \ref{MESP}, we  choose the dimension $k$ to be the rank of $C$, 
 which makes the ADMM algorithm very effective for low-rank covariance matrices. For the \ref{prob_linx} bound, the dimension is $n$ which ends up making the ADMM algorithm less competitive. For the \ref{bqp_original} bound, the dimension is $n+1$, but as we will see, the ADMM algorithm is competitive with the alternatives 
 for this relaxation. We note that for \ref{bqp_original}\,,
we also have an eigendecomposition to carry out for the $E$ update, but this can be done in parallel with the $Z$ update (see \S\ref{subsec:update_EZbqp3}).


\subsection{Experiments for the D-optimality problem}
We conducted experiments with four types of test instances for the ADMM algorithm described in \S\ref{sec:d-opt}, to compute the natural bound from \ref{prob} to \ref{dopt}, and compare the performance of the ADMM algorithm to KNITRO and MOSEK. SDPT3 did not perform well in these experiments, as we can see from the results in  Appendix \ref{app:performance}.

In the first experiment, following \cite*[Section 6.1]{PonteFampaLeeMPB}, we randomly generated normally-distributed elements for the  $n\times m$ full column-rank  matrices $A$, with mean $0$ and standard deviation $1$. For $m = 15,\dots,30$, we set $n := 10^{3}m$,  and $s := 2m$.

In the second experiment, we work with a subset of randomly-generated rows with respect to a ``full linear-response-surface model''. Generally, for a full linear model with 2 levels (coded as 0 and 1) and $F$ ``factors'', we have  $m=1+F$ and  $n=2^F$. Each row of $A$  has the form $v^\top:=(1;~\alpha^\top)$, with $\alpha\in\{0,1\}^F$. For our experiment, we set $i := 0,\dots,8$, and we define, for each $i$, $F := 19+i$, which leads to  $m = 20+i$. We set $n := (10 + 5i)\cdot 10^3$ (a subset of all possible rows) and $s:= 2m$. 

In the third experiment, following \cite*[Section 5.3]{hugedoptmohit_arxiv} (also see \cite{hugedoptmohit}), we work with a subset of randomly-generated rows with respect to a ``full quadratic-response-surface 
model''. In this case, for a full quadratic model with $L$ levels and $F$ ``factors'', we generally have $m=1+2F+\binom{F}{2}$ and   $n=L^F$. 
Each row of $A$  has the form
$
v^\top:=(1;~ 
\alpha_1,\ldots,\alpha_F\,;~
\alpha_1^2,\ldots,\alpha_F^2\,;~
\alpha_1 \alpha_2\,,\ldots ,\alpha_{F-1} \alpha_F)
$, and is identified by the levels in $\{0,1,\ldots,L-1\}$ of the factors $\alpha_1,\ldots,\alpha_F$\,.  For our experiment, we set $L:=3$ and $i := 0,\dots,8$.  For each $i$, we define $F := 19+i$ and we select $\binom{\left\lfloor{(F+1)/4}\right\rfloor}{2}$ pairs of factors (no squared term), which leads to $m = 1 + F + \binom{\left\lfloor{(F+1)/4}\right\rfloor}{2}$. We set $n := (10 + 5i)\cdot 10^3$ (a subset of all possible rows) and  $s:= 2m$. 
 
In the fourth experiment, we work with a real dataset, TICDATA2000.txt, which is the training data set that is part of the Insurance Company Benchmark (COIL 2000), from the
University of California Irvine (UCI) Machine Learning Repository; see \cite*{insurance}.
In our experiment, we worked with a  $5822 \times 60$ full column-rank matrix $A$
corresponding to the first 60 factors of that data set, and we set $s:=65,70,\dots,200$. 

\begin{figure}[!ht]
    \centering
    \subfigure[random instances] {
    \centering
    \includegraphics[width=0.47\linewidth]{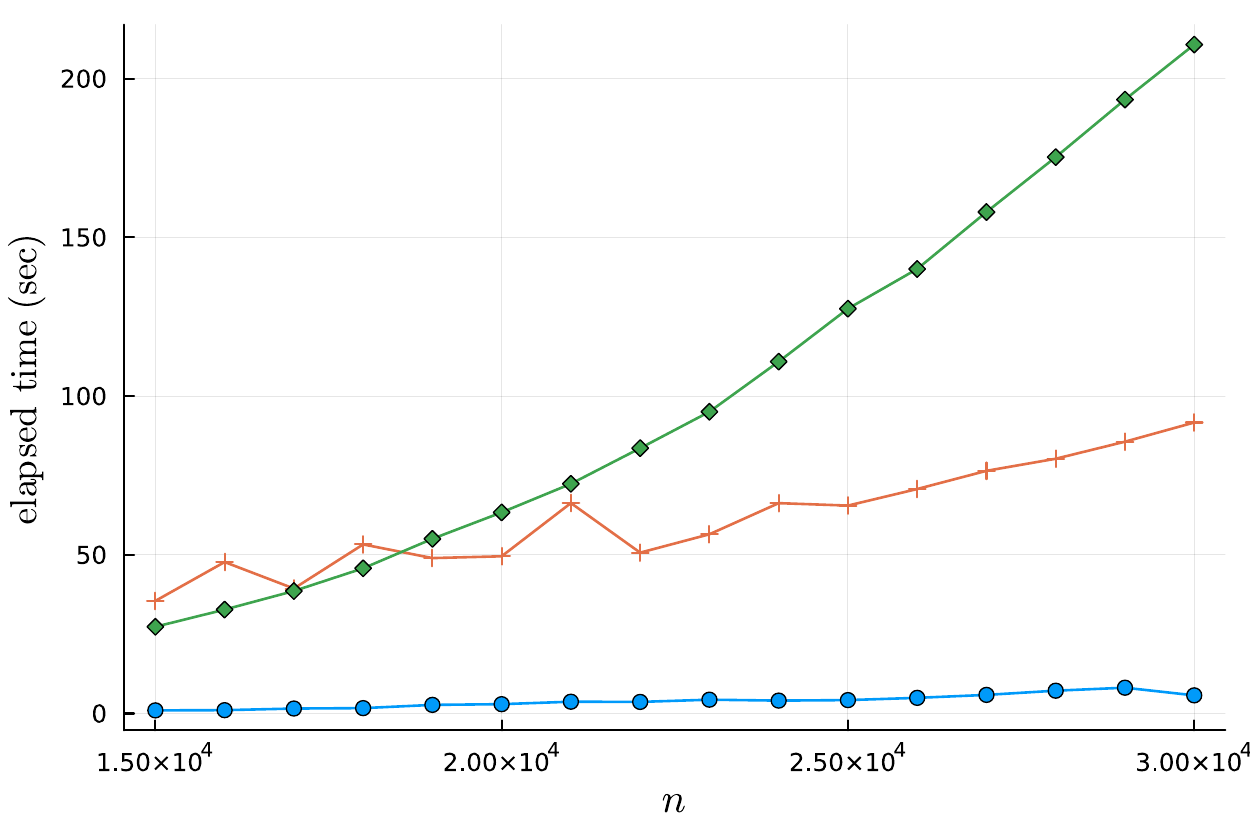}
    }
    \subfigure[linear-response model] {
    \centering
    \includegraphics[width=0.47\linewidth]{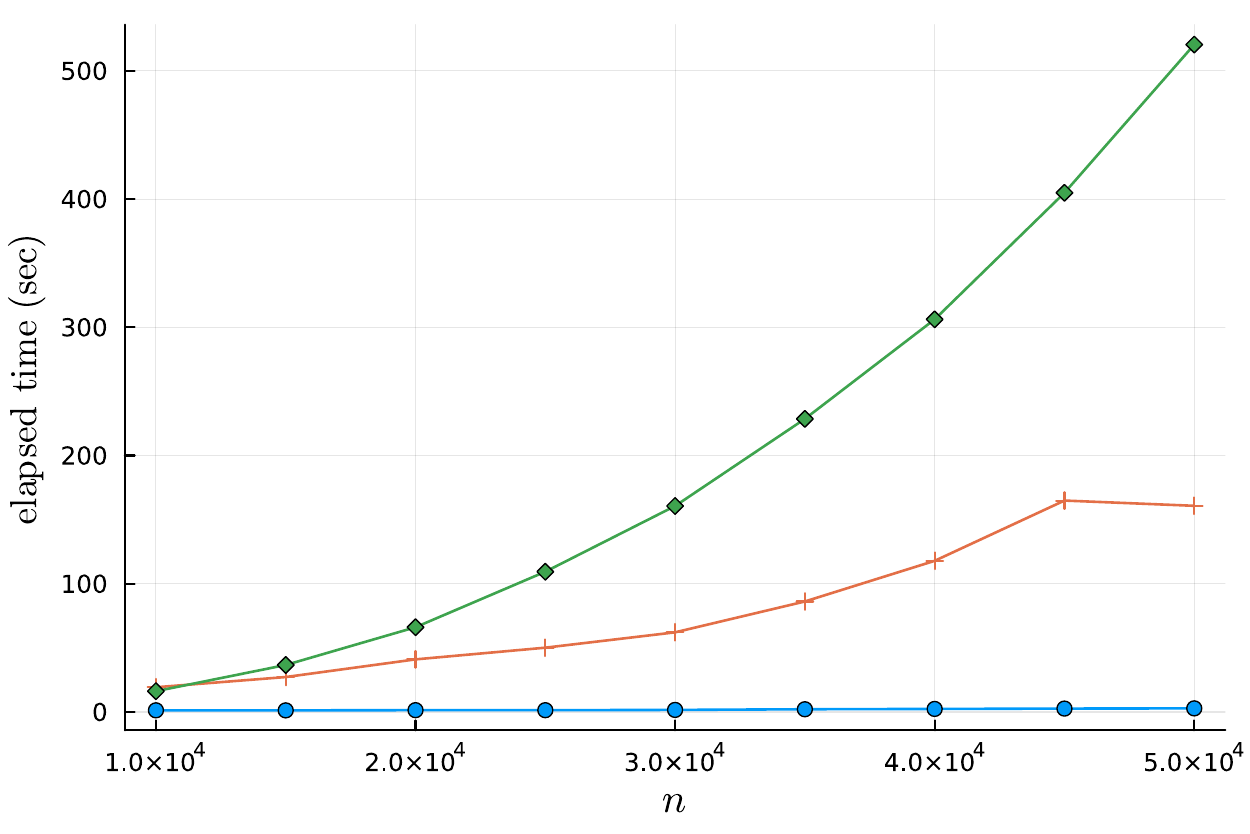}
    }
    \subfigure[quadratic-response model] {
    \centering
    \includegraphics[width=0.47\linewidth]{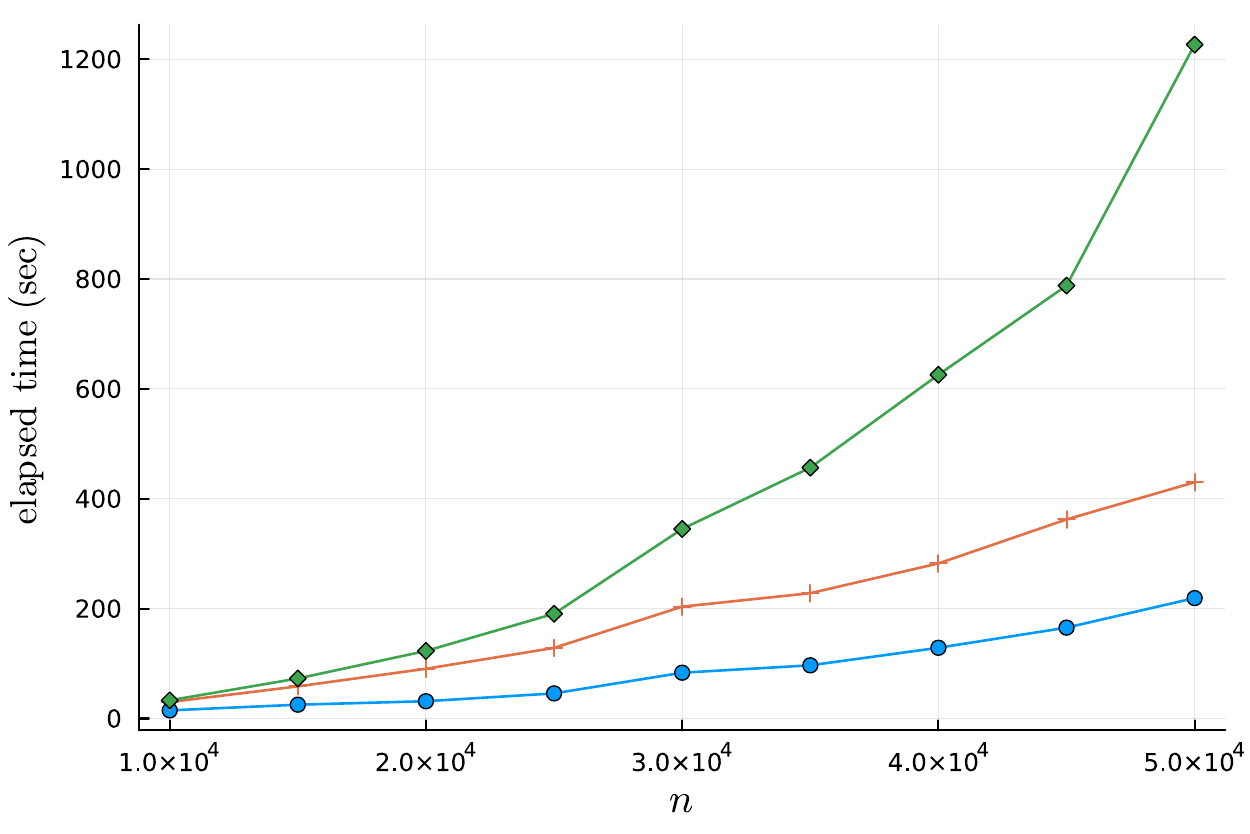}
    }
    \subfigure[real data set, $n=5822, m= 60$] {
    \centering
    \includegraphics[width=0.47\linewidth]{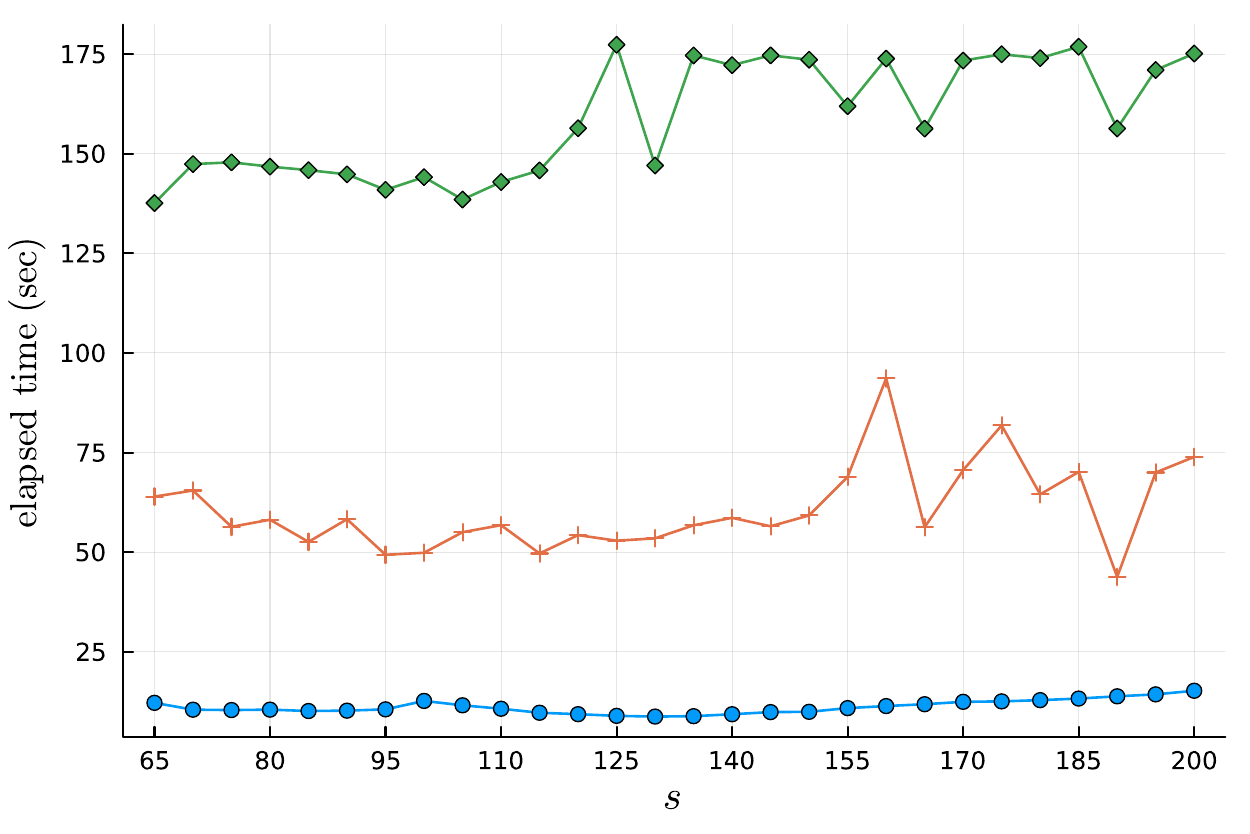}
    }\\
    \subfigure {
    \centering
    \includegraphics[width=0.525\linewidth]{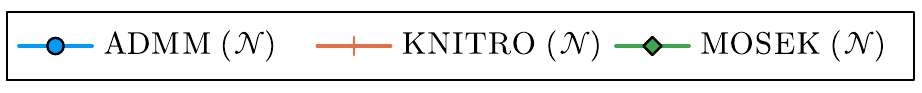}
    }
    \caption{Natural bound \ref{prob} for \ref{dopt}
    }
    \label{fig:dopt}
\end{figure}

In Figure \ref{fig:dopt}, we show the times to solve \ref{prob}, for the instances of the four experiments. We see that the ADMM algorithm for \ref{prob} performs very well in all of them, converging faster than KNITRO and MOSEK. We also observe that the times for the ADMM algorithm have a very stable behavior. Even for the quadratic-response model, where we see a larger increase in time with $n$, the increase is much smoother than for the solvers. 

\begin{figure}[!ht]
\setcounter{subfigure}{0}
    \centering
    \subfigure[random instances] {
    \centering
    \includegraphics[width=0.47\linewidth]{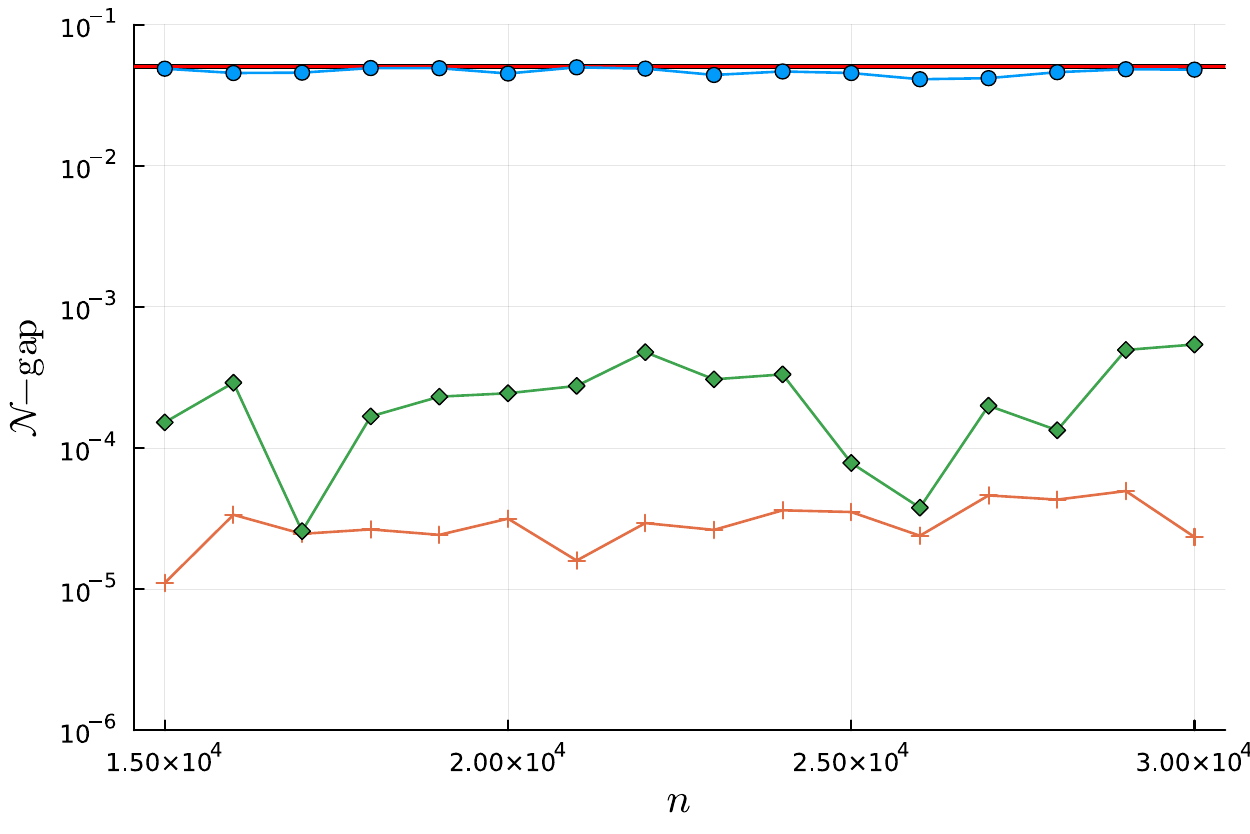}
    }
    \subfigure[linear-response model] {
    \centering
    \includegraphics[width=0.47\linewidth]{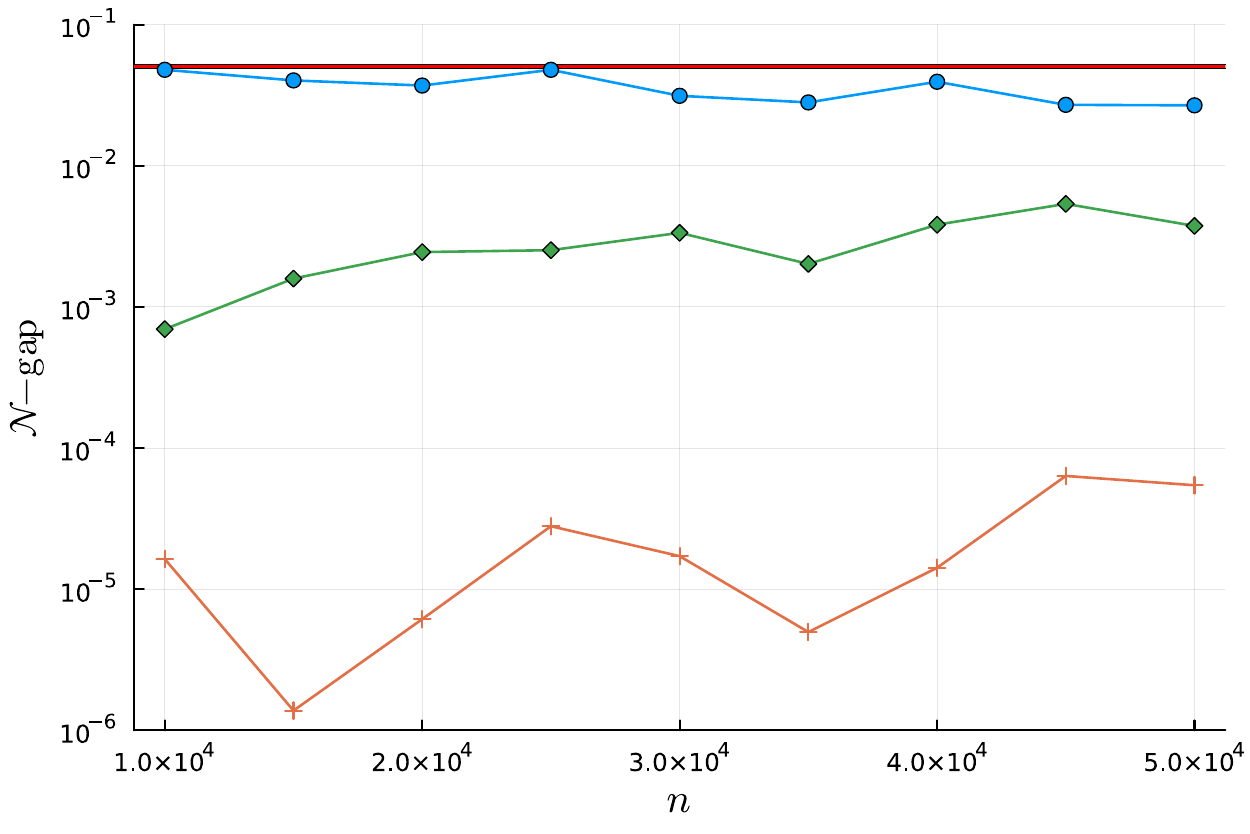}
    }

    \subfigure[quadratic-response model] {
    \centering
    \includegraphics[width=0.47\linewidth]{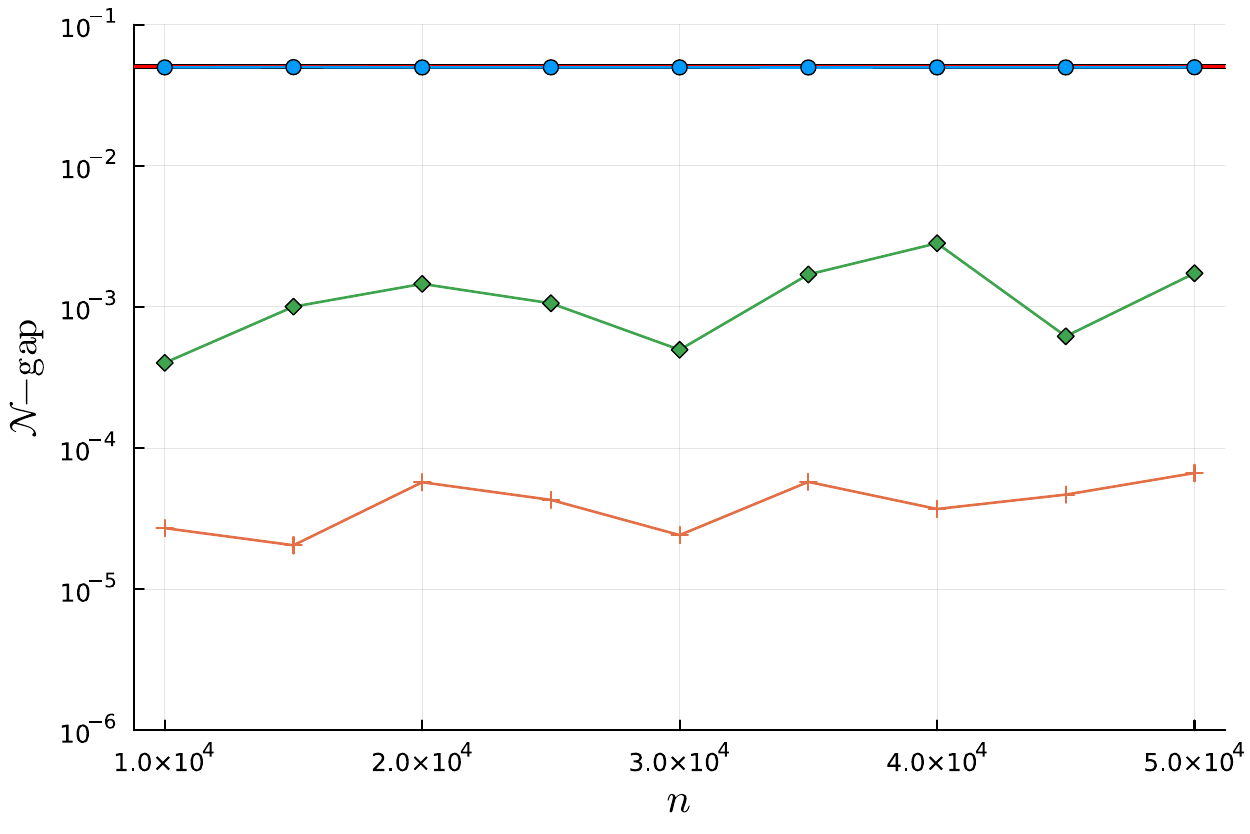}
    }
    \subfigure[real data set, $n=5822, m= 60$] {
    \centering
    \includegraphics[width=0.47\linewidth]{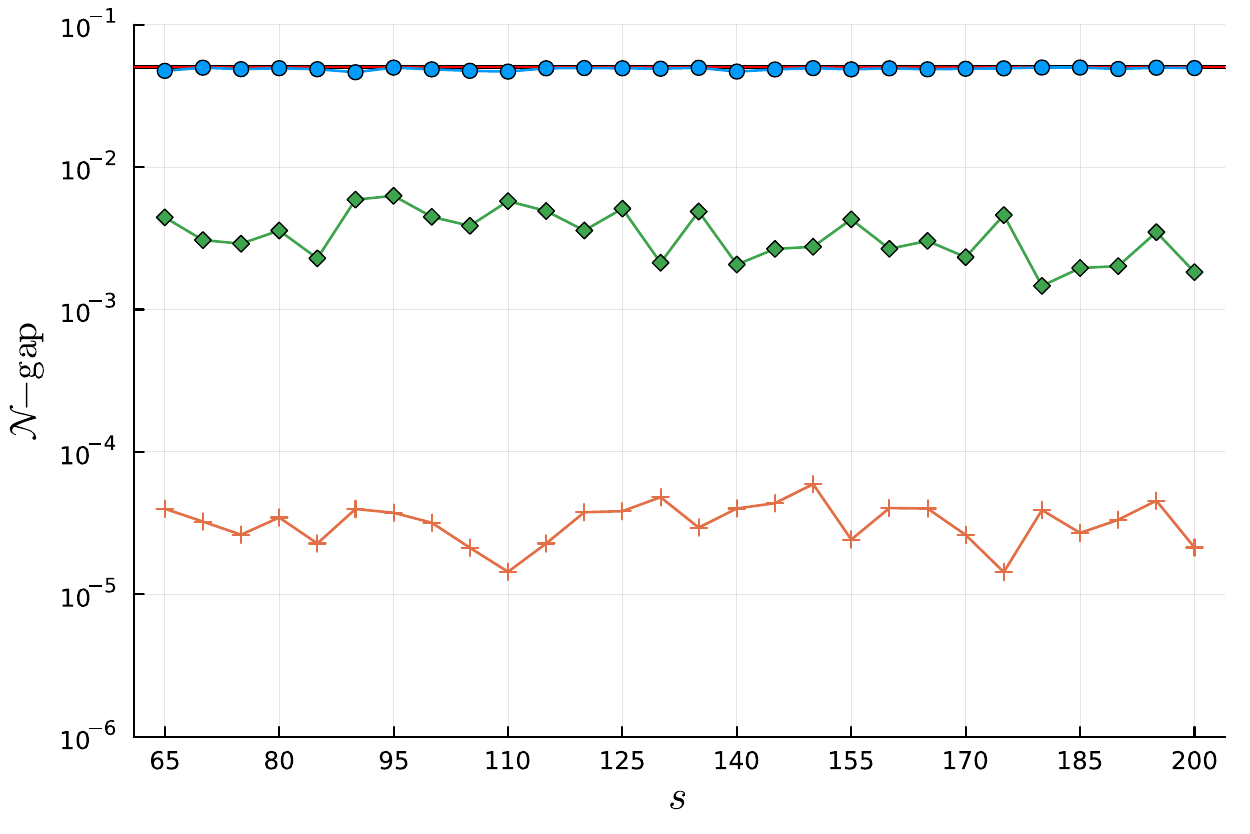}
    }\\
    \subfigure {
    \centering
    \includegraphics[width=0.70\linewidth]{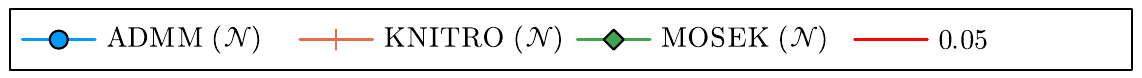}
    }\\
    \caption{Natural bound  \ref{prob} for \ref{dopt}}
    \label{fig:dopt_dual_gap}
\end{figure}

In Figure \ref{fig:dopt_dual_gap}, we show the dual gaps computed as previously described, from the solutions of the ADMM algorithm, MOSEK and KNITRO. We see that despite the parameter settings of the solvers seeking a $0.05$ optimality tolerance, the achieved differences between the dual and primal solution values are smaller. 
It is not surprising that the general-purposes solvers have this behavior, as they are aimed at constrained optimization, where a significant effort can be devoted to obtaining primal and dual feasibility, and once that is achieved, the gaps can turn out to be small.
Finally, we can see in Figure \ref{fig:dopt_opt_gap} in Appendix \ref{app:performance}, that $0.05$ is not significant when compared to the differences between the upper bounds and the best known solution values for the instances considered.

In Tables \ref{tab:random}--\ref{tab:real} of  Appendix \ref{app:performance}, we give the results that form the basis for Figures \ref{fig:dopt}--\ref{fig:dopt_dual_gap}, the $\rho$ values used for our ADMM, 
as well as (worse) results for additional solvers.

\subsection{Experiments for the Maximum-Entropy Sampling Problem}
We conducted experiments for the ADMM algorithm in \S\ref{subsec:ddfact}, to compute the factorization bound, and for the ADMM algorithm  in \S\ref{subsec:BQP}, to compute the BQP bound. We do not show results for the ADMM algorithm for \ref{prob_linx}\,. As we noted earlier, the bottleneck for the algorithm is the solution of the subproblem \eqref{eq:Zminlinxsubprob}, which makes our ADMM algorithm for \ref{prob_linx} not competitive.
Nevertheless, we decided to present the algorithm on \S\ref{subsec:linx}, in the hope that we can speed up the solution of the subproblem in future work.


\subsubsection{ADMM for the factorization bound}

We discuss two experiments to test the  ADMM algorithm  described in \S\ref{subsec:ddfact}, to compute the factorization bound from \ref{prob_ddfact} for \ref{MESP}.
For these experiments, we considered an $n=2000$ covariance matrix  with rank 949 based 
on Reddit data from  \cite*{Dey2018} and \cite*{Munmun}, and also used by 
\cite*{Weijun} and \cite*{ChenFampaLee_Fact}.

Before presenting our results, some observations should be made. We  first note that, for all instances tested, the inequality \eqref{conditionj} always holds, and therefore the integer $\hat\jmath$ considered in Lemma \ref{lem:theta_j} exists. Thus, we can  successfully solve  subproblems \eqref{eq:Zminddfactsubprob} with the closed-form solution presented in Theorem \ref{thm:updateGamma}. Nevertheless, if this were not the case, we could use an iterative algorithm to solve the subproblem for which $\hat\jmath$ could not be computed, for example, from KNITRO. 
Furthermore, we note that Proposition \ref{prop:weijun_supgrad} is defined for $Z \in \mathbb{S}_{+}$\,, and from Lemma \ref{lem:lambda_j} we may have $\lambda \not\geq 0$. In this case, we could project $\lambda$ onto the nonnegative orthant and then apply Theorem \ref{thm:updateGamma} to construct $Z^{t+1}$. 
However, in practice, when $\lambda$ has negative components (which are often quite small), we continue to construct $Z^{t+1}$ by applying Theorem \ref{thm:updateGamma}. This approach  worked better than projecting  $\lambda$ onto the nonnegative orthant and it did not impact the practical convergence of the ADMM algorithm.

\begin{figure}[hbtp]
	\centering
	\includegraphics[width=0.48\linewidth]{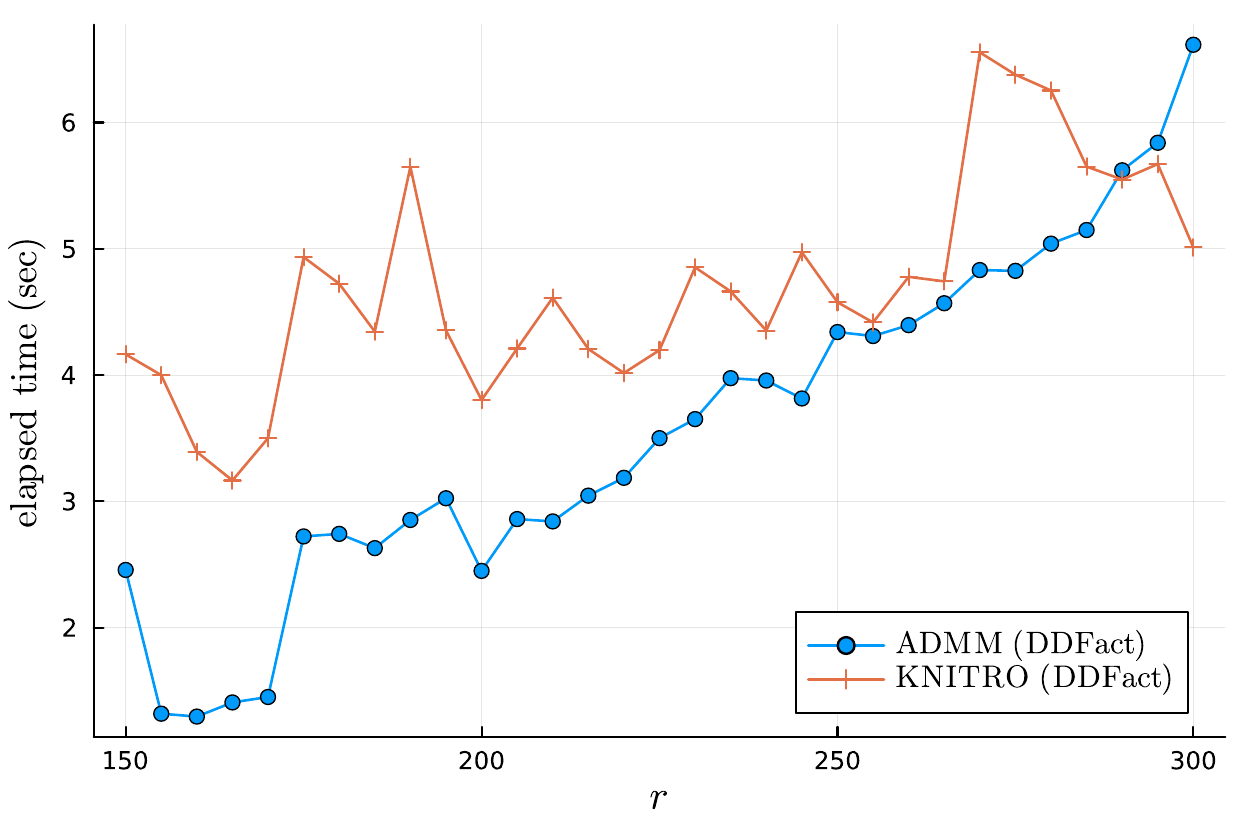}
	\includegraphics[width=0.48\linewidth]{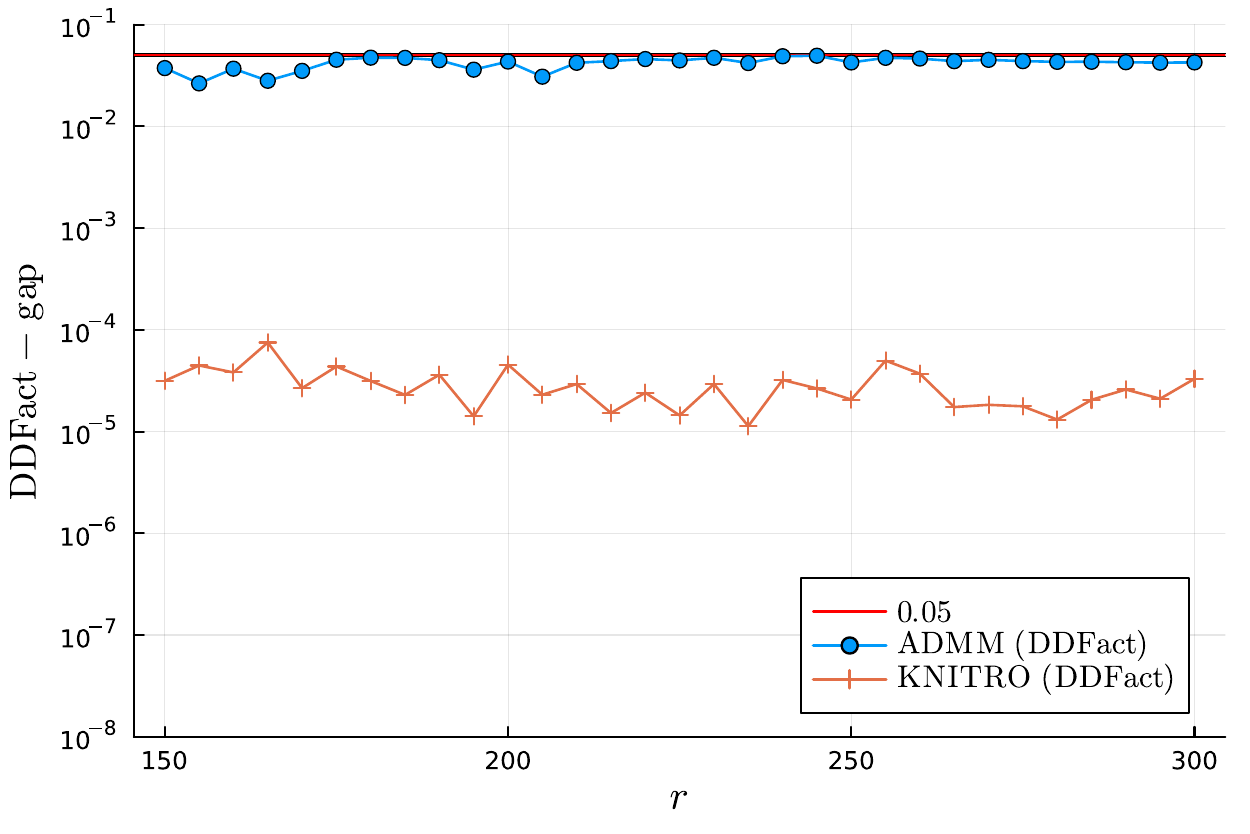}\\
	\caption{\ref{prob_ddfact} bound for \ref{MESP}, varying $r:=\rank(C)$ ($n=2000$, $s=140$)
    }\label{fig_LS_bound_r}
\end{figure}

In our first experiment, to analyze the performance impact of the rank of $C$, we constructed matrices with rank $r:=150,155,\ldots,300$, derived from the benchmark $n=2000$ covariance matrix  by selecting its $r$-largest principal components. For all $r$, we set $s:=140$. The results are in Figure \ref{fig_LS_bound_r}. In the first plot, we have the times for our ADMM algorithm and for KNITRO to solve \ref{prob_ddfact}.
We see that the ADMM algorithm is very efficient for \ref{prob_ddfact}. The vast majority of instances could be solved
faster than when KNITRO is applied. We can see that the ADMM algorithm takes advantage of the fact that the eigenvector decomposition required to update $Z$ (described in \S\ref{subsec:update_Zddfact}) is computed over a matrix of order $r:=\mbox{rank}(C)$, which is more efficient for smaller ranks. When the rank increases, this computation, which is a bottleneck of the ADMM algorithm, becomes heavier.

In the second plot of Figure \ref{fig_LS_bound_r}, we show the dual gaps computed from the solutions of the ADMM algorithm and KNITRO. As in Figure \ref{fig:dopt_dual_gap}, we see that although the KNITRO parameter settings seek an optimality tolerance of $0.05$, the differences achieved between the values of the dual and primal solution are smaller. The comment on Figure \ref{fig:dopt_dual_gap} could be repeated here.

\begin{figure}[hbtp]
	\centering
	\includegraphics[width=0.48\linewidth]{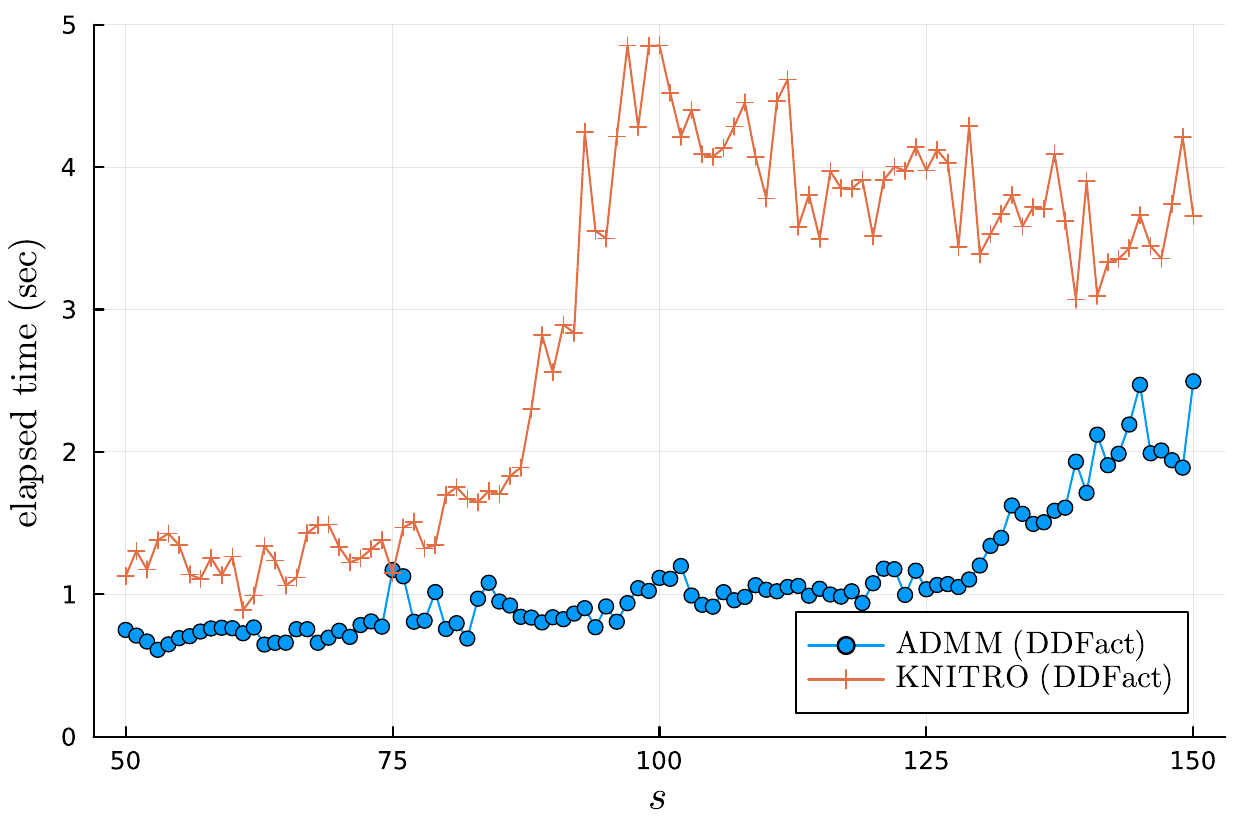}
    \includegraphics[width=0.48\linewidth]{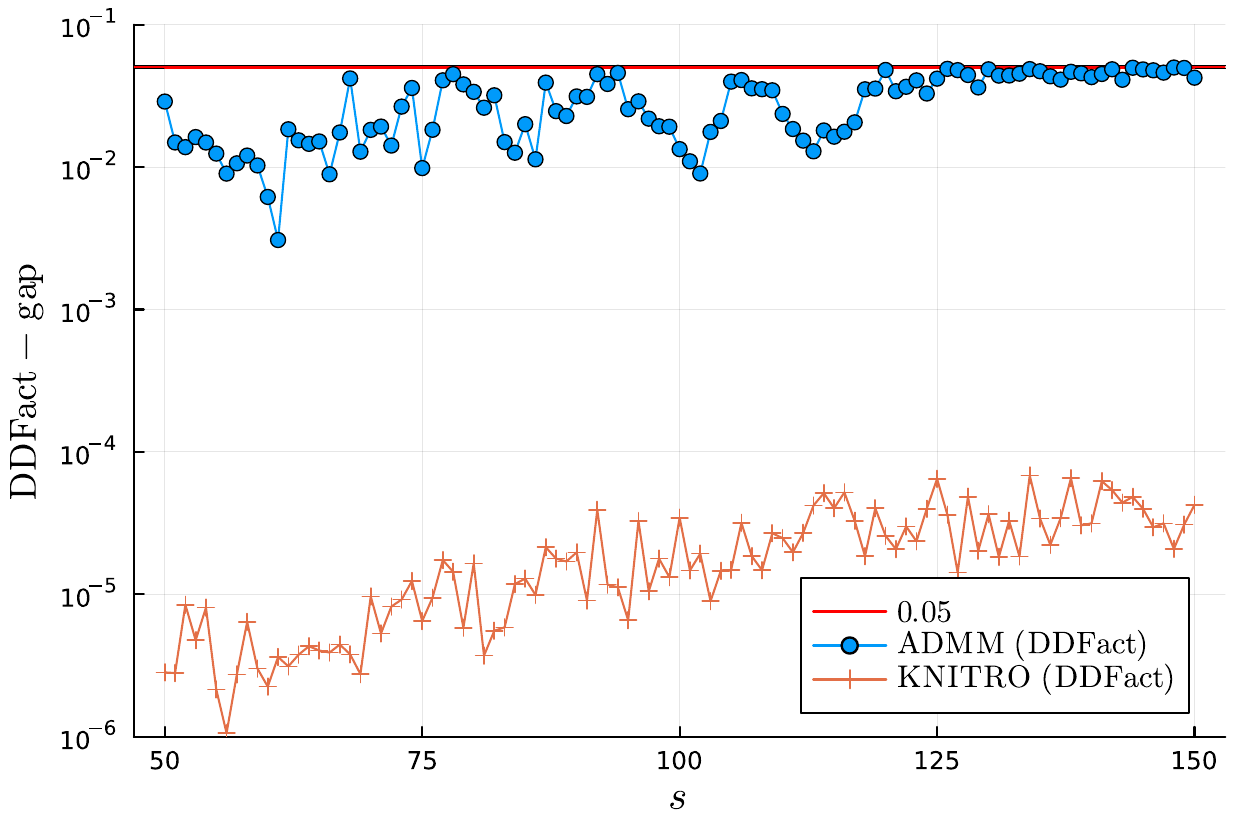}\\  
	
	\caption{\ref{prob_ddfact} bound for \ref{MESP}, varying $s$ ($n=2000$, $\rank(C)=150$)}\label{fig_LS_bound_s}
\end{figure}

In the second experiment, our aim is to analyze the impact of $s$ on the performance of the ADMM algorithm. In this case, we fix  $r:= 150$, i.e., we consider a matrix derived from the benchmark $n=2000$ covariance matrix  by selecting its $150$-largest principal components, and we set $s:=50,51,\dots,150$. In Figure \ref{fig_LS_bound_s}, we show results similar to those presented in Figure \ref{fig_LS_bound_r}, but now varying $s$ instead of $r$. Unlike what we see in Figure \ref{fig_LS_bound_r}, we now see a less significant impact of the increase in $s$ on the performance of the ADMM algorithm. It performs very well, with faster convergence than KNITRO for all instances. 

We conclude that, in general,  the ADMM algorithm is a very good method to compute the \ref{prob_ddfact} bound when the covariance matrix has a low rank. 

We observe in Figures \ref{fig_LS_bound_r} and \ref{fig_LS_bound_s}, varying both $r$ and $s$, that we generally have more stable computation times for our ADMM algorithm than for KNITRO, as we have observed in Figure \ref{fig:dopt} as well. 

Finally, we refer to  Figure \ref{fig_LS_bound_app_alone}  in  Appendix \ref{app:performance}, to confirm that $0.05$ is not significant when compared to the differences between the upper bounds and the best known solution values for the instances considered in the two experiments described above. It is also interesting to note from Figure \ref{fig_LS_bound_app} that, for the instances considered, \ref{prob_ddfact} gives a better bound than \ref{prob_linx}\,; additionally, we can report that the \ref{bqp_original} bound cannot be computed 
within the time limit for these instances, using any algorithm or software that we have tested.

In Tables \ref{tab:ddfact_r}-\ref{tab:ddfact_s2} of Appendix \ref{app:performance}, we give the results that form the basis for Figures \ref{fig_LS_bound_r}--\ref{fig_LS_bound_s}, the $\rho$ values used for our ADMM, 
as well as (worse) results for an additional solver.


\subsubsection{ADMM for the BQP bound}

We discuss two experiments to test the  ADMM algorithm  described in \S\ref{subsec:BQP}, to compute the BQP bound from \ref{bqp_original} for \ref{MESP}. We note that for the computation of the bounds, we first optimize the scaling parameter $\gamma$ (see \cite*{Mixing} regarding optimizing
the choice of $\gamma$).
Moreover, for the nonsingular benchmark   covariance matrix $C$ used in the experiments,  we compute the bounds considering the original relaxation and the complementary relaxation, and present only the results corresponding to the best. 

We compare the results for our ADMM algorithm with SDPT3, which performed better than MOSEK on this problem.

In our first experiment, we  use a benchmark covariance matrix of dimension $n=63$,  originally  
obtained from J. Zidek (University of British Columbia), coming from an application for re-designing an environmental monitoring network;
see \cite*{Guttorp-Le-Sampson-Zidek1993} and \cite*{HLW}. This  matrix has been used extensively in testing and developing algorithms for \ref{MESP}; see \cite*{KLQ,LeeConstrained,AFLW_Using,LeeWilliamsILP,HLW,AnstreicherLee_Masked,BurerLee,Anstreicher_BQP_entropy,Kurt_linx,Mixing,ChenFampaLee_Fact}. 

\begin{figure}[hbtp]
	\centering
	\includegraphics[width=0.48\linewidth]{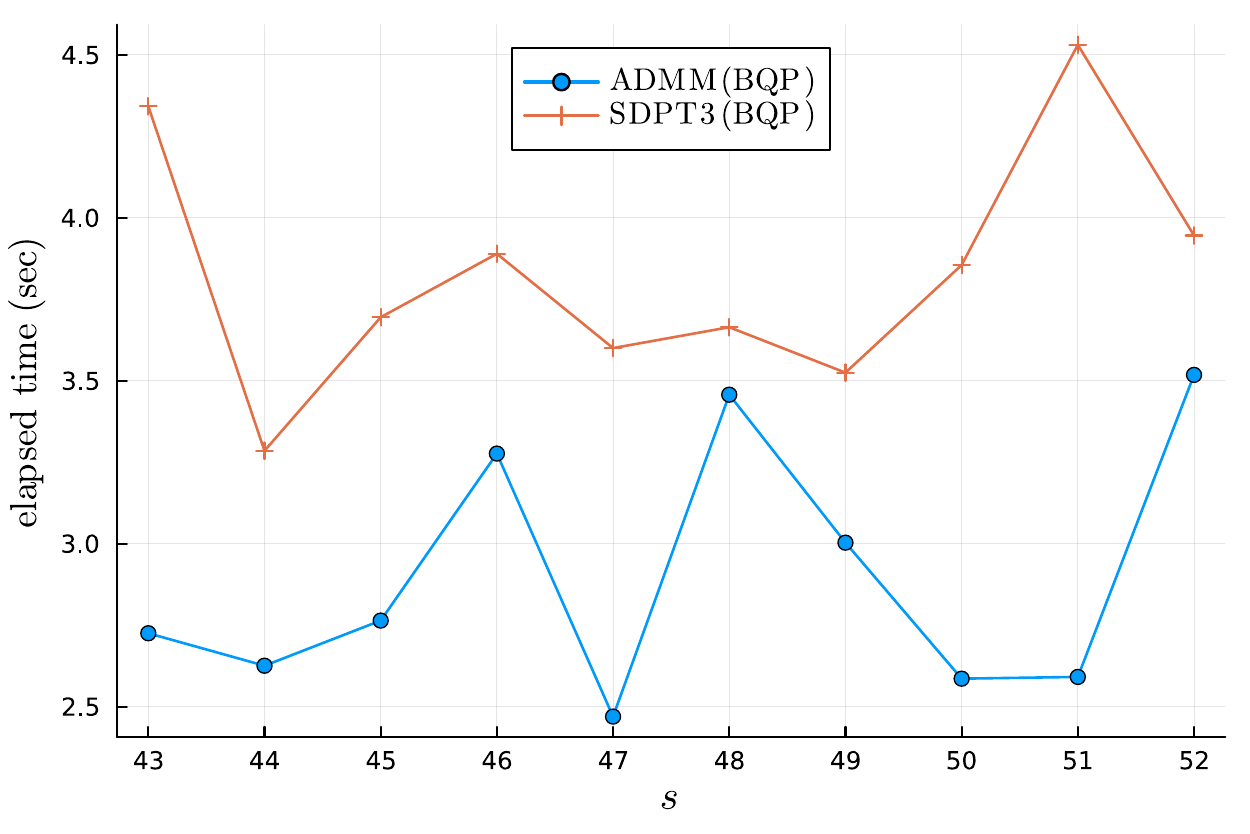}
    \includegraphics[width=0.48\linewidth]{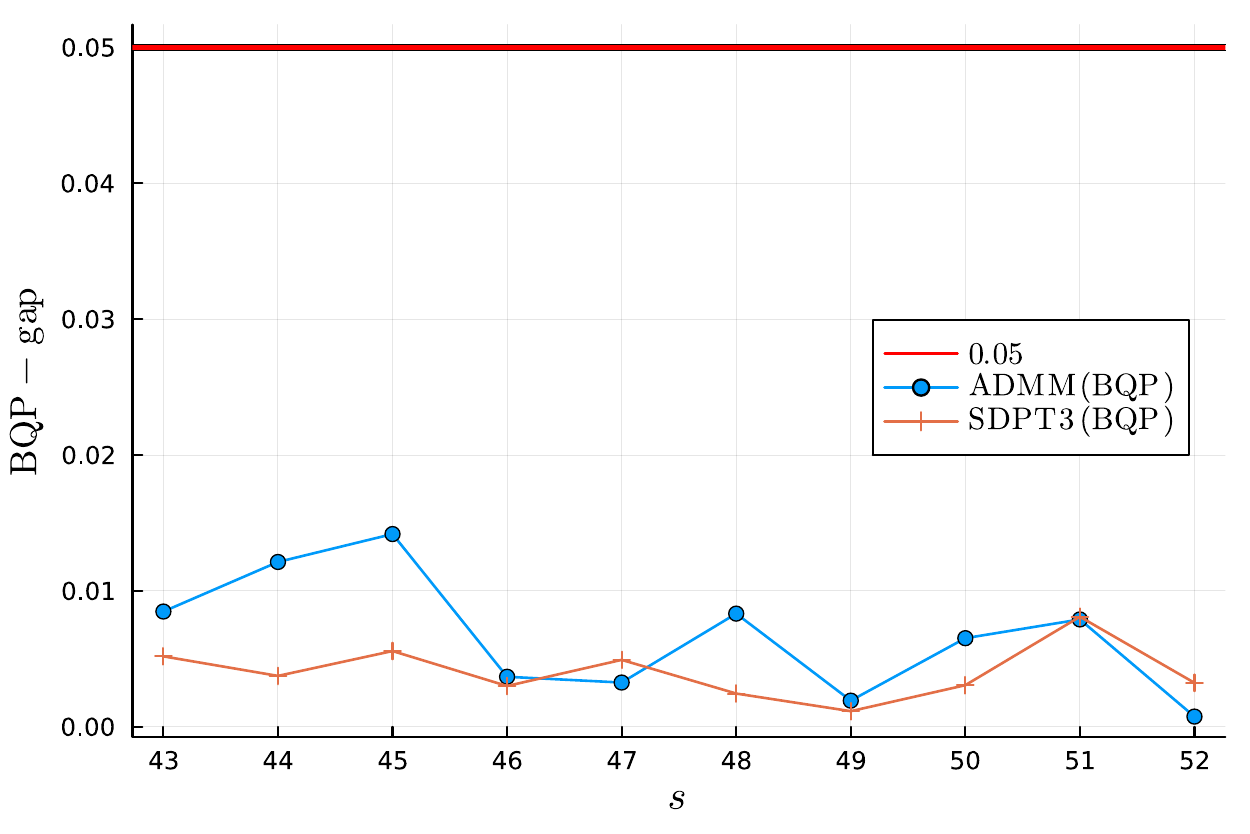}\\
	\includegraphics[width=0.48\linewidth]{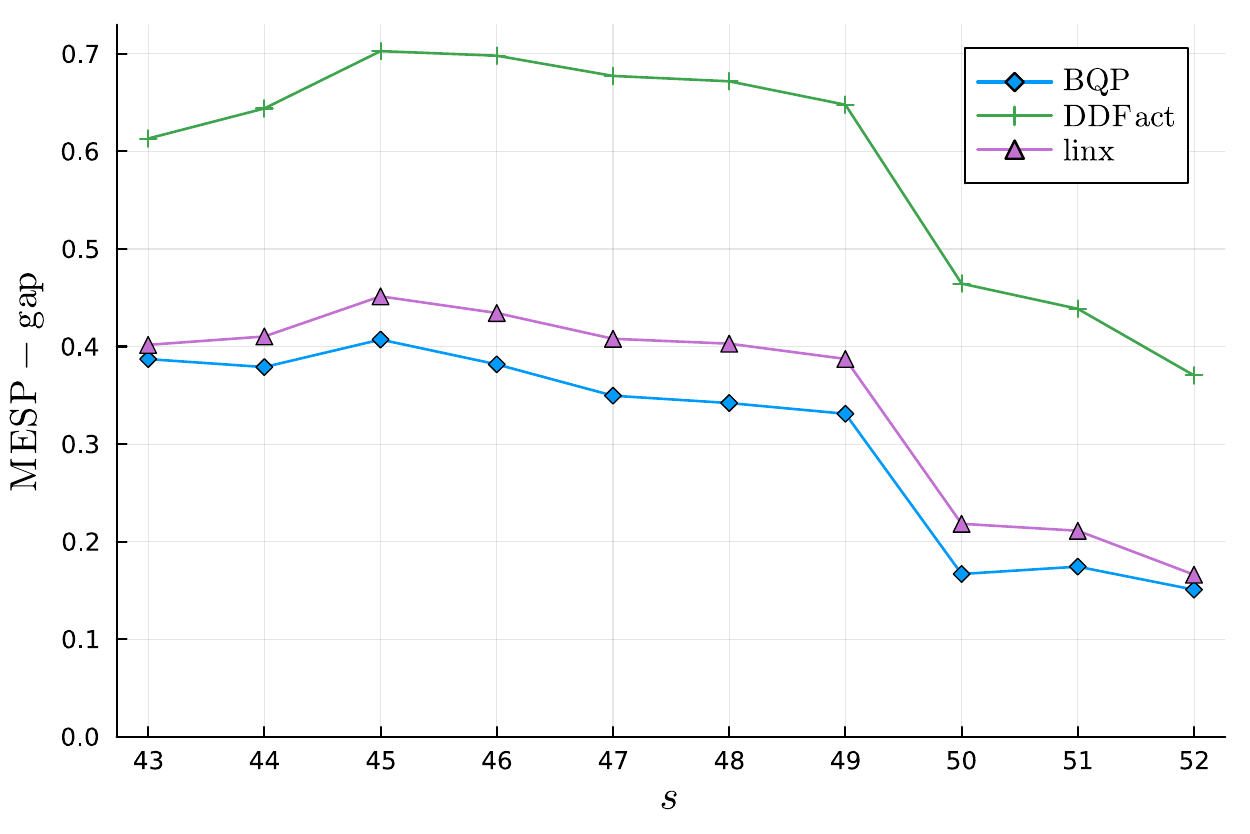}
	
	\caption{Behavior of the \ref{bqp_original} bound for \ref{MESP}, varying $s$ ($n=63$)}\label{fig_bqp_bound63}
\end{figure}

In Figure \ref{fig_bqp_bound63} we show results for $s=43,\ldots, 52$. We intentionally selected these values of $s$ to consider instances for which \ref{bqp_original} gives a better bound than \ref{prob_ddfact} and \ref{prob_linx}\,, motivating its consideration. 
In the first plot in Figure \ref{fig_bqp_bound63}, we  show the  times  to solve \ref{bqp_original}\,. We see that the ADMM algorithm for \ref{bqp_original} performs very well, converging faster than SDPT3 in all instances. In the second plot, we show the dual gaps computed as previously described, from the solutions of the ADMM algorithm and SDPT3. We see that, the dual gaps are smaller than the optimality tolerance of $0.05$. We saw this same behavior in Figure \ref{fig:dopt_dual_gap} for the solvers, but here we also see it for the ADMM algorithm. 
The reason is that, as mentioned above, due to the cost of computing dual solutions for \ref{bqp_original}\,, we only start computing them after 
many
iterations,
and for the considered instances, the dual gap was already smaller than $0.05$ at this point, leading to the advantage of needing to compute the dual solution only once.
Finally, on the third plot of Figure \ref{fig_bqp_bound63}, we see that the bound from \ref{bqp_original}  is a competitive bound for \ref{MESP}, which motivated the development of the ADMM algorithm to allow its computation for larger instances than the solvers can handle as we will address in the next experiment. 

In the second experiment, we use full-rank  principal submatrices of an order-$2000$ covariance matrix with rank $949$, based on Reddit data, used in \citep*{Weijun} and from  \citep*{Dey2018} (also see \citep*{Munmun}).  The submatrices selected have dimensions  $n=250,275,\ldots, 400$, and we set  $s := \lfloor n/2\rfloor$ in all test instances. To select the linear independent rows/columns of the order-$2000$ matrix,  we use the Matlab function nsub\footnote{\url{www.mathworks.com/matlabcentral/fileexchange/83638-linear-independent-rows-and-columns-generator}}  (see \cite*{FLPX2021} for details).

\begin{figure}[hbtp]
	\centering
	\includegraphics[width=0.48\linewidth]{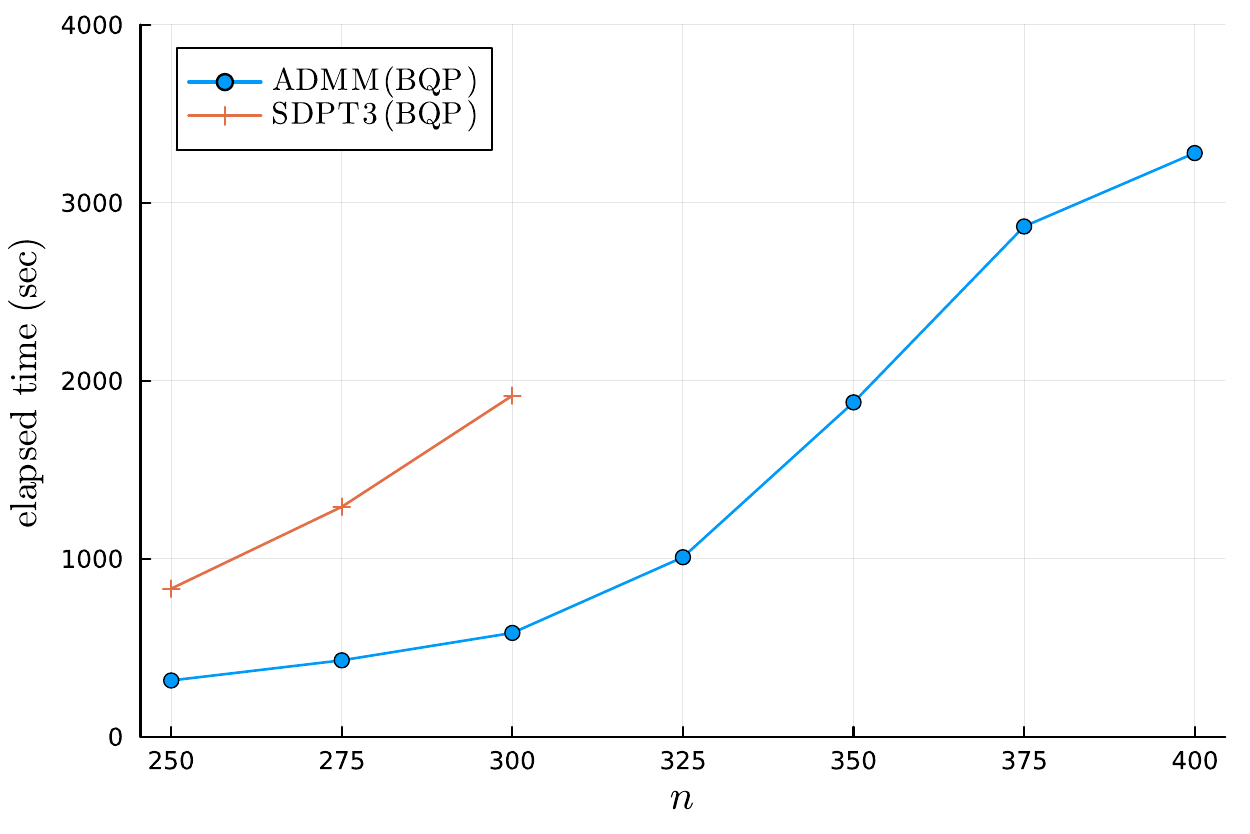}
	\includegraphics[width=0.48\linewidth]{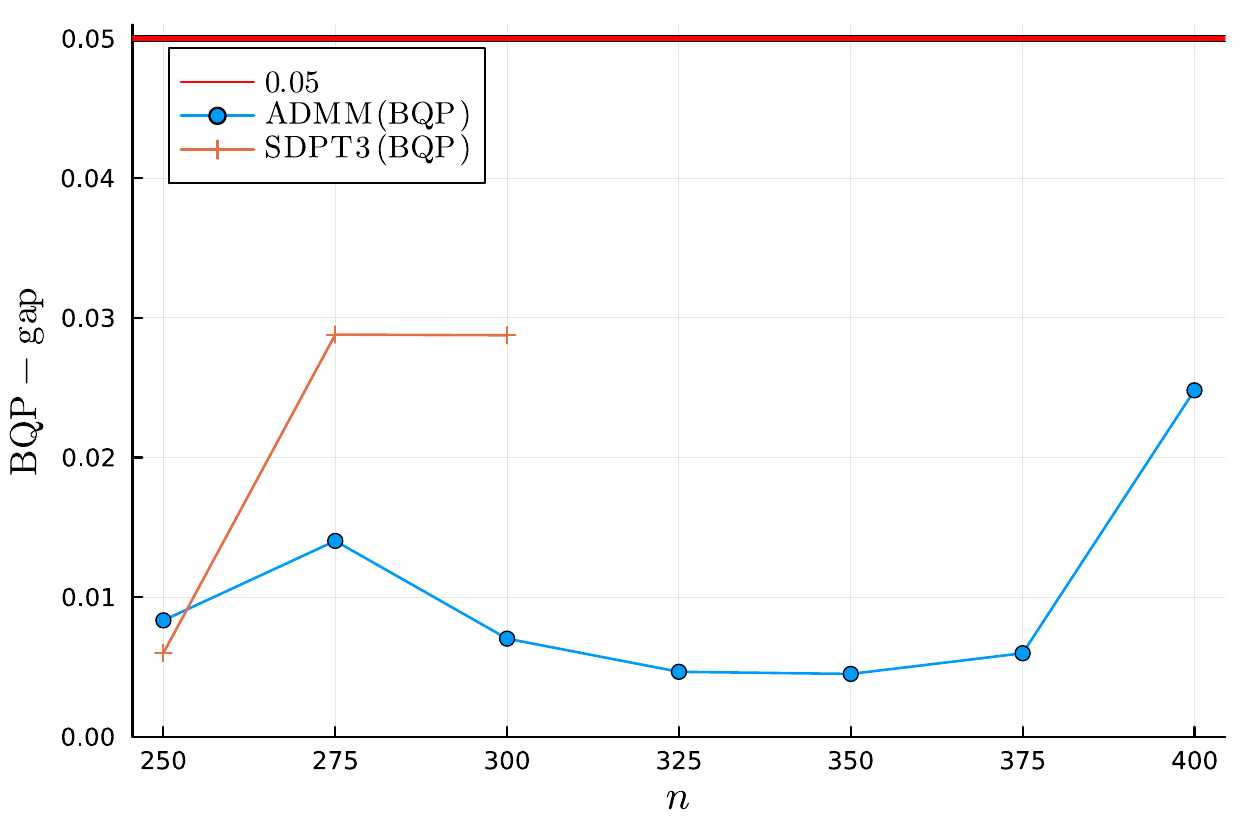}\\
    \includegraphics[width=0.48\linewidth]{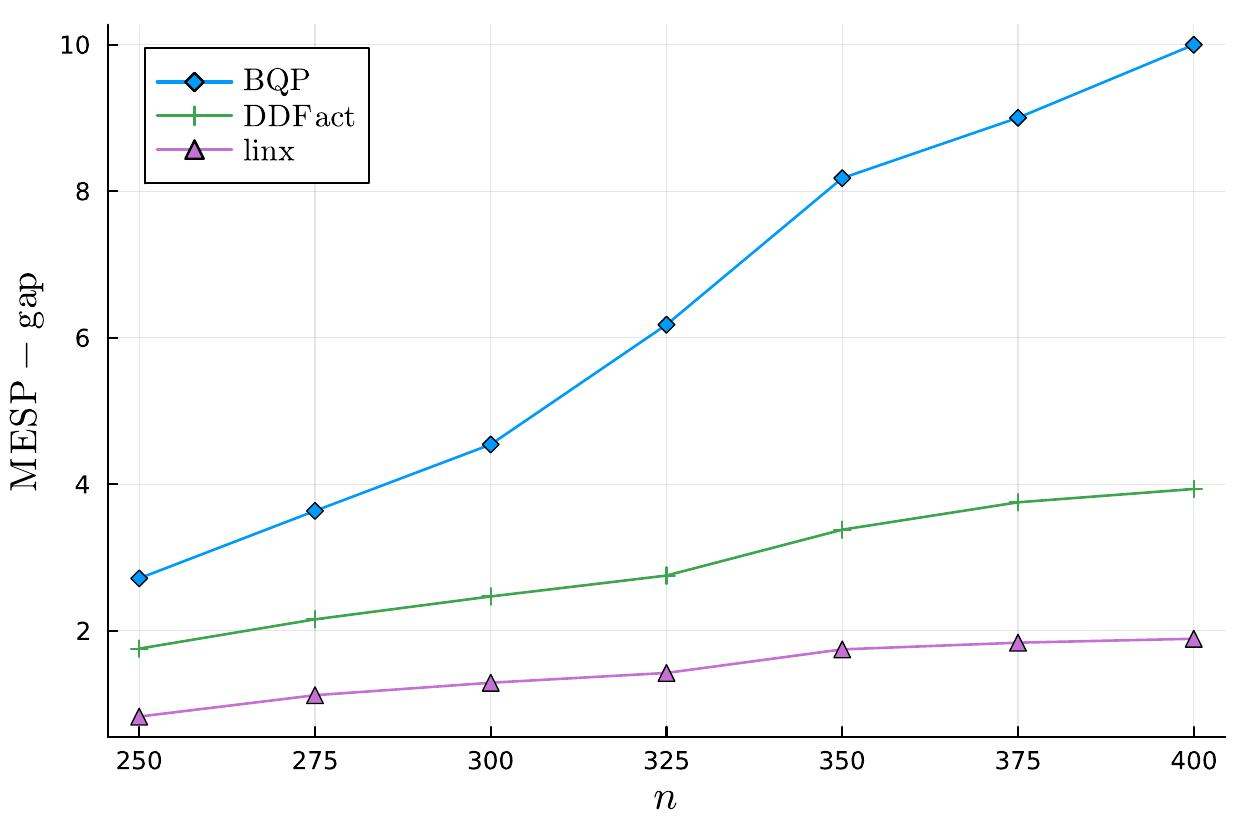}
	\caption{Behavior of the \ref{bqp_original} bound for \ref{MESP}, varying $n$, with $s := \lfloor n/2\rfloor$}\label{fig_bqp_bound2000}
\end{figure}

In Figure \ref{fig_bqp_bound2000}, we show the same statistics as shown in Figure \ref{fig_bqp_bound63} for this second experiment. Although the linx and factorization bounds are better than the BQP bound for these instances (see the third plot in Figure \ref{fig_bqp_bound2000}), it is still interesting to be able to solve \ref{bqp_original} and thus be able to investigate the BQP bound for them. We see that SDPT3 crashed due to lack of memory when $n>300$.  We note that we also tried to solve \ref{bqp_original} with MOSEK, but it crashed already for $n=250$ due to lack of memory.   

In Tables \ref{tab:bqpn63}--\ref{tab:bqp_n} of Appendix \ref{app:performance}, we give results that form most of the basis for Figures \ref{fig_bqp_bound63}--\ref{fig_bqp_bound2000}, the $\rho$ values used for our ADMM, 
as well as (worse) results for additional solvers.

\subsection{Impact of warm-starting ADMM within branch-and-bound}

To illustrate the potential impact of warm-starting ADMM within a B\&B framework, we consider the TICDATA200 instance of the D-optimality problem with $s:=150$ and perform a ``diving'' procedure based on the most fractional variable. Starting from the root node of the B\&B tree, at each diving step we compute the natural bound for the corresponding subproblem using ADMM, obtaining a solution $x^*$ with duality gap below $0.05$. We then select the most fractional component $x_j^*$\,, i.e., the one closest to $0.5$, and fix its value as it follows: if $x_j^*\leq 0.5$, we set $x_j:=0$, otherwise, we set $x_j:=1$.
Our diving procedure, aimed at simulating what can happen within B\&B, is
inspired by a successful ``diving heuristic'' for MINLO (see for example, \cite*{diving}).

The resulting subproblem, which incorporates all previous fixings together with the new one, is warm-started as follows. 
The initial dual variables $\Psi^0$ and $\delta^0$ are inherited directly from the parent node. To obtain $Z^0$, we consider the projection $\bar{x}$ of the parent node solution onto the set 
\[
\mathcal{X}:=\{x\in\mathbb{R}^n: 0\leq x\leq \mathbf{e}; ~x_i=0, ~\forall i\in\mathcal{I}^0; ~x_i=1, ~\forall i\in\mathcal{I}^1\}, 
\]
where $\mathcal{I}^0\subset N$ (resp., $\mathcal{I}^1\subset N$) denotes the set of indices of variables fixed at $0$ (resp., $1$) in the subproblem under consideration.
Given $\bar{x}$ and $\Psi^0$, $Z^0$ is then recomputed using the closed-form update (see Proposition~\ref{lem:closedformulaupdateZ}). 

This procedure is repeated $100$ times, simulating successive diving steps along the B\&B tree.
The ADMM scheme used to compute the natural bound for each subproblem follows the method described in Section~\ref{sec:d-opt}, with a single modification: in the $x$ update, instead of projecting onto $[0,1]^n$, we project onto $\mathcal{X}$.

\begin{figure}[!ht]
    \centering
    \subfigure 
    {
    \centering
    \includegraphics[width=0.47\linewidth]{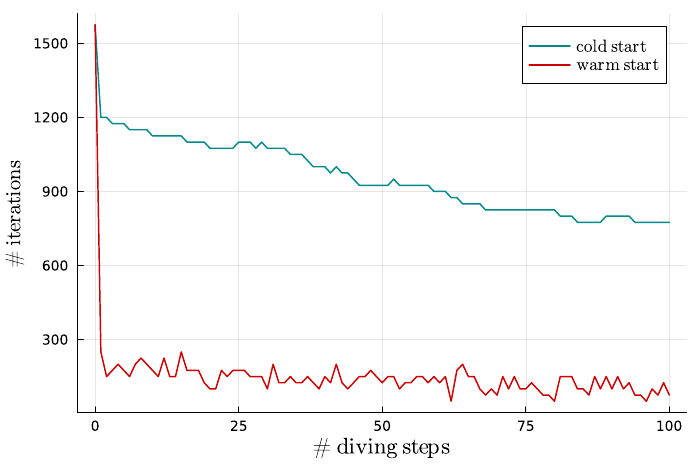}
    }
    \subfigure 
    {
    \centering
    \includegraphics[width=0.47\linewidth]{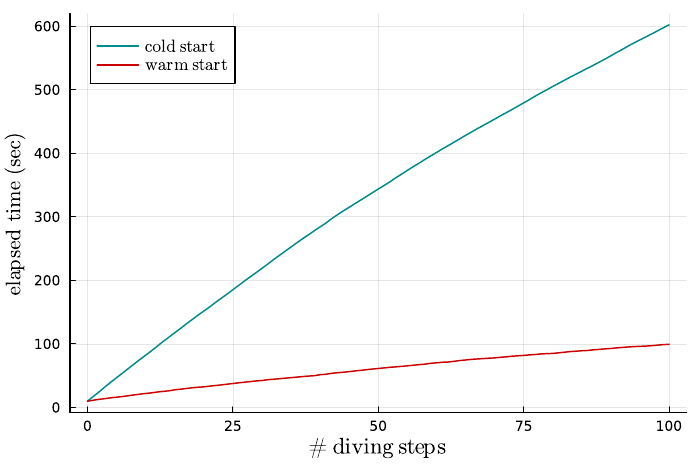}
    }
    \caption{Warm-start performance for the natural bound \ref{prob} for \ref{dopt} ($n=5822$, $m=60$, $s:=150$)
    }
    \label{fig:warm-start}
\end{figure}

We compare the performance of warm-started ADMM with a cold-start strategy, in which each subproblem arising in the diving procedure is initialized as described in Subsection \ref{subsec:implementation}. The results are reported in Figure \ref{fig:warm-start}. Warm-starting yields a substantial reduction in the number of ADMM iterations. While the root node requires $1575$ iterations for both approaches, subsequent warm-started subproblems consistently converge in $50$ to $250$ iterations. In contrast, the cold-start approach requires between $775$ and $1200$ iterations, with a gradual decrease as more variables are fixed and the subproblems become easier. Notably, the warm-started method quickly stabilizes at a low iteration count, whereas the cold-start method improves only progressively.

This difference is also reflected in the total computational time. After $100$ diving steps, the cumulative runtime is approximately $600$ seconds for the cold-start approach, compared to $100$ seconds for the warm-start strategy, corresponding to an overall speedup of about a factor of $6$.

\section{Next steps}\label{sec:concl}

Besides the bounds that we have considered,
there is also an effective (so-called) ``NLP bound'' for \ref{MESP}
(see \cite*{AFLW_Using} and \cite*[Section 3.5]{FL2022}).
But our ADMM approach 
would unfortunately lead to a non-convex subproblem, because for that bound, $\ldet$ acts on a nonlinear function of the problem variable $x\in\mathbb{R}^n$.
So we leave it as a challenge to develop a fast first-order method for calculating the NLP bound.
 We would like to mention, that for a particular well-known parameter choice,
the NLP bound is the so-called ``NLP-Id bound'',
and for that, using results of \cite*{MESP2DOPT}, we can calculate the NLP-Id bound for a given \ref{MESP} instance using the natural bound for a related 
\ref{dopt} instance. So, combining that 
with the ADMM for the natural bound in the present work,
we can calculate the NLP-Id bound for \ref{MESP} using an ADMM algorithm.

Our work develops tools that can be incorporated in  B\&B algorithms for \ref{dopt} and \ref{MESP}.
In that context, convex relaxations need to be solved to modest accuracy, and if we can re-solve quickly  based on ``parent'' solutions, then we have the possibility to handle a very large number of 
 B\&B subproblems. 
 We believe that our ADMM algorithms are very well suited for such a purpose. Because ADMM algorithms operate with subproblems that are unconstrained or simply-constrained, warm-starting based on parent solutions is usually quite simple. 
 On the other side, ADMM has parameters, notably the penalty parameter $\rho$, that might also need to be updated to get fast practical convergence.
 In this regard, we are heartened by two facts:
 (i) \cite*{Kurt_linx} and \cite*{Anstreicher_BQP_entropy}  were able to inherit and occasionally quickly update the scaling parameter 
 $\gamma$ for \ref{prob_linx} and \ref{bqp_original}\,, respectively, and (ii) we saw a lot of stability 
 for good choices of $\rho$ (and other parameters) in our experiments.
Additionally, we note that there are effective adaptive methods for updating $\rho$ in the course of running an ADMM; see, for example, 
\cite*{wohlberg2017} and \cite*[Section 3.4.1]{boyd2011distributed}. 
 Although the devil is in the details, overall, we are optimistic about the possibility of ADMM as a workhorse for approaching
\ref{dopt} and \ref{MESP}
with
  B\&B,  
  the most successful algorithm for exact solution of these problems.


\section*{Acknowledgments} 
 G. Ponte was supported in part by CNPq GM-GD scholarship 161501/2022-2. M. Fampa was supported in part by CNPq grant 307167/2022-4. 
J. Lee was supported in part by AFOSR grant FA9550-22-1-0172. 

\bibliography{Dopt_FLP}

\newpage


\appendix

\section{Appendix: Performance Analysis}\label{app:performance}

We present in figures the gaps between the upper bounds for \ref{dopt} and \ref{MESP}, computed by our ADMM algorithms, and lower bounds  computed by local-search heuristics   from  \cite*{PonteFampaLeeMPB} and \cite*{KLQ}, respectively. 

We present in tables detailed results from our comparisons between our ADMM algorithms developed for \ref{prob}, \ref{prob_ddfact} and \ref{bqp_original}\,, and general-purpose solvers commonly used in the literature for these kind of problems. We also present some comparisons to the two open-source Julia implementations of first-order methods, FrankWolfe.jl and COSMO.jl. 

We show in the tables the elapsed time required by the methods to solve our instances and the final dual gap, computed as described in \S\ref{sec:numexp}. In the last column, we also present the value of the penalty parameter $\rho$ used in our experiments. In every table the symbol `$*$' indicates that the method could not solve the instance in our time limit of one hour, or due to lack of memory. 


\subsection{0/1 D-optimality}

\begin{figure}[!ht]
\setcounter{subfigure}{0}
    \centering
    \subfigure[random instances] {
    \centering
    \includegraphics[width=0.47\linewidth]{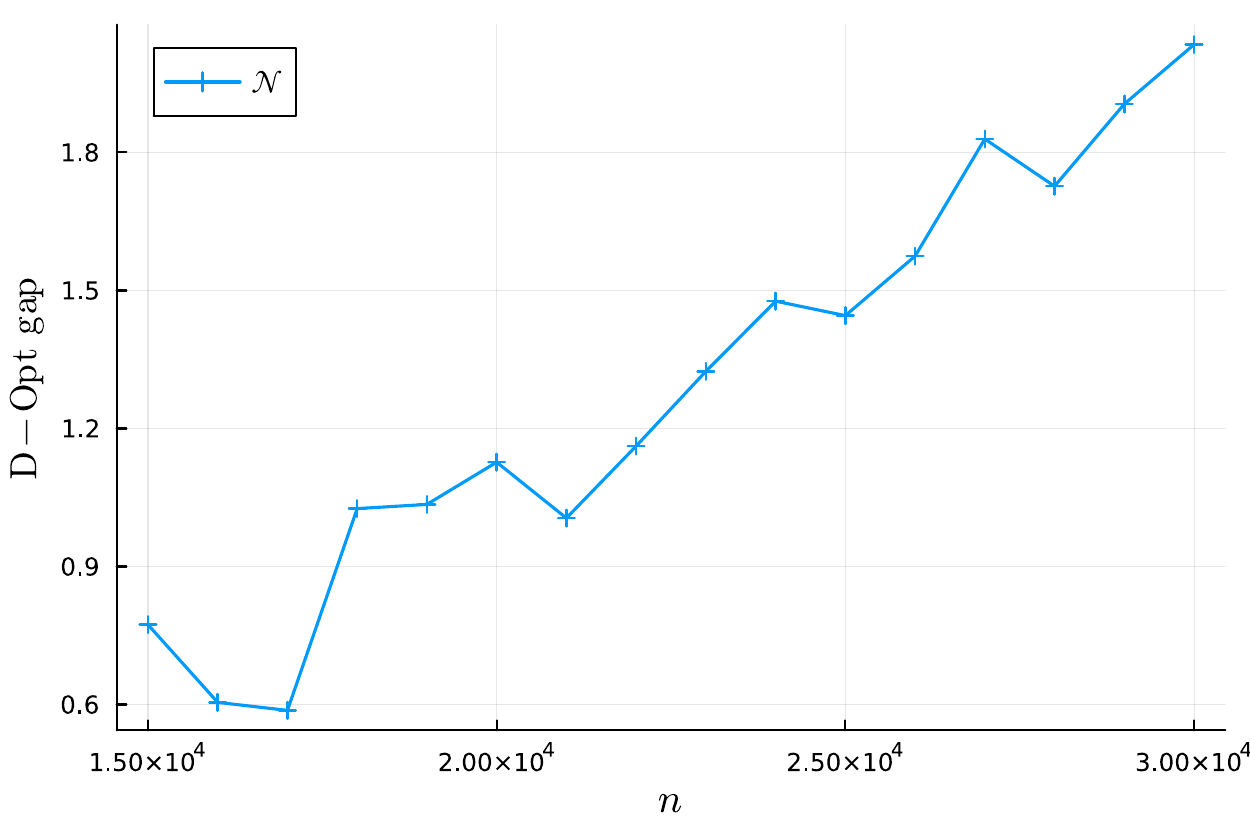}
    }
    \subfigure[linear-response model] {
    \centering
    \includegraphics[width=0.47\linewidth]{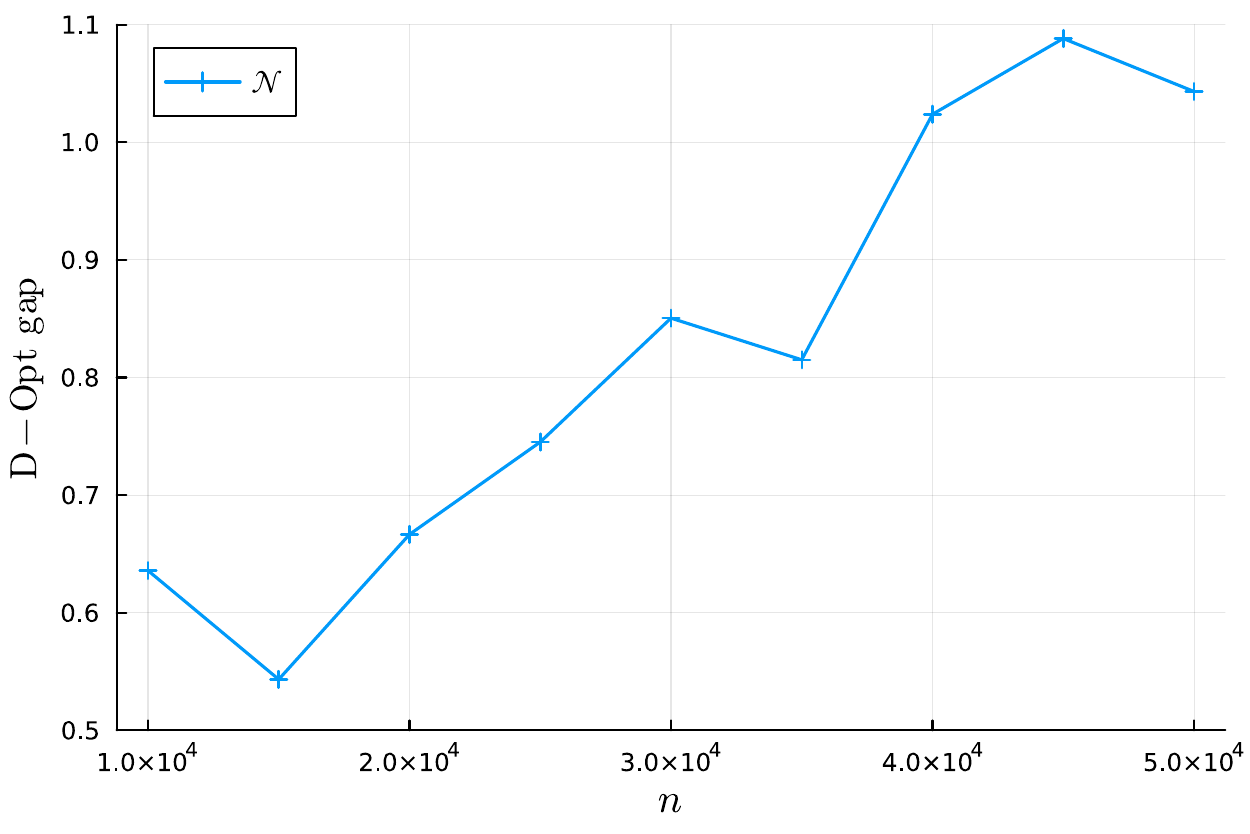}
    }
    \subfigure[quadratic-response model] {
    \centering
    \includegraphics[width=0.47\linewidth]{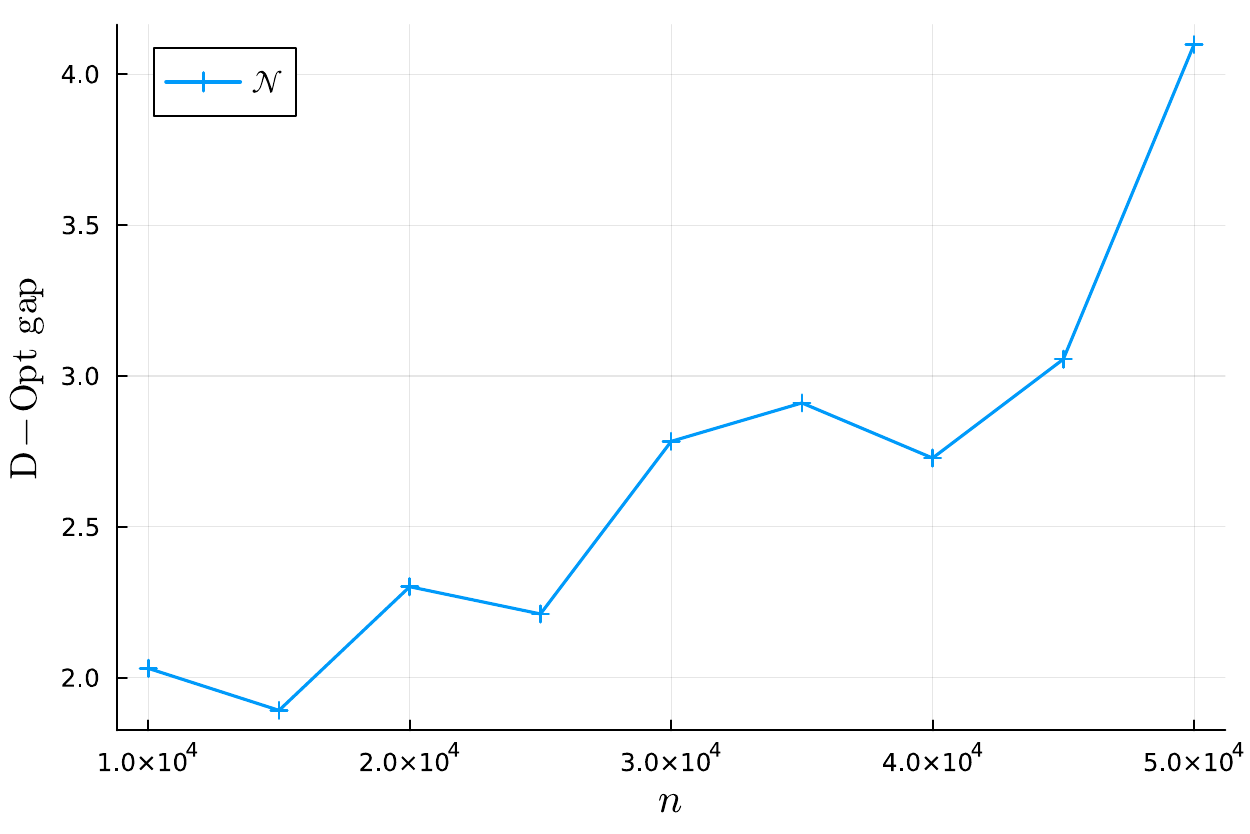}
    }
    \subfigure[real data set, $n=5822, m= 60$] {
    \centering
    \includegraphics[width=0.48\linewidth]{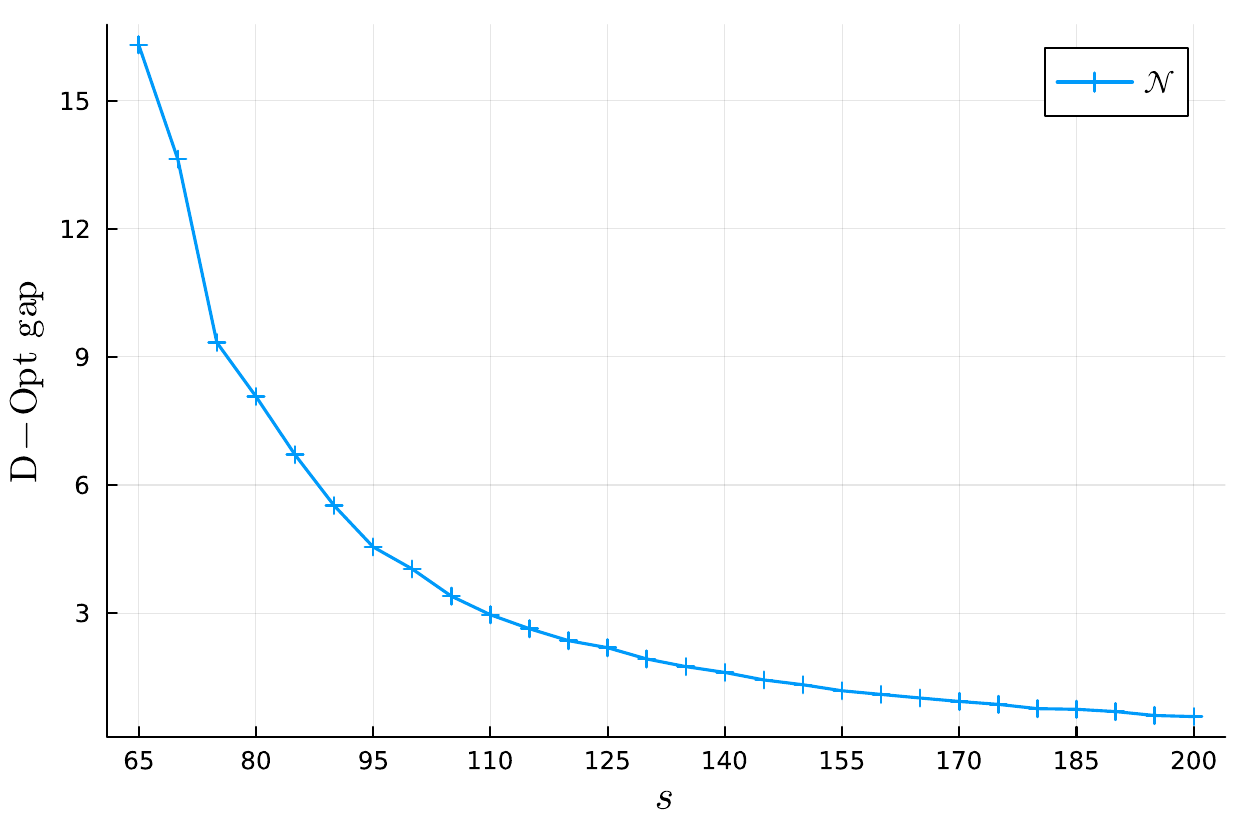}
    }
    \caption{Natural bound  \ref{prob} for \ref{dopt}
    }
    
    \label{fig:dopt_opt_gap}
\end{figure}

\begin{table}[!ht]
\footnotesize
\centering
\begin{tabular}{l|rrrrrr|rrrrrr|r}
                             & \multicolumn{6}{c|}{Elapsed time (sec)}                                                                                                                                                                    & \multicolumn{6}{c|}{Dual gap}                                                                                                                                                                              & \multicolumn{1}{c}{$\rho$}     \\
                             & \multicolumn{1}{c}{}           & \multicolumn{1}{c}{}             & \multicolumn{1}{c}{}            & \multicolumn{1}{c}{}            & \multicolumn{1}{c}{}            & \multicolumn{1}{c|}{\tiny Frank} & \multicolumn{1}{c}{}           & \multicolumn{1}{c}{}             & \multicolumn{1}{c}{}            & \multicolumn{1}{c}{}            & \multicolumn{1}{c}{}            & \multicolumn{1}{c|}{\tiny Frank} & \multicolumn{1}{c}{}           \\[-4pt]
\multicolumn{1}{c|}{$n,m,s$} & \multicolumn{1}{c}{\tiny ADMM} & \multicolumn{1}{c}{\tiny KNITRO} & \multicolumn{1}{c}{\tiny MOSEK} & \multicolumn{1}{c}{\tiny SDPT3} & \multicolumn{1}{c}{\tiny COSMO} & \multicolumn{1}{c|}{\tiny Wolfe} & \multicolumn{1}{c}{\tiny ADMM} & \multicolumn{1}{c}{\tiny KNITRO} & \multicolumn{1}{c}{\tiny MOSEK} & \multicolumn{1}{c}{\tiny SDPT3} & \multicolumn{1}{c}{\tiny COSMO} & \multicolumn{1}{c|}{\tiny Wolfe} & \multicolumn{1}{c}{\tiny ADMM} \\ \hline
15000,15,30                  & 1.0                            & 35.5                             & 27.3                            & 272.1                           & 147.2                           & 316.4                            & 4.9e-02                        & 1.1e-05                          & 1.5e-04                         & 1.1e-03                         & 2.4e-01                         & 5.1e-02                          & 2.5e-04                        \\
16000,16,32                  & 1.1                            & 47.7                             & 32.7                            & 329.2                           & 168.5                           & 442.7                            & 4.5e-02                        & 3.4e-05                          & 2.9e-04                         & 9.7e-04                         & 2.1e-01                         & 4.9e-02                          & 2.5e-04                        \\
17000,17,34                  & 1.6                            & 39.4                             & 38.6                            & 443.0                           & 257.0                           & 445.6                            & 4.6e-02                        & 2.5e-05                          & 2.6e-05                         & 1.7e-03                         & 2.1e-01                         & 5.6e-02                          & 2.5e-04                        \\
18000,18,36                  & 1.7                            & 53.3                             & 45.7                            & 510.2                           & 275.4                           & 615.7                            & 4.9e-02                        & 2.7e-05                          & 1.7e-04                         & 4.3e-04                         & 2.1e-01                         & 4.8e-02                          & 2.5e-04                        \\
19000,19,38                  & 2.7                            & 49.0                             & 55.0                            & 656.0                           & 328.1                           & 637.0                            & 4.9e-02                        & 2.4e-05                          & 2.3e-04                         & 8.3e-04                         & 2.7e-01                         & 5.3e-02                          & 2.5e-04                        \\
20000,20,40                  & 2.9                            & 49.5                             & 63.4                            & 815.7                           & 415.7                           & 728.2                            & 4.5e-02                        & 3.2e-05                          & 2.4e-04                         & 7.8e-03                         & 3.1e-01                         & 5.4e-02                          & 1.0e-04                        \\
21000,21,42                  & 3.7                            & 66.3                             & 72.4                            & *                               & 486.6                           & 892.3                            & 5.0e-02                        & 1.6e-05                          & 2.8e-04                         & *                               & 3.3e-01                         & 5.6e-02                          & 1.0e-04                        \\
22000,22,44                  & 3.7                            & 50.7                             & 83.6                            & *                               & 500.5                           & 1076.5                           & 4.9e-02                        & 2.9e-05                          & 4.8e-04                         & *                               & 2.9e-01                         & 4.9e-02                          & 1.0e-04                        \\
23000,23,46                  & 4.4                            & 56.5                             & 95.1                            & *                               & 515.0                           & 1150.3                           & 4.4e-02                        & 2.6e-05                          & 3.1e-04                         & *                               & 3.2e-01                         & 5.2e-02                          & 1.0e-04                        \\
24000,24,48                  & 4.1                            & 66.3                             & 110.9                           & *                               & 693.6                           & 1296.7                           & 4.7e-02                        & 3.6e-05                          & 3.3e-04                         & *                               & 3.0e-01                         & 5.1e-02                          & 1.0e-04                        \\
25000,25,50                  & 4.2                            & 65.5                             & 127.5                           & *                               & 656.4                           & 1548.2                           & 4.5e-02                        & 3.5e-05                          & 7.8e-05                         & *                               & 3.5e-01                         & 5.3e-02                          & 1.0e-04                        \\
26000,26,52                  & 4.9                            & 70.7                             & 140.1                           & *                               & 836.9                           & 1836.3                           & 4.1e-02                        & 2.4e-05                          & 3.8e-05                         & *                               & 3.4e-01                         & 5.1e-02                          & 1.0e-04                        \\
27000,27,54                  & 5.9                            & 76.5                             & 158.0                           & *                               & 792.2                           & 1925.9                           & 4.2e-02                        & 4.6e-05                          & 2.0e-04                         & *                               & 3.2e-01                         & 5.3e-02                          & 1.0e-04                        \\
28000,28,56                  & 7.2                            & 80.3                             & 175.3                           & *                               & 1297.8                          & 2188.5                           & 4.6e-02                        & 4.3e-05                          & 1.3e-04                         & *                               & 3.1e-01                         & 5.6e-02                          & 5.0e-05                        \\
29000,29,58                  & 8.1                            & 85.6                             & 193.4                           & *                               & 1115.0                          & 2709.7                           & 4.8e-02                        & 5.0e-05                          & 5.0e-04                         & *                               & 3.6e-01                         & 5.4e-02                          & 5.0e-05                        \\
30000,30,60                  & 5.7                            & 91.7                             & 210.7                           & *                               & 918.4                           & 2468.7                           & 4.8e-02                        & 2.3e-05                          & 5.4e-04                         & *                               & 3.1e-01                         & 5.9e-02                          & 5.0e-05                       
\end{tabular}
\caption{Random instances: \ref{dopt}}\label{tab:random}
\end{table}

\begin{table}[!ht]
\centering
\footnotesize
\begin{tabular}{l|rrrrrr|rrrrrr|r}
                             & \multicolumn{6}{c|}{Elapsed time (sec)}                                                                                                                                                                    & \multicolumn{6}{c|}{Dual gap}                                                                                                                                                                              & \multicolumn{1}{c}{$\rho$}     \\
                             & \multicolumn{1}{c}{}           & \multicolumn{1}{c}{}             & \multicolumn{1}{c}{}            & \multicolumn{1}{c}{}            & \multicolumn{1}{c}{}            & \multicolumn{1}{c|}{\tiny Frank} & \multicolumn{1}{c}{}           & \multicolumn{1}{c}{}             & \multicolumn{1}{c}{}            & \multicolumn{1}{c}{}            & \multicolumn{1}{c}{}            & \multicolumn{1}{c|}{\tiny Frank} & \multicolumn{1}{c}{}           \\
\multicolumn{1}{c|}{$n,m,s$} & \multicolumn{1}{c}{\tiny ADMM} & \multicolumn{1}{c}{\tiny KNITRO} & \multicolumn{1}{c}{\tiny MOSEK} & \multicolumn{1}{c}{\tiny SDPT3} & \multicolumn{1}{c}{\tiny COSMO} & \multicolumn{1}{c|}{\tiny Wolfe} & \multicolumn{1}{c}{\tiny ADMM} & \multicolumn{1}{c}{\tiny KNITRO} & \multicolumn{1}{c}{\tiny MOSEK} & \multicolumn{1}{c}{\tiny SDPT3} & \multicolumn{1}{c}{\tiny COSMO} & \multicolumn{1}{c|}{\tiny Wolfe} & \multicolumn{1}{c}{\tiny ADMM} \\ \hline
10000,20,40                  & 1.3                            & 19.3                             & 16.2                            & 90.9                            & 25.4                            & 40.4                             & 4.8e-02                        & 1.6e-05                          & 7.0e-04                         & 9.40e-05                        & 6.5e-04                         & 5.0e-02                          & 2.5e-02                        \\
15000,21,42                  & 1.2                            & 27.2                             & 36.6                            & 222.3                           & 54.7                            & 52.7                             & 4.0e-02                        & 1.4e-06                          & 1.6e-03                         & 2.55e-04                        & 2.9e-04                         & 5.1e-02                          & 2.5e-02                        \\
20000,22,44                  & 1.4                            & 41.0                             & 66.1                            & 386.3                           & 96.1                            & 102.1                            & 3.7e-02                        & 6.1e-06                          & 2.4e-03                         & 3.36e-04                        & 3.0e-04                         & 5.5e-02                          & 2.5e-02                        \\
25000,23,46                  & 1.4                            & 50.1                             & 109.4                           & 759.0                           & 153.1                           & 147.0                            & 4.8e-02                        & 2.8e-05                          & 2.5e-03                         & 4.21e-05                        & 7.9e-05                         & 5.9e-02                          & 2.5e-02                        \\
30000,24,48                  & 1.6                            & 62.1                             & 160.6                           & *                               & 224.0                           & 205.2                            & 3.1e-02                        & 1.7e-05                          & 3.3e-03                         & *                               & 5.0e-04                         & 5.0e-02                          & 2.5e-02                        \\
35000,25,50                  & 2.1                            & 86.2                             & 228.6                           & *                               & 316.6                           & 278.0                            & 2.8e-02                        & 5.0e-06                          & 2.0e-03                         & *                               & 2.4e-04                         & 6.5e-02                          & 2.5e-02                        \\
40000,26,52                  & 2.4                            & 117.9                            & 306.2                           & *                               & 423.9                           & 353.1                            & 3.9e-02                        & 1.4e-05                          & 3.8e-03                         & *                               & 3.0e-04                         & 6.5e-02                          & 2.5e-02                        \\
45000,27,54                  & 2.6                            & 164.9                            & 404.9                           & *                               & 549.9                           & 428.4                            & 2.7e-02                        & 6.3e-05                          & 5.4e-03                         & *                               & 1.2e-03                         & 8.2e-02                          & 2.5e-02                        \\
50000,28,56                  & 2.8                            & 160.8                            & 520.5                           & *                               & 701.3                           & 482.5                            & 2.7e-02                        & 5.5e-05                          & 3.7e-03                         & *                               & 8.8e-04                         & 7.9e-02                          & 2.5e-02                       
\end{tabular}
\caption{Linear-response model: \ref{dopt}}\label{tab:linmodel}
\end{table}

\begin{table}[!ht]
\centering
\footnotesize
\begin{tabular}{l|rrrrrr|rrrrrr|r}
                             & \multicolumn{6}{c|}{Elapsed time (sec)}                                                                                                                                                                    & \multicolumn{6}{c|}{Dual gap}                                                                                                                                                                              & \multicolumn{1}{c}{$\rho$}     \\
                             & \multicolumn{1}{c}{}           & \multicolumn{1}{c}{}             & \multicolumn{1}{c}{}            & \multicolumn{1}{c}{}            & \multicolumn{1}{c}{}            & \multicolumn{1}{c|}{\tiny Frank} & \multicolumn{1}{c}{}           & \multicolumn{1}{c}{}             & \multicolumn{1}{c}{}            & \multicolumn{1}{c}{}            & \multicolumn{1}{c}{}            & \multicolumn{1}{c|}{\tiny Frank} & \multicolumn{1}{c}{}           \\[-4pt]
\multicolumn{1}{c|}{$n,m,s$} & \multicolumn{1}{c}{\tiny ADMM} & \multicolumn{1}{c}{\tiny KNITRO} & \multicolumn{1}{c}{\tiny MOSEK} & \multicolumn{1}{c}{\tiny SDPT3} & \multicolumn{1}{c}{\tiny COSMO} & \multicolumn{1}{c|}{\tiny Wolfe} & \multicolumn{1}{c}{\tiny ADMM} & \multicolumn{1}{c}{\tiny KNITRO} & \multicolumn{1}{c}{\tiny MOSEK} & \multicolumn{1}{c}{\tiny SDPT3} & \multicolumn{1}{c}{\tiny COSMO} & \multicolumn{1}{c|}{\tiny Wolfe} & \multicolumn{1}{c}{\tiny ADMM} \\ \hline
10000,30,60                  & 14.8                           & 29.9                             & 33.0                            & 164.4                           & 247.6                           & 747.0                            & 5.0e-02                        & 2.7e-05                          & 4.0e-04                         & 1.6e-03                         & 8.6e-02                         & 5.5e-02                          & 7.0e-04                        \\
15000,31,62                  & 25.1                           & 58.5                             & 72.8                            & 357.1                           & 315.5                           & 1420.3                           & 5.0e-02                        & 2.0e-05                          & 1.0e-03                         & 1.3e-02                         & 1.4e-01                         & 4.9e-02                          & 7.0e-04                        \\
20000,32,64                  & 31.4                           & 90.5                             & 122.9                           & 821.4                           & 567.1                           & 2476.2                           & 5.0e-02                        & 5.7e-05                          & 1.5e-03                         & 3.4e-03                         & 1.9e-01                         & 5.3e-02                          & 7.0e-04                        \\
25000,33,66                  & 45.7                           & 128.5                            & 190.6                           & *                               & 841.5                           & *                                & 5.0e-02                        & 4.3e-05                          & 1.1e-03                         & *                               & 2.3e-01                         & *                                & 6.0e-04                        \\
30000,39,78                  & 83.3                           & 203.5                            & 345.0                           & *                               & 2086.2                          & *                                & 5.0e-02                        & 2.4e-05                          & 5.0e-04                         & *                               & 4.2e-01                         & *                                & 6.0e-04                        \\
35000,40,80                  & 96.7                           & 228.1                            & 456.5                           & *                               & 1884.2                          & *                                & 5.0e-02                        & 5.8e-05                          & 1.7e-03                         & *                               & 4.4e-01                         & *                                & 5.0e-04                        \\
40000,41,82                  & 128.7                          & 282.2                            & 625.7                           & *                               & 3154.0                          & *                                & 5.0e-02                        & 3.7e-05                          & 2.8e-03                         & *                               & 5.0e-01                         & *                                & 5.0e-04                        \\
45000,42,84                  & 165.3                          & 362.3                            & 788.0                           & *                               & 3374.0                          & *                                & 5.0e-02                        & 4.7e-05                          & 6.2e-04                         & *                               & 5.5e-01                         & *                                & 4.0e-04                        \\
50000,49,98                  & 219.1                          & 430.0                            & 1226.7                          & *                               & *                               & *                                & 5.0e-02                        & 6.7e-05                          & 1.7e-03                         & *                               & *                               & *                                & 4.0e-04                       
\end{tabular}
\caption{Quadratic-response model: \ref{dopt}}\label{tab:quadmodel}
\end{table}

\begin{table}[!ht]
\centering
\footnotesize
\begin{tabular}{r|rrrrrr|rrrrrr|r}
                         & \multicolumn{6}{c|}{Elapsed time (sec)}                                                                                                                                                                    & \multicolumn{6}{c|}{Dual gap}                                                                                                                                                                              & \multicolumn{1}{c}{$\rho$}     \\
                         & \multicolumn{1}{c}{}           & \multicolumn{1}{c}{}             & \multicolumn{1}{c}{}            & \multicolumn{1}{c}{}            & \multicolumn{1}{c}{}            & \multicolumn{1}{c|}{\tiny Frank} & \multicolumn{1}{c}{}           & \multicolumn{1}{c}{}             & \multicolumn{1}{c}{}            & \multicolumn{1}{c}{}            & \multicolumn{1}{c}{}            & \multicolumn{1}{c|}{\tiny Frank} & \multicolumn{1}{c}{}           \\[-4pt]
\multicolumn{1}{c|}{$s$} & \multicolumn{1}{c}{\tiny ADMM} & \multicolumn{1}{c}{\tiny KNITRO} & \multicolumn{1}{c}{\tiny MOSEK} & \multicolumn{1}{c}{\tiny SDPT3} & \multicolumn{1}{c}{\tiny COSMO} & \multicolumn{1}{c|}{\tiny Wolfe} & \multicolumn{1}{c}{\tiny ADMM} & \multicolumn{1}{c}{\tiny KNITRO} & \multicolumn{1}{c}{\tiny MOSEK} & \multicolumn{1}{c}{\tiny SDPT3} & \multicolumn{1}{c}{\tiny COSMO} & \multicolumn{1}{c|}{\tiny Wolfe} & \multicolumn{1}{c}{\tiny ADMM} \\ \hline
65                       & 12.2                           & 63.9                             & 137.6                           & 280.4                           & *                               & 1910.2                           & 4.7e-02                        & 4.0e-05                          & 4.4e-03                         & 4.2e-03                         & *                               & 5.4e-02                          & 3.0e-03                        \\
70                       & 10.5                           & 65.5                             & 147.4                           & 276.0                           & *                               & 1617.8                           & 5.0e-02                        & 3.2e-05                          & 3.1e-03                         & 1.9e-03                         & *                               & 4.7e-02                          & 3.0e-03                        \\
75                       & 10.4                           & 56.4                             & 147.8                           & 277.2                           & *                               & 1462.7                           & 4.9e-02                        & 2.6e-05                          & 2.9e-03                         & 1.6e-03                         & *                               & 5.2e-02                          & 3.0e-03                        \\
80                       & 10.5                           & 58.1                             & 146.8                           & 273.7                           & *                               & 1295.3                           & 4.9e-02                        & 3.5e-05                          & 3.6e-03                         & 1.7e-03                         & *                               & 5.6e-02                          & 2.0e-03                        \\
85                       & 10.1                           & 52.6                             & 145.9                           & 281.9                           & *                               & 1203.3                           & 4.9e-02                        & 2.3e-05                          & 2.3e-03                         & 2.6e-03                         & *                               & 5.5e-02                          & 2.0e-03                        \\
90                       & 10.2                           & 58.3                             & 144.8                           & 284.5                           & *                               & 1040.8                           & 4.6e-02                        & 4.0e-05                          & 5.9e-03                         & 3.3e-03                         & *                               & 5.5e-02                          & 2.0e-03                        \\
95                       & 10.5                           & 49.3                             & 140.9                           & 283.2                           & *                               & 971.2                            & 5.0e-02                        & 3.7e-05                          & 6.3e-03                         & 1.0e-02                         & *                               & 5.0e-02                          & 2.0e-03                        \\
100                      & 12.7                           & 49.8                             & 144.1                           & 289.2                           & *                               & 932.8                            & 4.8e-02                        & 3.2e-05                          & 4.5e-03                         & 4.5e-03                         & *                               & 4.9e-02                          & 1.0e-03                        \\
105                      & 11.6                           & 55.0                             & 138.5                           & 283.9                           & *                               & 847.0                            & 4.7e-02                        & 2.1e-05                          & 3.9e-03                         & 8.4e-03                         & *                               & 4.8e-02                          & 1.0e-03                        \\
110                      & 10.7                           & 56.8                             & 142.9                           & 280.6                           & *                               & 768.5                            & 4.7e-02                        & 1.4e-05                          & 5.8e-03                         & 1.4e-02                         & *                               & 5.1e-02                          & 1.0e-03                        \\
115                      & 9.7                            & 49.7                             & 145.8                           & 270.4                           & *                               & 695.6                            & 4.9e-02                        & 2.3e-05                          & 4.9e-03                         & 1.0e-02                         & *                               & 5.2e-02                          & 1.0e-03                        \\
120                      & 9.3                            & 54.3                             & 156.4                           & 280.9                           & *                               & 655.4                            & 5.0e-02                        & 3.8e-05                          & 3.6e-03                         & 2.3e-02                         & *                               & 5.0e-02                          & 1.0e-03                        \\
125                      & 8.9                            & 52.9                             & 177.3                           & 292.0                           & *                               & 604.8                            & 4.9e-02                        & 3.8e-05                          & 5.1e-03                         & 1.2e-03                         & *                               & 5.1e-02                          & 1.0e-03                        \\
130                      & 8.7                            & 53.5                             & 147.0                           & 271.0                           & *                               & 538.9                            & 4.9e-02                        & 4.8e-05                          & 2.1e-03                         & 3.4e-02                         & *                               & 5.2e-02                          & 1.0e-03                        \\
135                      & 8.8                            & 56.7                             & 174.7                           & 280.6                           & *                               & 521.9                            & 5.0e-02                        & 2.9e-05                          & 4.9e-03                         & 2.1e-02                         & *                               & 5.1e-02                          & 1.0e-03                        \\
140                      & 9.3                            & 58.6                             & 172.2                           & 273.1                           & *                               & 483.2                            & 4.7e-02                        & 4.0e-05                          & 2.1e-03                         & 1.5e-02                         & *                               & 5.6e-02                          & 1.0e-03                        \\
145                      & 9.9                            & 56.5                             & 174.7                           & 279.2                           & *                               & 441.5                            & 4.8e-02                        & 4.4e-05                          & 2.7e-03                         & 1.1e-02                         & *                               & 5.1e-02                          & 1.0e-03                        \\
150                      & 9.9                            & 59.2                             & 173.6                           & 281.4                           & *                               & 429.2                            & 4.9e-02                        & 5.9e-05                          & 2.8e-03                         & 1.1e-02                         & *                               & 5.0e-02                          & 1.0e-03                        \\
155                      & 10.9                           & 68.9                             & 161.9                           & 285.5                           & *                               & 426.7                            & 4.8e-02                        & 2.4e-05                          & 4.3e-03                         & 1.4e-02                         & *                               & 4.9e-02                          & 1.0e-03                        \\
160                      & 11.4                           & 93.6                             & 173.9                           & 282.9                           & *                               & 398.8                            & 4.9e-02                        & 4.0e-05                          & 2.7e-03                         & 1.5e-02                         & *                               & 4.9e-02                          & 1.0e-03                        \\
165                      & 11.8                           & 56.3                             & 156.3                           & 269.9                           & *                               & 376.5                            & 4.9e-02                        & 4.0e-05                          & 3.0e-03                         & 2.8e-02                         & *                               & 4.9e-02                          & 1.0e-03                        \\
170                      & 12.4                           & 70.6                             & 173.4                           & 279.7                           & *                               & 353.8                            & 4.9e-02                        & 2.6e-05                          & 2.3e-03                         & 3.0e-03                         & *                               & 4.8e-02                          & 1.0e-03                        \\
175                      & 12.5                           & 81.8                             & 175.0                           & 285.0                           & *                               & 330.1                            & 4.9e-02                        & 1.4e-05                          & 4.6e-03                         & 1.0e-02                         & *                               & 4.7e-02                          & 1.0e-03                        \\
180                      & 12.9                           & 64.5                             & 174.0                           & 283.4                           & *                               & 326.8                            & 5.0e-02                        & 3.9e-05                          & 1.5e-03                         & 1.3e-02                         & *                               & 5.2e-02                          & 1.0e-03                        \\
185                      & 13.3                           & 70.1                             & 176.8                           & 285.1                           & *                               & 319.9                            & 5.0e-02                        & 2.7e-05                          & 2.0e-03                         & 1.7e-02                         & *                               & 4.8e-02                          & 1.0e-03                        \\
190                      & 13.8                           & 43.8                             & 156.3                           & 291.7                           & *                               & 293.9                            & 4.9e-02                        & 3.3e-05                          & 2.0e-03                         & 1.6e-02                         & *                               & 5.4e-02                          & 1.0e-03                        \\
195                      & 14.3                           & 70.0                             & 171.0                           & 286.9                           & *                               & 312.8                            & 5.0e-02                        & 4.5e-05                          & 3.5e-03                         & 2.1e-02                         & *                               & 5.0e-02                          & 1.0e-03                        \\
200                      & 15.2                           & 73.9                             & 175.1                           & 286.8                           & *                               & 272.7                            & 5.0e-02                        & 2.1e-05                          & 1.8e-03                         & 1.9e-02                         & *                               & 5.7e-02                          & 1.0e-03                       
\end{tabular}
\caption{Real instance $n=5822, m=60$: \ref{dopt}}\label{tab:real}
\end{table}


\FloatBarrier

\subsection{MESP}

\begin{figure}[!ht]
\setcounter{subfigure}{0}
    \centering
    \subfigure[ varying $r:=\rank(C)$ ($s=140$)] {
    \centering
    \includegraphics[width=0.47\linewidth]{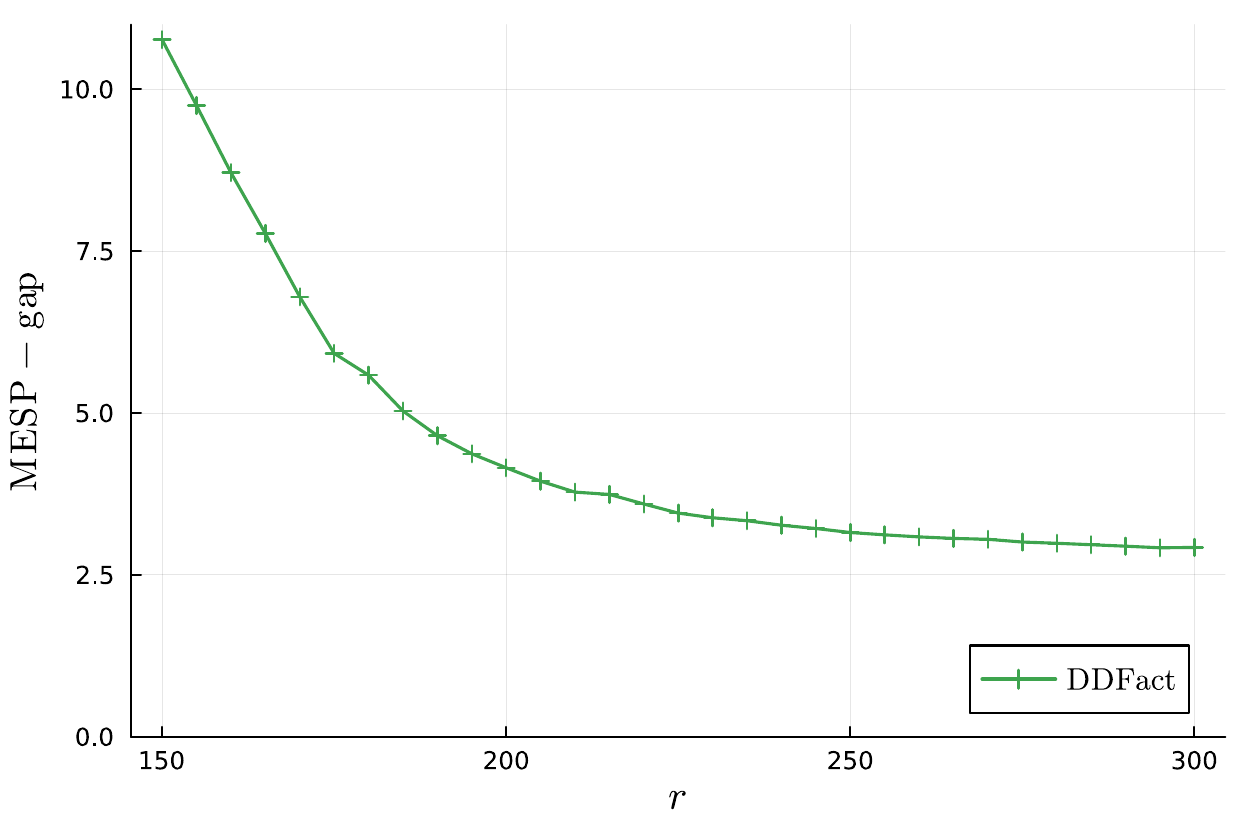}
    }
    \subfigure[varying $s$ ($\rank(C)=150$)] {
    \centering
    \includegraphics[width=0.47\linewidth]{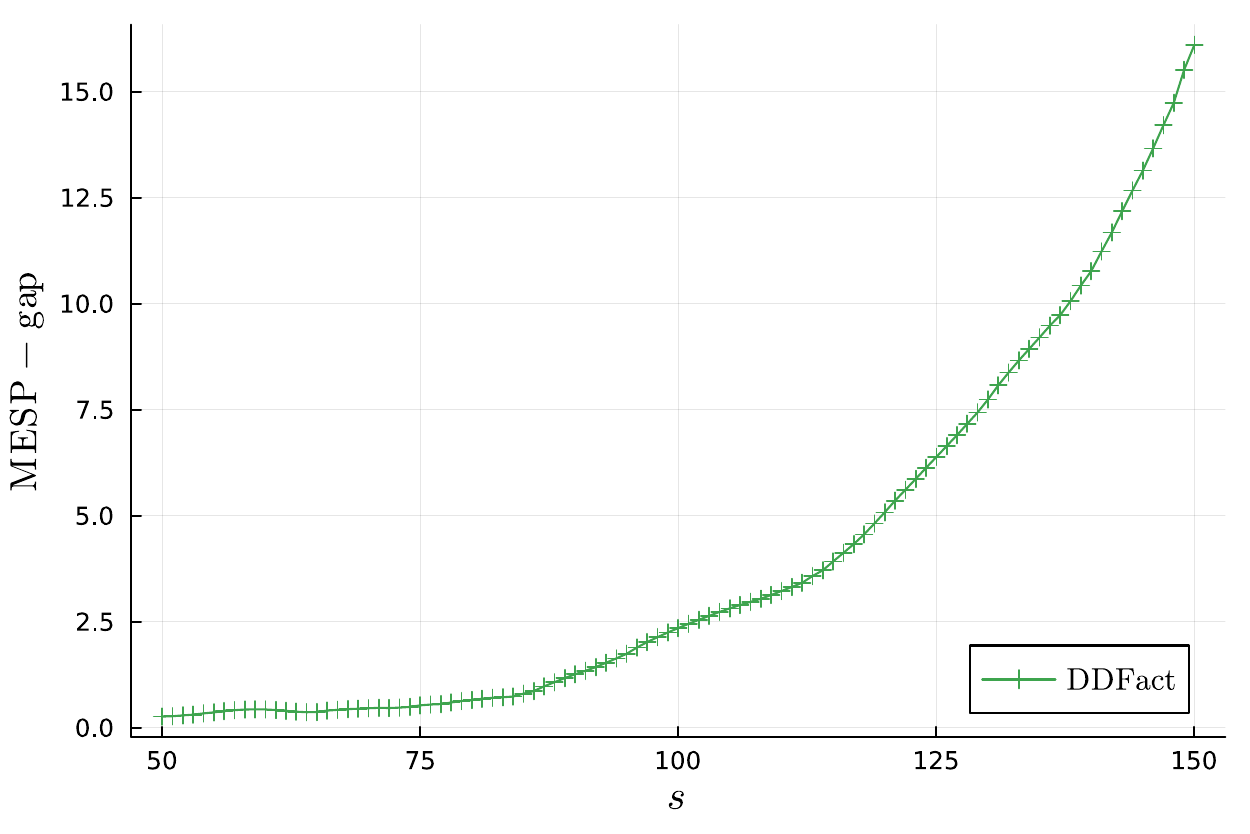}
    }\\
    \caption{\ref{prob_ddfact} bound for \ref{MESP}  ($n=2000$)}
    \label{fig_LS_bound_app_alone}
\end{figure}

\begin{figure}[!ht]
\setcounter{subfigure}{0}
    \centering
    \subfigure[ varying $r:=\rank(C)$ ($s=140$)] {
    \centering
    \includegraphics[width=0.47\linewidth]{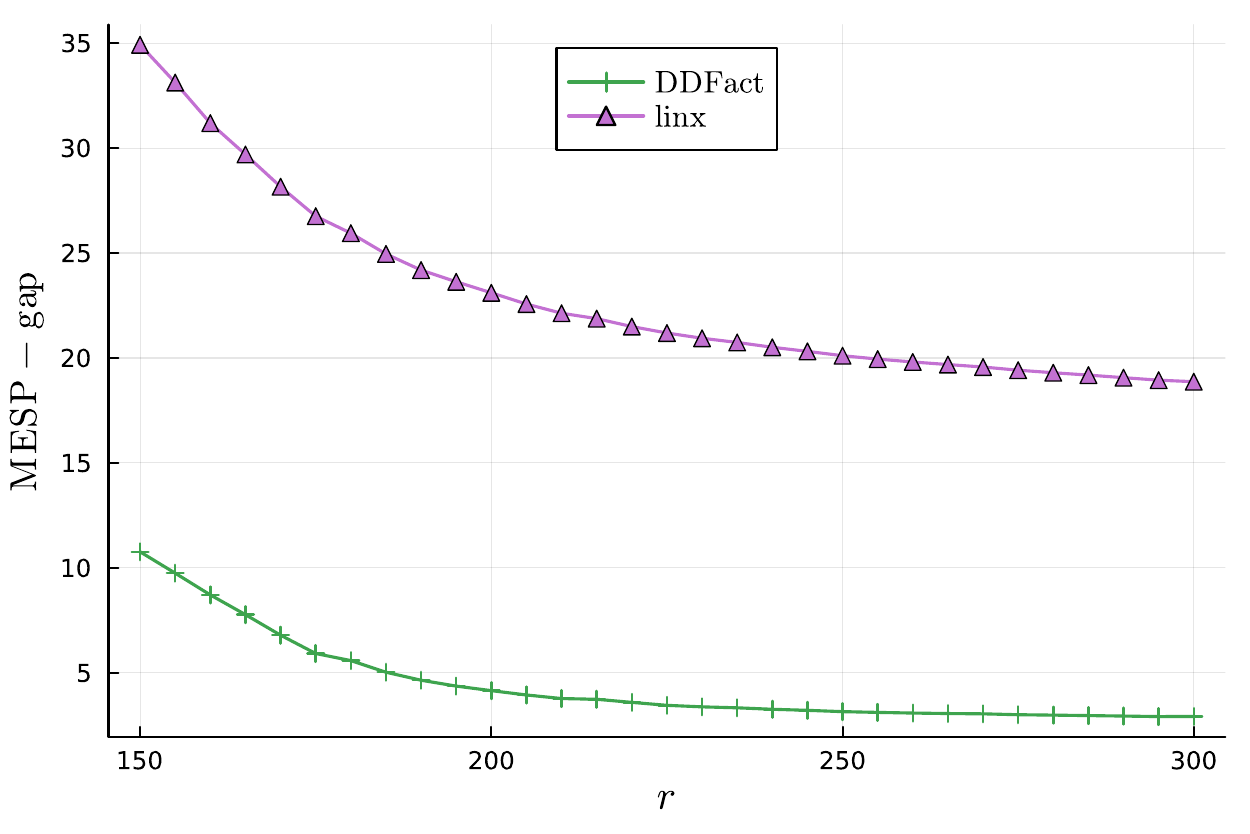}
    }
    \subfigure[varying $s$ ($\rank(C)=150$)] {
    \centering
    \includegraphics[width=0.47\linewidth]{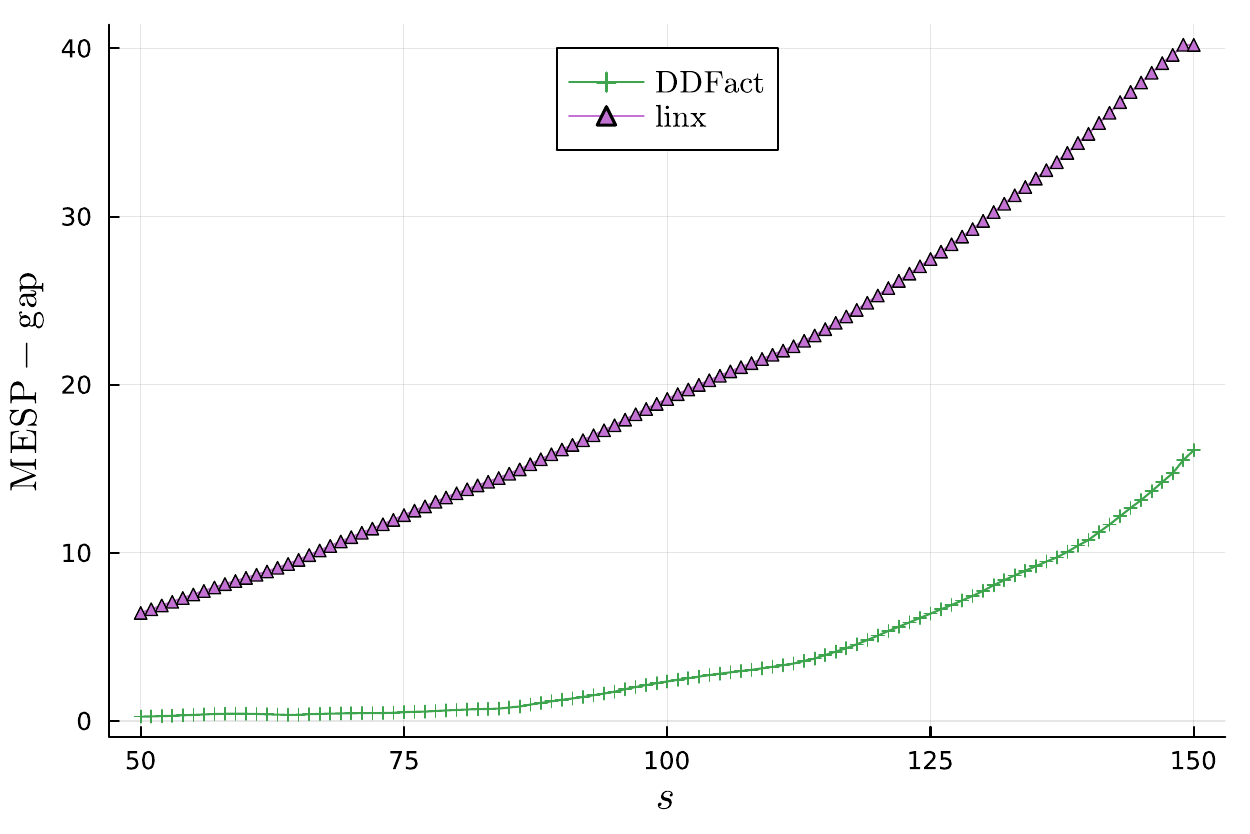}
    }\\
    \caption{\ref{prob_ddfact} and \ref{prob_linx} bound for \ref{MESP}  ($n=2000$)}
    \label{fig_LS_bound_app}
\end{figure}

\begin{table}[!ht]
\centering
\footnotesize
\begin{tabular}{l|ccc|rrr|r}
\multicolumn{1}{c|}{}    & \multicolumn{3}{c|}{Elapsed time (sec)}                                                              & \multicolumn{3}{c|}{Dual gap}                                                                        & \multicolumn{1}{c}{$\rho$}     \\
\multicolumn{1}{c|}{}    & \multicolumn{1}{c}{}           & \multicolumn{1}{c}{}             & \multicolumn{1}{c|}{\tiny Frank} & \multicolumn{1}{c}{}           & \multicolumn{1}{c}{}             & \multicolumn{1}{c|}{\tiny Frank} & \multicolumn{1}{c}{}           \\[-4pt]
\multicolumn{1}{c|}{$r$} & \multicolumn{1}{c}{\tiny ADMM} & \multicolumn{1}{c}{\tiny KNITRO} & \multicolumn{1}{c|}{\tiny Wolfe} & \multicolumn{1}{c}{\tiny ADMM} & \multicolumn{1}{c}{\tiny KNITRO} & \multicolumn{1}{c|}{\tiny Wolfe} & \multicolumn{1}{c}{\tiny ADMM} \\ \hline
150                      & 2.46                           & 4.16                             & 454.69                           & 3.7e-02                        & 3.2e-05                          & 5.4e-02                          & 2.0e-03                        \\
155                      & 1.32                           & 4.00                             & 458.58                           & 2.6e-02                        & 4.5e-05                          & 5.3e-02                          & 2.0e-03                        \\
160                      & 1.29                           & 3.39                             & 414.88                           & 3.7e-02                        & 3.8e-05                          & 5.3e-02                          & 2.0e-03                        \\
165                      & 1.41                           & 3.16                             & 447.05                           & 2.8e-02                        & 7.5e-05                          & 5.4e-02                          & 2.0e-03                        \\
170                      & 1.45                           & 3.50                             & 434.54                           & 3.5e-02                        & 2.7e-05                          & 5.4e-02                          & 3.2e-03                        \\
175                      & 2.72                           & 4.93                             & 446.80                           & 4.5e-02                        & 4.4e-05                          & 5.0e-02                          & 3.2e-03                        \\
180                      & 2.74                           & 4.72                             & 428.92                           & 4.7e-02                        & 3.1e-05                          & 5.0e-02                          & 3.2e-03                        \\
185                      & 2.63                           & 4.34                             & 397.94                           & 4.7e-02                        & 2.3e-05                          & 5.2e-02                          & 3.2e-03                        \\
190                      & 2.85                           & 5.65                             & 401.22                           & 4.5e-02                        & 3.6e-05                          & 5.1e-02                          & 3.2e-03                        \\
195                      & 3.02                           & 4.35                             & 378.29                           & 3.6e-02                        & 1.4e-05                          & 5.8e-02                          & 3.2e-03                        \\
200                      & 2.45                           & 3.80                             & 381.76                           & 4.3e-02                        & 4.6e-05                          & 7.7e-02                          & 3.2e-03                        \\
205                      & 2.86                           & 4.21                             & 371.57                           & 3.1e-02                        & 2.3e-05                          & 5.2e-02                          & 3.5e-03                        \\
210                      & 2.84                           & 4.61                             & 382.16                           & 4.2e-02                        & 2.9e-05                          & 5.0e-02                          & 3.5e-03                        \\
215                      & 3.04                           & 4.21                             & 378.98                           & 4.4e-02                        & 1.5e-05                          & 5.2e-02                          & 3.5e-03                        \\
220                      & 3.19                           & 4.02                             & 367.74                           & 4.6e-02                        & 2.4e-05                          & 5.2e-02                          & 3.5e-03                        \\
225                      & 3.50                           & 4.20                             & 379.23                           & 4.4e-02                        & 1.4e-05                          & 5.4e-02                          & 3.5e-03                        \\
230                      & 3.65                           & 4.85                             & 384.66                           & 4.7e-02                        & 2.9e-05                          & 5.6e-02                          & 3.5e-03                        \\
235                      & 3.98                           & 4.66                             & 392.23                           & 4.2e-02                        & 1.1e-05                          & 5.3e-02                          & 3.5e-03                        \\
240                      & 3.96                           & 4.35                             & 381.85                           & 4.9e-02                        & 3.2e-05                          & 5.9e-02                          & 3.5e-03                        \\
245                      & 3.81                           & 4.97                             & 404.05                           & 4.9e-02                        & 2.7e-05                          & 5.4e-02                          & 3.5e-03                        \\
250                      & 4.34                           & 4.58                             & 398.35                           & 4.2e-02                        & 2.1e-05                          & 5.6e-02                          & 3.5e-03                        \\
255                      & 4.31                           & 4.42                             & 403.14                           & 4.7e-02                        & 5.0e-05                          & 5.6e-02                          & 3.5e-03                        \\
260                      & 4.40                           & 4.78                             & 404.06                           & 4.6e-02                        & 3.7e-05                          & 5.4e-02                          & 3.5e-03                        \\
265                      & 4.57                           & 4.74                             & 406.68                           & 4.4e-02                        & 1.7e-05                          & 5.6e-02                          & 3.5e-03                        \\
270                      & 4.83                           & 6.56                             & 414.16                           & 4.5e-02                        & 1.8e-05                          & 5.7e-02                          & 3.5e-03                        \\
275                      & 4.83                           & 6.38                             & 418.91                           & 4.4e-02                        & 1.8e-05                          & 5.2e-02                          & 3.5e-03                        \\
280                      & 5.04                           & 6.25                             & 422.48                           & 4.3e-02                        & 1.3e-05                          & 7.2e-02                          & 3.5e-03                        \\
285                      & 5.15                           & 5.65                             & 430.03                           & 4.3e-02                        & 2.1e-05                          & 6.2e-02                          & 3.5e-03                        \\
290                      & 5.62                           & 5.55                             & 432.12                           & 4.3e-02                        & 2.6e-05                          & 5.0e-02                          & 3.5e-03                        \\
295                      & 5.84                           & 5.67                             & 453.87                           & 4.2e-02                        & 2.1e-05                          & 5.1e-02                          & 3.5e-03                        \\
300                      & 6.62                           & 5.01                             & 449.03                           & 4.3e-02                        & 3.3e-05                          & 6.8e-02                          & 3.5e-03       
\end{tabular}
\caption{\ref{prob_ddfact} bound for \ref{MESP}, varying $r:=\rank(C)$ ($n=2000$, $s=140$)}\label{tab:ddfact_r}
\end{table}

\begin{table}[!ht]
\centering
\footnotesize
\begin{tabular}{r|ccc|rrr|r}
\multicolumn{1}{c|}{}    & \multicolumn{3}{c|}{Elapsed time (sec)}                                                              & \multicolumn{3}{c|}{Dual gap}                                                                        & \multicolumn{1}{c}{$\rho$}     \\
\multicolumn{1}{c|}{}    & \multicolumn{1}{c}{}           & \multicolumn{1}{c}{}             & \multicolumn{1}{c|}{\tiny Frank} & \multicolumn{1}{c}{}           & \multicolumn{1}{c}{}             & \multicolumn{1}{c|}{\tiny Frank} & \multicolumn{1}{c}{}           \\[-4pt]
\multicolumn{1}{c|}{$r$} & \multicolumn{1}{c}{\tiny ADMM} & \multicolumn{1}{c}{\tiny KNITRO} & \multicolumn{1}{c|}{\tiny Wolfe} & \multicolumn{1}{c}{\tiny ADMM} & \multicolumn{1}{c}{\tiny KNITRO} & \multicolumn{1}{c|}{\tiny Wolfe} & \multicolumn{1}{c}{\tiny ADMM} \\ \hline
50                       & 0.75                           & 1.13                             & 21.25                            & 2.9e-02                        & 2.8e-06                          & 4.8e-02                          & 1.25e-03                       \\
51                       & 0.71                           & 1.30                             & 21.63                            & 1.5e-02                        & 2.8e-06                          & 4.8e-02                          & 1.25e-03                       \\
52                       & 0.67                           & 1.17                             & 21.03                            & 1.4e-02                        & 8.4e-06                          & 4.5e-02                          & 1.25e-03                       \\
53                       & 0.61                           & 1.38                             & 23.29                            & 1.6e-02                        & 4.8e-06                          & 4.9e-02                          & 1.25e-03                       \\
54                       & 0.65                           & 1.43                             & 22.72                            & 1.5e-02                        & 8.1e-06                          & 4.8e-02                          & 1.25e-03                       \\
55                       & 0.69                           & 1.35                             & 22.32                            & 1.2e-02                        & 2.2e-06                          & 5.9e-02                          & 1.25e-03                       \\
56                       & 0.71                           & 1.14                             & 22.38                            & 9.0e-03                        & 1.1e-06                          & 5.0e-02                          & 1.25e-03                       \\
57                       & 0.74                           & 1.11                             & 23.86                            & 1.1e-02                        & 2.7e-06                          & 4.7e-02                          & 1.25e-03                       \\
58                       & 0.76                           & 1.26                             & 24.92                            & 1.2e-02                        & 6.4e-06                          & 4.5e-02                          & 1.25e-03                       \\
59                       & 0.77                           & 1.13                             & 24.43                            & 1.0e-02                        & 3.0e-06                          & 4.9e-02                          & 1.25e-03                       \\
60                       & 0.76                           & 1.27                             & 22.45                            & 6.2e-03                        & 2.3e-06                          & 5.2e-02                          & 1.25e-03                       \\
61                       & 0.73                           & 0.89                             & 22.88                            & 3.1e-03                        & 3.6e-06                          & 4.8e-02                          & 1.25e-03                       \\
62                       & 0.77                           & 0.99                             & 22.40                            & 1.8e-02                        & 3.1e-06                          & 5.1e-02                          & 1.25e-03                       \\
63                       & 0.65                           & 1.34                             & 23.83                            & 1.5e-02                        & 3.8e-06                          & 4.8e-02                          & 1.25e-03                       \\
64                       & 0.66                           & 1.24                             & 23.17                            & 1.5e-02                        & 4.3e-06                          & 4.4e-02                          & 1.25e-03                       \\
65                       & 0.66                           & 1.06                             & 23.71                            & 1.5e-02                        & 4.0e-06                          & 4.4e-02                          & 1.25e-03                       \\
66                       & 0.76                           & 1.12                             & 22.84                            & 8.9e-03                        & 3.9e-06                          & 4.6e-02                          & 1.25e-03                       \\
67                       & 0.76                           & 1.43                             & 22.87                            & 1.7e-02                        & 4.4e-06                          & 4.9e-02                          & 1.25e-03                       \\
68                       & 0.66                           & 1.49                             & 23.33                            & 4.2e-02                        & 3.8e-06                          & 5.4e-02                          & 1.25e-03                       \\
69                       & 0.70                           & 1.49                             & 24.14                            & 1.3e-02                        & 2.8e-06                          & 5.0e-02                          & 1.25e-03                       \\
70                       & 0.74                           & 1.33                             & 24.01                            & 1.8e-02                        & 9.7e-06                          & 5.1e-02                          & 1.25e-03                       \\
71                       & 0.70                           & 1.22                             & 24.35                            & 1.9e-02                        & 5.3e-06                          & 4.9e-02                          & 1.25e-03                       \\
72                       & 0.78                           & 1.25                             & 25.06                            & 1.4e-02                        & 8.2e-06                          & 4.9e-02                          & 1.25e-03                       \\
73                       & 0.81                           & 1.32                             & 24.64                            & 2.7e-02                        & 9.2e-06                          & 5.0e-02                          & 1.25e-03                       \\
74                       & 0.77                           & 1.38                             & 25.96                            & 3.6e-02                        & 1.2e-05                          & 4.6e-02                          & 1.25e-03                       \\
75                       & 1.17                           & 1.15                             & 26.05                            & 9.8e-03                        & 6.5e-06                          & 6.1e-02                          & 1.25e-03                       \\
76                       & 1.13                           & 1.47                             & 27.49                            & 1.8e-02                        & 9.5e-06                          & 4.6e-02                          & 1.25e-03                       \\
77                       & 0.81                           & 1.51                             & 26.83                            & 4.1e-02                        & 1.7e-05                          & 5.7e-02                          & 1.25e-03                       \\
78                       & 0.81                           & 1.32                             & 29.38                            & 4.5e-02                        & 1.4e-05                          & 5.5e-02                          & 1.25e-03                       \\
79                       & 1.02                           & 1.35                             & 37.53                            & 3.8e-02                        & 5.8e-06                          & 5.9e-02                          & 1.25e-03                       \\
80                       & 0.76                           & 1.70                             & 42.83                            & 3.4e-02                        & 1.7e-05                          & 5.4e-02                          & 1.25e-03                       \\
81                       & 0.80                           & 1.75                             & 39.81                            & 2.6e-02                        & 3.7e-06                          & 6.6e-02                          & 1.25e-03                       \\
82                       & 0.69                           & 1.67                             & 43.18                            & 3.2e-02                        & 5.6e-06                          & 6.0e-02                          & 1.25e-03                       \\
83                       & 0.97                           & 1.65                             & 46.16                            & 1.5e-02                        & 5.9e-06                          & 6.4e-02                          & 1.25e-03                       \\
84                       & 1.08                           & 1.73                             & 52.41                            & 1.3e-02                        & 1.2e-05                          & 6.7e-02                          & 1.25e-03                       \\
85                       & 0.95                           & 1.70                             & 83.30                            & 2.0e-02                        & 1.3e-05                          & 5.2e-02                          & 1.25e-03                       \\
86                       & 0.92                           & 1.83                             & 137.36                           & 1.1e-02                        & 9.9e-06                          & 4.9e-02                          & 1.25e-03                       \\
87                       & 0.84                           & 1.89                             & 137.15                           & 3.9e-02                        & 2.1e-05                          & 6.6e-02                          & 1.25e-03                       \\
88                       & 0.84                           & 2.30                             & 162.54                           & 2.5e-02                        & 1.8e-05                          & 5.4e-02                          & 1.25e-03                       \\
89                       & 0.80                           & 2.82                             & 167.03                           & 2.3e-02                        & 1.7e-05                          & 6.1e-02                          & 1.25e-03                       \\
90                       & 0.84                           & 2.56                             & 173.27                           & 3.1e-02                        & 2.0e-05                          & 5.2e-02                          & 1.25e-03                       \\
91                       & 0.83                           & 2.89                             & 170.28                           & 3.1e-02                        & 9.0e-06                          & 8.4e-02                          & 1.25e-03                       \\
92                       & 0.86                           & 2.84                             & 178.27                           & 4.5e-02                        & 3.9e-05                          & 7.4e-02                          & 1.25e-03                       \\
93                       & 0.90                           & 4.25                             & 180.67                           & 3.8e-02                        & 1.2e-05                          & 5.2e-02                          & 1.25e-03                       \\
94                       & 0.77                           & 3.55                             & 187.70                           & 4.6e-02                        & 1.1e-05                          & 5.3e-02                          & 1.25e-03                       \\
95                       & 0.91                           & 3.50                             & 186.12                           & 2.5e-02                        & 6.6e-06                          & 6.0e-02                          & 1.25e-03                       \\
96                       & 0.81                           & 4.22                             & 191.57                           & 2.9e-02                        & 3.3e-05                          & 4.9e-02                          & 1.25e-03                       \\
97                       & 0.94                           & 4.85                             & 198.83                           & 2.2e-02                        & 1.1e-05                          & 6.7e-02                          & 1.25e-03                       \\
98                       & 1.04                           & 4.28                             & 204.86                           & 1.9e-02                        & 1.8e-05                          & 5.6e-02                          & 1.25e-03                       \\
99                       & 1.02                           & 4.85                             & 207.31                           & 1.9e-02                        & 1.3e-05                          & 5.3e-02                          & 1.25e-03
\end{tabular}
\caption{\ref{prob_ddfact} bound for \ref{MESP}, varying $s$ ($n=2000$, $\rank(C)=150$) - Part I}\label{tab:ddfact_s1}
\end{table}

\begin{table}[!ht]
\centering
\footnotesize
\begin{tabular}{r|ccc|rrr|r}
\multicolumn{1}{c|}{}    & \multicolumn{3}{c|}{Elapsed time (sec)}                                                              & \multicolumn{3}{c|}{Dual gap}                                                                        & \multicolumn{1}{c}{$\rho$}     \\
\multicolumn{1}{c|}{}    & \multicolumn{1}{c}{}           & \multicolumn{1}{c}{}             & \multicolumn{1}{c|}{\tiny Frank} & \multicolumn{1}{c}{}           & \multicolumn{1}{c}{}             & \multicolumn{1}{c|}{\tiny Frank} & \multicolumn{1}{c}{}           \\[-4pt]
\multicolumn{1}{c|}{$r$} & \multicolumn{1}{c}{\tiny ADMM} & \multicolumn{1}{c}{\tiny KNITRO} & \multicolumn{1}{c|}{\tiny Wolfe} & \multicolumn{1}{c}{\tiny ADMM} & \multicolumn{1}{c}{\tiny KNITRO} & \multicolumn{1}{c|}{\tiny Wolfe} & \multicolumn{1}{c}{\tiny ADMM} \\ \hline
100                      & 1.12                           & 4.85                             & 218.01                           & 1.3e-02                        & 3.5e-05                          & 5.3e-02                          & 1.25e-03                       \\
101                      & 1.11                           & 4.52                             & 209.47                           & 1.1e-02                        & 1.5e-05                          & 5.8e-02                          & 1.25e-03                       \\
102                      & 1.20                           & 4.21                             & 223.35                           & 9.0e-03                        & 1.9e-05                          & 5.4e-02                          & 1.25e-03                       \\
103                      & 0.99                           & 4.40                             & 218.94                           & 1.8e-02                        & 8.9e-06                          & 5.3e-02                          & 1.25e-03                       \\
104                      & 0.93                           & 4.09                             & 226.80                           & 2.1e-02                        & 1.5e-05                          & 6.4e-02                          & 1.25e-03                       \\
105                      & 0.91                           & 4.07                             & 228.41                           & 4.0e-02                        & 1.5e-05                          & 5.0e-02                          & 1.25e-03                       \\
106                      & 1.01                           & 4.13                             & 227.71                           & 4.1e-02                        & 3.2e-05                          & 6.8e-02                          & 1.25e-03                       \\
107                      & 0.96                           & 4.29                             & 242.74                           & 3.6e-02                        & 1.9e-05                          & 6.1e-02                          & 1.25e-03                       \\
108                      & 0.98                           & 4.45                             & 234.64                           & 3.5e-02                        & 1.5e-05                          & 5.9e-02                          & 1.25e-03                       \\
109                      & 1.06                           & 4.07                             & 240.93                           & 3.5e-02                        & 2.7e-05                          & 5.5e-02                          & 1.25e-03                       \\
110                      & 1.03                           & 3.78                             & 251.14                           & 2.4e-02                        & 2.5e-05                          & 5.3e-02                          & 1.25e-03                       \\
111                      & 1.02                           & 4.46                             & 235.05                           & 1.8e-02                        & 2.0e-05                          & 5.9e-02                          & 1.25e-03                       \\
112                      & 1.05                           & 4.61                             & 262.75                           & 1.5e-02                        & 2.7e-05                          & 5.3e-02                          & 1.25e-03                       \\
113                      & 1.06                           & 3.58                             & 259.55                           & 1.3e-02                        & 4.2e-05                          & 7.0e-02                          & 1.25e-03                       \\
114                      & 0.99                           & 3.81                             & 270.63                           & 1.8e-02                        & 5.1e-05                          & 5.2e-02                          & 1.25e-03                       \\
115                      & 1.04                           & 3.50                             & 276.27                           & 1.6e-02                        & 4.0e-05                          & 5.3e-02                          & 1.25e-03                       \\
116                      & 1.00                           & 3.97                             & 292.88                           & 1.8e-02                        & 5.2e-05                          & 5.5e-02                          & 1.25e-03                       \\
117                      & 0.98                           & 3.85                             & 313.38                           & 2.1e-02                        & 3.3e-05                          & 5.2e-02                          & 1.25e-03                       \\
118                      & 1.02                           & 3.85                             & 321.75                           & 3.5e-02                        & 1.9e-05                          & 5.5e-02                          & 1.25e-03                       \\
119                      & 0.94                           & 3.91                             & 325.93                           & 3.6e-02                        & 4.0e-05                          & 5.9e-02                          & 1.25e-03                       \\
120                      & 1.08                           & 3.51                             & 327.84                           & 4.8e-02                        & 2.6e-05                          & 5.1e-02                          & 1.25e-03                       \\
121                      & 1.18                           & 3.91                             & 348.77                           & 3.4e-02                        & 2.1e-05                          & 5.3e-02                          & 1.25e-03                       \\
122                      & 1.18                           & 4.00                             & 359.50                           & 3.7e-02                        & 3.0e-05                          & 5.5e-02                          & 1.25e-03                       \\
123                      & 1.00                           & 3.97                             & 370.74                           & 4.1e-02                        & 2.4e-05                          & 5.0e-02                          & 1.25e-03                       \\
124                      & 1.17                           & 4.14                             & 379.31                           & 3.3e-02                        & 4.0e-05                          & 5.1e-02                          & 1.25e-03                       \\
125                      & 1.04                           & 3.98                             & 383.83                           & 4.2e-02                        & 6.5e-05                          & 4.9e-02                          & 1.25e-03                       \\
126                      & 1.07                           & 4.12                             & 370.85                           & 4.9e-02                        & 3.6e-05                          & 5.4e-02                          & 1.25e-03                       \\
127                      & 1.07                           & 4.03                             & 372.88                           & 4.8e-02                        & 1.4e-05                          & 5.2e-02                          & 1.25e-03                       \\
128                      & 1.05                           & 3.44                             & 400.67                           & 4.4e-02                        & 4.8e-05                          & 5.2e-02                          & 1.25e-03                       \\
129                      & 1.10                           & 4.29                             & 396.80                           & 3.6e-02                        & 2.0e-05                          & 5.9e-02                          & 1.25e-03                       \\
130                      & 1.20                           & 3.39                             & 412.03                           & 4.9e-02                        & 3.6e-05                          & 5.7e-02                          & 1.25e-03                       \\
131                      & 1.34                           & 3.53                             & 414.45                           & 4.4e-02                        & 1.8e-05                          & 5.4e-02                          & 1.25e-03                       \\
132                      & 1.40                           & 3.67                             & 400.70                           & 4.4e-02                        & 3.3e-05                          & 5.1e-02                          & 1.25e-03                       \\
133                      & 1.62                           & 3.80                             & 396.95                           & 4.5e-02                        & 1.8e-05                          & 5.1e-02                          & 1.25e-03                       \\
134                      & 1.56                           & 3.58                             & 420.78                           & 4.9e-02                        & 6.8e-05                          & 5.2e-02                          & 1.25e-03                       \\
135                      & 1.49                           & 3.72                             & 420.50                           & 4.7e-02                        & 3.4e-05                          & 5.0e-02                          & 1.25e-03                       \\
136                      & 1.51                           & 3.71                             & 423.31                           & 4.3e-02                        & 2.2e-05                          & 5.0e-02                          & 1.25e-03                       \\
137                      & 1.59                           & 4.09                             & 409.57                           & 4.1e-02                        & 3.4e-05                          & 5.2e-02                          & 1.25e-03                       \\
138                      & 1.61                           & 3.62                             & 430.51                           & 4.7e-02                        & 6.5e-05                          & 5.6e-02                          & 1.25e-03                       \\
139                      & 1.93                           & 3.07                             & 430.46                           & 4.5e-02                        & 3.0e-05                          & 5.9e-02                          & 1.25e-03                       \\
140                      & 1.71                           & 3.90                             & 427.55                           & 4.3e-02                        & 3.1e-05                          & 5.4e-02                          & 1.25e-03                       \\
141                      & 2.12                           & 3.10                             & 428.55                           & 4.5e-02                        & 6.2e-05                          & 5.1e-02                          & 5.25e-03                       \\
142                      & 1.91                           & 3.33                             & 416.84                           & 4.9e-02                        & 5.4e-05                          & 5.5e-02                          & 5.25e-03                       \\
143                      & 1.99                           & 3.35                             & 446.17                           & 4.1e-02                        & 4.4e-05                          & 5.0e-02                          & 5.25e-03                       \\
144                      & 2.19                           & 3.43                             & 460.96                           & 5.0e-02                        & 4.8e-05                          & 5.5e-02                          & 5.25e-03                       \\
145                      & 2.47                           & 3.66                             & 490.63                           & 4.8e-02                        & 4.0e-05                          & 5.1e-02                          & 5.25e-03                       \\
146                      & 1.99                           & 3.44                             & 456.84                           & 4.8e-02                        & 3.0e-05                          & 5.1e-02                          & 5.25e-03                       \\
147                      & 2.01                           & 3.36                             & 482.71                           & 4.6e-02                        & 3.2e-05                          & 5.1e-02                          & 5.25e-03                       \\
148                      & 1.94                           & 3.74                             & 498.64                           & 5.0e-02                        & 2.1e-05                          & 5.3e-02                          & 5.25e-03                       \\
149                      & 1.89                           & 4.21                             & 485.76                           & 5.0e-02                        & 3.1e-05                          & 5.3e-02                          & 5.25e-03                       \\
150                      & 2.50                           & 3.66                             & 520.17                           & 4.2e-02                        & 4.2e-05                          & 5.1e-02                          & 5.25e-03                      
\end{tabular}
\caption{\ref{prob_ddfact} bound for \ref{MESP}, varying $s$ ($n=2000$, $\rank(C)=150$) - Part II}\label{tab:ddfact_s2}
\end{table}

\begin{table}[!ht]
\centering
\footnotesize
\begin{tabular}{r|rrrr|rrrr|r}
                         & \multicolumn{4}{c|}{Elapsed time (sec)}                                                                                               & \multicolumn{4}{c|}{Dual gap}                                                                                                         & \multicolumn{1}{c}{$\rho$}     \\
\multicolumn{1}{c|}{$s$} & \multicolumn{1}{c}{\tiny ADMM} & \multicolumn{1}{c}{\tiny SDPT3} & \multicolumn{1}{c}{\tiny MOSEK} & \multicolumn{1}{c|}{\tiny COSMO} & \multicolumn{1}{c}{\tiny ADMM} & \multicolumn{1}{c}{\tiny SDPT3} & \multicolumn{1}{c}{\tiny MOSEK} & \multicolumn{1}{c|}{\tiny COSMO} & \multicolumn{1}{c}{\tiny ADMM} \\ \hline
43                       & 2.7                            & 4.3                             & 21.5                            & *                                & 8.5e-03                        & 5.2e-03                         & 3.3e-07                         & *                                & 1.25e-01                       \\
44                       & 2.6                            & 3.3                             & 21.0                            & *                                & 1.2e-02                        & 3.7e-03                         & 4.7e-07                         & *                                & 1.25e-01                       \\
45                       & 2.8                            & 3.7                             & 20.9                            & *                                & 1.4e-02                        & 5.6e-03                         & 7.3e-07                         & *                                & 1.20e-01                       \\
46                       & 3.3                            & 3.9                             & 21.5                            & *                                & 3.7e-03                        & 3.0e-03                         & 5.4e-07                         & *                                & 1.20e-01                       \\
47                       & 2.5                            & 3.6                             & 19.7                            & *                                & 3.3e-03                        & 4.9e-03                         & 1.9e-07                         & *                                & 1.20e-01                       \\
48                       & 3.5                            & 3.7                             & 20.4                            & *                                & 8.3e-03                        & 2.4e-03                         & 3.9e-07                         & *                                & 1.20e-01                       \\
49                       & 3.0                            & 3.5                             & 20.7                            & *                                & 1.9e-03                        & 1.1e-03                         & 5.0e-07                         & *                                & 1.20e-01                       \\
50                       & 2.6                            & 3.9                             & 15.0                            & *                                & 6.5e-03                        & 3.1e-03                         & 5.2e-07                         & *                                & 1.20e-01                       \\
51                       & 2.6                            & 4.5                             & 19.3                            & *                                & 7.9e-03                        & 8.1e-03                         & 8.1e-07                         & *                                & 1.20e-01                       \\
52                       & 3.5                            & 3.9                             & 18.9                            & *                                & 7.5e-04                        & 3.2e-03                         & 5.6e-07                         & *                                & 1.20e-01                      
\end{tabular}
\caption{\ref{bqp_original} bound for \ref{MESP}, varying $s$ ($n=63$)}\label{tab:bqpn63}
\end{table}

\begin{table}[!ht]
\centering
\footnotesize
\begin{tabular}{r|rrrr|rrrr|r}
                         & \multicolumn{4}{c|}{Elapsed time (sec)}                                                                                               & \multicolumn{4}{c|}{Dual gap}                                                                                                         & \multicolumn{1}{c}{$\rho$}     \\
\multicolumn{1}{c|}{$n$} & \multicolumn{1}{c}{\tiny ADMM} & \multicolumn{1}{c}{\tiny SDPT3} & \multicolumn{1}{c}{\tiny MOSEK} & \multicolumn{1}{c|}{\tiny COSMO} & \multicolumn{1}{c}{\tiny ADMM} & \multicolumn{1}{c}{\tiny SDPT3} & \multicolumn{1}{c}{\tiny MOSEK} & \multicolumn{1}{c|}{\tiny COSMO} & \multicolumn{1}{c}{\tiny ADMM} \\ \hline
250                      & 316.1                          & 831.9                           & *                               & *                                & 8.3e-03                        & 6.0e-03                         & *                               & *                                & 5.0e-02                        \\
275                      & 429.2                          & 1291.7                          & *                               & *                                & 1.4e-02                        & 2.9e-02                         & *                               & *                                & 5.0e-02                        \\
300                      & 583.2                          & 1916.5                          & *                               & *                                & 7.0e-03                        & 2.9e-02                         & *                               & *                                & 5.0e-02                        \\
325                      & 1008.6                         & *                               & *                               & *                                & 4.6e-03                        & *                               & *                               & *                                & 5.0e-02                        \\
350                      & 1878.6                         & *                               & *                               & *                                & 4.5e-03                        & *                               & *                               & *                                & 4.0e-02                        \\
375                      & 2866.9                         & *                               & *                               & *                                & 6.0e-03                        & *                               & *                               & *                                & 4.0e-02                        \\
400                      & 3279.3                         & *                               & *                               & *                                & 2.5e-02                        & *                               & *                               & *                                & 4.0e-02                       
\end{tabular}
\caption{\ref{bqp_original} bound for \ref{MESP}, varying $n$, with $s := \lfloor n/2\rfloor$}\label{tab:bqp_n}
\end{table}

\FloatBarrier


\section{Appendix: Computation of dual-feasible solutions}\label{sec:appdual}

In the following, we show how to construct dual-feasible solutions to \ref{prob}, \ref{prob_ddfact}, \ref{prob_linx}, and \ref{bqp_original}\,, from primal-feasible solutions. More details can be found in \cite*[Section 2]{PonteFampaLeeMPB} for \ref{prob}, in  \cite*[Section 3.4.4.1]{FL2022} for \ref{prob_ddfact}, in   \cite*[Section 3.3.4.1]{FL2022} for \ref{prob_linx}, and in  \cite*[Section 3.6.4]{FL2022} for  \ref{bqp_original}).


\subsection{The natural bound}

Using similar techniques as \cite*[Section 3.3.2]{FL2022}, we formulate the Lagrangian dual  of \ref{prob} as (see also \cite*[Section 2]{PonteFampaLeeMPB})
\begin{align}\label{dualnat}\tag{Du-$\mathcal{N}$}
\begin{array}{lll}
&\min &-\ldet \Psi   + \nu^\top \mathbf{e} + \delta s - {m},\\
&\text{s.t.} 
&\diag(A\Psi A^\top)  - \nu - \delta\mathbf{e} \leq 0,\\
&&\Psi \succ 0,\nu \geq 0.
\end{array}
\end{align}

Next, we show how to construct a closed-form feasible solution of \ref{dualnat} from a feasible solution $\hat x$ of \ref{prob} such that $A^\top\Diag(\hat x)A\!\in \!\mathbb{S}^m_{++}$\,,
 with the goal of having a small duality gap. 
 We define $\hat\Psi:= (A^\top \Diag(\hat x) A)^{-1}$. The minimum   gap between the objective value of  \ref{prob} at $\hat x$ and the objective value  of  \ref{dualnat} at feasible solutions $(\hat\Psi,\hat\nu,\hat\delta)$, is the optimal value of the linear program
\begin{align}\label{eq:g_theta}\tag{$G(\hat\Psi)$}
\begin{array}{rrl}
&\min &  \nu^\top \mathbf{e} + \delta s,\\
&\text{s.t.} 
&  \nu + \delta\mathbf{e} \geq \diag(A\hat\Psi A^\top)  ,\\
&&\nu \geq 0.
\end{array}
\end{align}
To obtain an optimal solution of  \ref{eq:g_theta},   we consider its  optimality conditions 
\begin{equation}\label{kktnatural}
\begin{array}{l}
    \diag(A\hat \Psi A^\top)   \leq \nu + \delta\mathbf{e},~\nu\geq 0,\\
    \mathbf{e}^\top x = s,~0\leq x \leq \mathbf{e},\\
     \nu^\top \mathbf{e} + \delta s = \diag(A\hat \Psi A^\top)^\top x.
\end{array}
\end{equation}

It is possible to verify that the following solution satisfies \eqref{kktnatural}.
 
\begin{align*}
&    \delta^*:=\diag(A\hat\Psi A^\top)_{\sigma(s)}~,\\
&\nu^*_{\sigma(\ell)}:=\left\{\begin{array}{ll}
\diag(A\hat\Psi A^\top)_{\sigma(\ell)}- \delta^*~,&\mbox{ for } 1\leq \ell\leq s;\\
0,&\mbox{ for } s< \ell\leq n,\\
\end{array}\right.\\
&x^*_{\sigma(\ell)}:=\left\{\begin{array}{ll}
1,&\mbox{ for } 1\leq \ell\leq s;\\
0,&\mbox{ for } s< \ell\leq n,\\
\end{array}\right.
\end{align*}
where $\sigma$ is the permutation  of the indices in $N$, such that $\diag(A\hat\Psi A^\top)_{\sigma(1)} \geq \dots \geq \diag(A\hat\Psi A^\top)_{\sigma(n)}$~.

Finally, $(\hat{\Psi},\nu^*,\delta^*)$ is the constructed dual-feasible solution to \ref{prob}.

From the optimality conditions for  \ref{prob}, we can see that if $\hat x$ is an optimal solution, then we have   $ \Psi=(A^\top \Diag(\hat x) A)^{-1}$ in an optimal solution of \ref{dualnat}. Therefore,  an optimal solution of  \ref{eq:g_theta} gives optimal values of the remaining variables  $(\nu,\delta)$ of \ref{dualnat}. In this case, due to strong duality for \ref{prob} and \ref{dualnat}, the optimal objective value of  \ref{mingapproblinx} is equal to zero.


\subsection{The DDFact bound}

The Lagrangian dual  of \ref{prob_ddfact} is (see \cite*[Section 3.4.2]{FL2022}, for a detailed derivation of the dual formulation)
\begin{align}\tag{DFact}\label{DFact}
\begin{array}{ll}
\min~ 
&\displaystyle - \sum_{\ell=k-s+1}^k \!\!\!\!\log\left(\lambda_{\ell} \left(\Psi\right)\right)
+ \nu^\top \mathbf{e}  +\delta s - s\\[4pt]
     \mbox{s.t.}
     & \diag(F \Psi F^\top)  - \nu   - \delta\mathbf{e}\leq 0,\\[4pt]
&\ \Psi\succ 0,  ~\nu\geq 0.
\end{array}
\end{align}

Next, we show how to construct a feasible solution of \ref{DFact} from a feasible solution $\hat x$ of \ref{prob_ddfact}, with the goal of producing a small gap. 

We  consider the spectral 
decomposition $F^\top \Diag(\hat x)F=\sum_{\ell=1}^{k} \hat \lambda_\ell \hat u_\ell \hat u_\ell^\top~,$
with $\hat \lambda_1\geq\hat \lambda_2\geq\cdots\geq \hat \lambda_{\hat r}>\hat \lambda_{\hat{r}+1}=\cdots=\hat \lambda_k=0$. 
We  define 
$\hat{\Psi}:=\sum_{\ell=1}^{k} {\hat \beta}_\ell \hat{u}_\ell \hat{u}_\ell^\top$~,
where 
\begin{equation}\label{betaepsilon}
\hat{\beta}_\ell:=\left\{
\begin{array}{ll}
        \textstyle 1/\hat{\lambda}_\ell~,      
       &\mbox{ for }1\leq \ell\leq \hat{\iota};\\ 
     1/\hat{\delta},&\mbox{ for }\hat{\iota}<\ell\leq \hat{r};\\ 
     (1+\epsilon)/\hat{\delta},&\mbox{ for }\hat{r}<\ell\leq k,
\end{array}\right.
\end{equation}
where $\epsilon>0$, and $\hat{\iota}$ is the unique integer defined  in Lemma \ref{Ni13} for $\lambda_\ell=\hat{\lambda}_\ell$~, and
$\hat \delta:=\frac{1}{s-\hat \iota}\sum_{\ell=\hat \iota+1}^{k}\hat \lambda_\ell
$~. As shown in \cite*[Section 2]{PonteFampaLeeMPB}, the smaller the value of $\epsilon$, the smaller the gap. In our computational experiments, taking rounding errors into account, we set $\epsilon = 0$.

The minimum duality gap between the objective value of  \ref{prob_ddfact} computed at $\hat x$ and the objective value of  \ref{DFact} computed at feasible solutions
of the form $(\hat\Psi,\nu,\delta)$,
is the optimal value of the linear program 

\begin{align}\label{mingapprobFact}\tag{$G(\hat\Psi)$}
\begin{array}{ll}
\min& 
 \nu^\top \mathbf{e}  +\delta s \\[3pt]
 \mbox{s.t.}&  \nu   + \delta\mathbf{e}\geq   \diag(F \hat \Psi F^\top) ,\\[3pt]
& \nu\geq 0.
\end{array}
\end{align}

\noindent  Following the same development of the previous subsection, we can verify that the following solution is optimal for \ref{mingapprobFact}.
\begin{align*}
&    \delta^*:=\diag(F \hat \Psi F^\top)_{\sigma(s)}~,\\
&\nu^*_{\sigma(\ell)}:=\left\{\begin{array}{ll}
\diag(F \hat \Psi F^\top)_{\sigma(\ell)}-\delta^* ~,&\mbox{ for } 1\leq \ell\leq s;\\
0,&\mbox{ for } s< \ell\leq n,\\
\end{array}\right.\\
&x^*_{\sigma(\ell)}:=\left\{\begin{array}{ll}
1,&\mbox{ for } 1\leq \ell\leq s;\\
0,&\mbox{ for } s< \ell\leq n\\
\end{array}\right.
\end{align*}
where $\sigma$ is the permutation  of the indices in $N$, such that $\diag(F \hat \Psi F^\top)_{\sigma(1)}\geq \cdots\geq \diag(F \hat \Psi F^\top)_{\sigma(n)}$~.

Finally, $(\hat{\Psi},\nu^*,\delta^*)$ is the constructed dual-feasible solution to \ref{prob_ddfact}.

We note that the choice of $\hat\Psi$ is motivated by the fact that, if $\hat{x}$ is an optimal solution of \ref{prob_ddfact}, then with this choice of $\hat{\Psi}$, the dual-feasible solution constructed is optimal (see Theorem 21 \cite*{GMESPalgorithmica} for a proof of this result for the more general problem GMESP).


\subsection{The linx bound}

The Lagrangian dual  of \ref{prob_linx} is (see \cite*[Section 3.3.2]{FL2022}, for a detailed derivation of the dual formulation) 
\begin{align}\label{Dlinx}\tag{Dlinx$_\gamma$}
\begin{array}{ll}
\min&-\frac{1}{2}\ldet (2\Psi) + \Trace(\Psi) + \nu^\top\mathbf{e}  + \delta s - n/2 \\[4pt]
     &\mbox{subject to:}\\[4pt]
&\diag(\gamma C\Psi C-\Psi)     -\nu  -  \delta\mathbf{e} \leq 0,\\[4pt]
&\Psi\succ 0, ~\nu\geq 0.
\end{array}
\end{align}

Next, we show how  to construct a feasible solution of \ref{Dlinx} from a feasible solution $\hat x$ of \ref{prob_linx} such that $L(\hat x)\succ 0$, with the goal of producing a small gap. 

We  define  $\hat{\Psi}:= \frac{1}{2} (L(\hat x))^{-1}$, and we see that $\frac{1}{2}\ldet L(\hat x) = -\frac{1}{2}\ldet (2\hat \Psi)$.
The minimum duality gap between the objective value of \ref{prob_linx} computed at $\hat x$ and the objective value of \ref{Dlinx} computed at feasible solutions of the form $(\hat\Psi,\nu,\delta)$,
is the optimal value of the linear program

\begin{align}\label{mingapproblinx}\tag{$G(\hat\Psi)$}
\begin{array}{cl}
\min~& 
 \nu^\top\mathbf{e}  + \delta s \\
     \mbox{s.t.}& 
         \nu  + \delta\mathbf{e} \geq \diag( \gamma C\hat\Psi C - \hat\Psi) ,\\
&\nu\geq 0.
\end{array}
\end{align}

\noindent  Following the same development of the previous subsection, we can verify that the following solution is optimal for \ref{mingapproblinx}.
 
\begin{align*}
&    \delta^*:=\diag(\gamma C\hat\Psi C - \hat\Psi)_{\sigma(s)}~,\\
&\nu^*_{\sigma(\ell)}:=\left\{\begin{array}{ll}
\diag(\gamma C\hat\Psi C - \hat\Psi)_{\sigma(\ell)}- \delta^*~,&\mbox{ for } 1\leq \ell\leq s;\\
0,&\mbox{ for } s< \ell\leq n,\\
\end{array}\right.\\
&x^*_{\sigma(\ell)}:=\left\{\begin{array}{ll}
1,&\mbox{ for } 1\leq \ell\leq s;\\
0,&\mbox{ for } s< \ell\leq n,\\
\end{array}\right.
\end{align*}
where $\sigma$ is the permutation  of the indices in $N$, such that $\diag(\gamma C\hat\Psi C - \hat\Psi)_{\sigma(1)}\geq \cdots\geq \diag(\gamma C\hat\Psi C - \hat\Psi)_{\sigma(n)}$~.

Finally, $(\hat{\Psi},\nu^*,\delta^*)$ is the constructed dual-feasible solution to \ref{prob_linx}\,.

From the optimality conditions for  \ref{prob_linx}\,, we can see that if $\hat x$ is an optimal solution, then we have   $ \Psi=\frac{1}{2}(L(\hat{x}))^{-1}$ in an optimal solution of \ref{Dlinx}\,. Therefore,  an optimal solution of  \ref{mingapproblinx} gives optimal values of the remaining variables  $(\nu,\delta)$ of \ref{Dlinx}\.. In this case, due to strong duality for \ref{prob_linx} and \ref{Dlinx}\,, the optimal objective value of  \ref{mingapproblinx} is equal to zero.


\subsection{The BQP bound}

The Lagrangian dual  of \ref{bqp_original} is (see \cite*[Section 3.6.2]{FL2022}, for a detailed derivation of the dual formulation) 
\begin{align}\tag{DBQP$_\gamma$}\label{DBQP}
\begin{array}{ll}
&\min -\ldet(\Psi) +\Trace(\Psi)  + \omega^\top g - (n+1) \\
     &\quad\mbox{subject to:}\\
&\quad\Psi\circ \tilde{C} 
     - \textstyle \sum_{\ell=1}^{2n+2} \omega_\ell G_\ell \preceq 0,\\
&\quad \Psi\succ 0,
\end{array}
\end{align}
where $\tilde{C}$, $G_\ell$ for $\ell=1,\ldots,2n+2$, and $g := (g_1,\ldots,g_{2n+2})^\top$ are defined in the same way as in problem \eqref{prob:bqp}.

Next, 
 we show how to construct a feasible solution of \ref{DBQP} from a feasible solution $(\hat x,\hat X)$ of \ref{bqp_original} such that $\gamma C\circ \hat{X} + \Diag(\mathbf{e}-\hat{x})\succ 0$, with the goal of producing a small gap. 
We  define  
\begin{equation*}
\hat{\tilde{X}}:=\left(\begin{array}{c l}
1&\hat x^\top \\    \hat x &\hat X\end{array}\right)
\end{equation*}
and
$\hat{\Psi}:= (\tilde{C}\circ \hat{\tilde{X}}+I_{n+1})^{-1}$.
The minimum duality gap between $(\hat x,\hat X)$ in \ref{bqp_original} and feasible solutions
of \ref{DBQP} of the form $(\hat\Psi,\omega)$
is the optimal value of the semidefinite program 

\begin{align}\label{mingapprobBQP}\tag{$G(\hat\Psi)$}
\begin{array}{ll}
\min & 
 \omega^\top g\\[5pt]
 \mbox{s.t.}&
 \textstyle  \sum_{\ell=1}^{2n+2} \omega_\ell G_\ell \succeq \hat \Psi\circ \tilde{C}.
\end{array}
\end{align}

Finally, $(\hat{\Psi},\omega^*)$ is the constructed dual-feasible solution to \ref{bqp_original}\,, where $\omega^*$ is an optimal solution to \ref{mingapprobBQP}.

From the optimality conditions for  \ref{bqp_original}, we can see that if $(\hat x,\hat X)$ is an optimal solution, then we have   $\Psi=(\tilde C\circ{\hat{\tilde X}} +I_{n+1})^{-1}$ in an optimal solution of \ref{DBQP}\,. Therefore,  an optimal solution of  \ref{mingapprobBQP} gives the optimal values of the remaining variable  $\omega$ of \ref{DBQP}\,. In this case, due to strong duality for \ref{bqp_original} and \ref{DBQP}\,, the optimal objective value of  \ref{mingapprobBQP} is equal to zero.

\end{document}